\documentclass[10pt, a4paper]{amsart}

\usepackage{hyperref}
\usepackage[all]{xy} 

\newcommand\AAA{\mathbb{A}}
\newcommand\KK{\mathbb{K}}
\newcommand\TT{\mathcal{T}}
\newcommand\ZZ{\mathbb{Z}}
\newcommand\QQ{\mathbb{Q}}
\newcommand\PP{\mathbb{P}}
\newcommand\LL{\mathcal{L}}
\newcommand\xx{\mathbf{x}}
\newcommand\yy{\mathbf{y}}
\newcommand\Pone{{\PP^1}}
\newcommand\Ptwo{{\PP^2}}
\newcommand\Pthree{{\PP^3}}
\newcommand\Pfour{{\PP^4}}
\newcommand\Pfive{{\PP^5}}
\newcommand\Psix{{\PP^6}}
\newcommand{\tS}{{\smash{\widetilde S}}}
\newcommand{\rto}{{\dashrightarrow}}
\newcommand{\ADE}{{\mathbf{ADE}}}
\newcommand{\Aone}{{\mathbf A}_1}
\newcommand{\Atwo}{{\mathbf A}_2}
\newcommand{\Athree}{{\mathbf A}_3}
\newcommand{\Afour}{{\mathbf A}_4}
\newcommand{\Afive}{{\mathbf A}_5}
\newcommand{\Asix}{{\mathbf A}_6}
\newcommand{\Aseven}{{\mathbf A}_7}
\newcommand{\Aeight}{{\mathbf A}_8}
\newcommand{\Dfour}{{\mathbf D}_4}
\newcommand{\Dfive}{{\mathbf D}_5}
\newcommand{\Dsix}{{\mathbf D}_6}
\newcommand{\Dseven}{{\mathbf D}_7}
\newcommand{\Deight}{{\mathbf D}_8}
\newcommand{\Esix}{{\mathbf E}_6}
\newcommand{\Eseven}{{\mathbf E}_7}
\newcommand{\Eeight}{{\mathbf E}_8}
\newcommand\dpbox[2]{#2}
\newcommand{\toric}{\text{toric}}
\newcommand{\onerel}{\text{$1$ relation}}
\newcommand{\tworel}{\text{$\ge 2$ relations}}
\newcommand{\ex}[1]{*+<5pt>[o][F]{E_{#1}}}
\newcommand{\li}[1]{*+<3pt>[F]{E_{#1}}}
\DeclareMathOperator\rk{rk}
\DeclareMathOperator\Pic{Pic}
\DeclareMathOperator\Cox{Cox}
\DeclareMathOperator\Bl{Bl}
\DeclareMathOperator{\Spec}{Spec}
\DeclareMathOperator{\Hom}{Hom}
\newcommand{\Tns}{{T_{\mathrm{NS}}}}
\newcommand{\Gm}{\mathbb{G}_{\mathrm{m}}}
\newcommand\cl[1]{\bar{#1}}
\newcommand{\coxcase}[2]{\medskip\noindent\textbf{Type #2. }}

\newcommand{\coxembed}[2]{Anticanonical model $\pi: \tS \to S \subset #1$:
  \[\begin{split}#2.\end{split}\]
}

\newcommand{\coxbirat}[2]{Birational map (inducing $\rho : \tS \to \Ptwo$): 
  $\phi: S \rto \Ptwo,\ \xx \mapsto (#2)$.\\
}

\newcommand\coxbiratproj[3]{Projection from $#3$ (inducing $\rho : \tS \to
  \Ptwo$): $\phi: S \rto \Ptwo,\ \xx \mapsto (#2)$.\\
}

\newcommand\coxbiratinv[2]{Birational map (inducing $\rho : \tS \to \Ptwo$):
  $\phi: S \rto \Ptwo,\ \xx \mapsto (#1)$, inverse
  \begin{equation*}
    \yy \mapsto (#2).
  \end{equation*}
}

\newcommand\coxclass[2]{\cl{E}_{#1}&=#2}
\newcommand\coxcurve[2]{\pi(E_{#1})=\{#2\}}
\newcommand\coxproj[2]{\rho(E_{#1})=\{#2=0\}}
\newcommand\coxprojp[2]{\rho(E_{#1})=\{#2\}} 

\newcommand\coxcc[1]{#1,\displaybreak[0]\\}
\newcommand\coxcd[1]{#1.}

\newcommand\coxc[4]{\coxclass{#1}{#2},\ \coxproj{#1}{#4},\ \coxcurve{#1}{#3}}
\newcommand\coxco[4]{\coxclass{#1}{#2},\ \coxproj{#1}{#4}}
\newcommand\coxcs[4]{\coxclass{#1}{#2},\ \coxproj{#1}{#4},\\\coxcurves{#1}{#3}}
\newcommand\coxcss[4]{\coxclass{#1}{#2},\\\coxprojs{#1}{#4},\ \coxcurve{#1}{#3}}
\newcommand\coxcp[4]{\coxclass{#1}{#2},\ \coxprojp{#1}{#4},\ \coxcurve{#1}{#3}}
\newcommand\coxcop[4]{\coxclass{#1}{#2},\ \coxprojp{#1}{#4}}

\newcommand\coxim[4]{\coxcc{\coxc{#1}{#2}{#3}{#4}}} 
\newcommand\coximd[4]{\coxcd{\coxc{#1}{#2}{#3}{#4}}} 
\newcommand\coxims[4]{\coxcc{\coxcs{#1}{#2}{#3}{#4}}} 
\newcommand\coximss[4]{\coxcc{\coxcss{#1}{#2}{#3}{#4}}} 
\newcommand\coximds[4]{\coxcd{\coxcs{#1}{#2}{#3}{#4}}} 
\newcommand\coximdss[4]{\coxcd{\coxcss{#1}{#2}{#3}{#4}}} 
\newcommand\coximp[4]{\coxcc{\coxcp{#1}{#2}{#3}{#4}}} 
\newcommand\coximpd[4]{\coxcd{\coxcp{#1}{#2}{#3}{#4}}} 
\newcommand\coximo[4]{\coxcc{\coxco{#1}{#2}{#3}{#4}}} 
\newcommand\coximod[4]{\coxcd{\coxco{#1}{#2}{#3}{#4}}} 
\newcommand\coximop[4]{\coxcc{\coxcop{#1}{#2}{#3}{#4}}} 
\newcommand\coximopd[4]{\coxcd{\coxcop{#1}{#2}{#3}{#4}}} 

\newcommand\dynkin[1]{\xymatrix@R=0.05in @C=0.05in{#1}}

\newcommand{\coxdynkin}[1]{Extended Dynkin diagram:
  \begin{equation*}
    \dynkin{#1}
  \end{equation*}
}

\newcommand{\coxdynkinp}[3]{Extended Dynkin diagram 
  ($\pi(#1)=#2$):
  \begin{equation*}
    \dynkin{#3}
  \end{equation*}
}

\newcommand{\coxdynkinpp}[2]{Extended Dynkin diagram (#1):
  \begin{equation*}
    \dynkin{#2}
  \end{equation*}
}

\newcommand{\coxonerel}[3]{Cox ring: generators
  $#1$ with relation (degree $#3$)
  \[#2.\]}

\newcommand{\coxanti}[1]{Anticanonical sections $\pi^*(x_i)$:
  \begin{align*}
    #1.\end{align*}}

\newcommand{\coxantie}[2]{Anticanonical sections $\pi^*(x_i)$:
  \begin{gather*}
    #1,\\
    #2.
  \end{gather*}}

\newcommand{\coxtwoone}[4]{Singularity: $#1$ in $#2$ with $#3$:
\begin{align*}#4\end{align*}}

\newcommand{\coxtwotwo}[7]{Singularities: $#1$ in $#2$ with $#3$;
  $#4$ in $#5$ with $#6$:
\begin{align*}#7\end{align*}}

\newcommand{\coxtwomulti}[4]{Singularities: $#1$ in $#2$  with #3,
  respectively:
\begin{align*}#4\end{align*}}

\newcommand{\coxonem}[2]{$(-1)$-curves: $#1$ with \begin{align*}#2\end{align*}}

\newcommand{\coxone}[2]{$(-1)$-curve: $#1$ with \begin{align*}#2\end{align*}}

\newcommand{\coxother}[2]{Extra generators: #1
  with \begin{align*}#2\end{align*}}

\newcommand\OO{\mathcal{O}}
\newcommand\e{\eta}
\newcommand\inj{\hookrightarrow}

\newtheorem{theorem}{Theorem}
\newtheorem{lemma}[theorem]{Lemma}
\newtheorem{prop}[theorem]{Proposition} 

\theoremstyle{definition}
\newtheorem{definition}[theorem]{Definition} 
\newtheorem{remark}[theorem]{Remark}
\newtheorem{example}[theorem]{Example}
\newtheorem{step}{Step}
\newtheorem{result}{Result}
\newtheorem*{ack}{Acknowledgements}

\numberwithin{equation}{section}

\begin{document}

\title[Del Pezzo surfaces whose universal torsors are hypersurfaces]
  {Singular del Pezzo surfaces\\ whose universal torsors are hypersurfaces}

\author{Ulrich Derenthal}

\address{Mathematisches Institut, Ludwig-Maximilians-Universit\"at M\"unchen, 
  Theresienstr. 39, 80333 M\"un\-chen, Germany}

\email{ulrich.derenthal@mathematik.uni-muenchen.de}

\date{July 22, 2013}


\begin{abstract}
  We classify all generalized del Pezzo surfaces (i.e., minimal
  desingularizations of singular del Pezzo surfaces containing only
  rational double points) whose universal torsors are open subsets of
  hypersurfaces in affine space.  Equivalently, their Cox rings are
  polynomial rings with exactly one relation.  For all $30$ types with
  this property, we describe the Cox rings in detail.

  These explicit descriptions can be applied to study Manin's
  conjecture on the asymptotic behavior of the number of rational
  points of bounded height for singular del Pezzo surfaces, using the
  universal torsor method.
\end{abstract}

\maketitle

\tableofcontents

\section{Introduction}

Universal torsors were introduced by Colliot-Th\'el\`ene and Sansuc to study
the Hasse principle and weak approximation for rational points on del Pezzo
surfaces and other geometrically rational varieties \cite{MR54:2657},
\cite{MR0447246}, \cite{MR0447250}, \cite{MR605344}, \cite{MR89f:11082}.

Universal torsors are also applied in the context of Manin's conjecture
\cite{MR89m:11060}, \cite{MR1032922} on the asymptotic behavior of the
number of rational points of bounded height on Fano varieties. This
goes back to Salberger \cite{MR1679841} and Peyre
\cite{MR1679842}. After the proof of Manin's conjecture for toric
varieties (which include several del Pezzo surfaces) by Batyrev and
Tschinkel \cite{MR1620682}, Salberger gave a different proof using
universal torsors \cite{MR1679841}.

Subsequently, Salberger's approach to Manin's conjecture via universal
torsors was applied to many non-toric varieties, in particular del
Pezzo surfaces. In each case, one needs a very precise understanding
of universal torsors, in particular of defining equations. In this
article, we focus on universal torsors that have a presentation as
Zariski open subsets of affine hypersurfaces.

The first examples, universal torsors over minimal desingularizations of
cubic surfaces with singularities of types $\Esix$ resp.\ $\Dfour$, were
worked out by Hassett and Tschinkel \cite{MR2029868}.  This lead to the first
complete proof of Manin's conjecture for a non-toric singular cubic surface
over $\QQ$ in \cite{MR2332351} and was a starting point in the program to
study Manin's conjecture for all types of split singular del Pezzo surfaces over
number fields and function fields.

For an overview of further progress towards Manin's conjecture for del
Pezzo surfaces over $\QQ$, we refer to \cite{MR2362193},
\cite{MR2290499}, \cite{MR2559866}. If Manin's conjecture is known for
a del Pezzo surface of a certain \emph{type} (see
Definition~\ref{def:type} and Remark~\ref{rem:type_classification}), we
give a reference to the proof in
Tables~\ref{tab:degree_789}--\ref{tab:degree_1}.

Del Pezzo surfaces of some types are toric or equivariant compactifications of
the additive group $\mathbb{G}_{\mathrm{a}}^2$ (classified in
\cite{MR2753646}), so Manin's conjecture follows from \cite{MR1620682} or
\cite{MR1906155}, without using universal torsors.  We observe that most other
proofs of Manin's conjecture are concerned with one particular example of a
split singular del Pezzo surface such that a universal torsor of its minimal
desingularization is an open subset of an affine hypersurface. So far, all
proofs for examples of this kind are based on the description of universal
torsors given in \cite{MR2029868} or in (preliminary versions of) our work
that we present here.

Bourqui \cite{MR2809202} proved Manin's conjecture for toric varieties
over function fields using techniques analogous to
\cite{MR1620682}. Additionally, he developed universal torsor methods
over function fields to prove Manin's conjecture for some non-toric
del Pezzo surfaces \cite{MR2573192}, \cite{arXiv:1205.3573}. His work
also uses our explicit description of universal torsors of del Pezzo
surfaces.

\medskip

From here, we work over an algebraically closed field of
characteristic $0$. We use the following terminology. A normal
projective surface with ample anticanonical class is called
\emph{ordinary del Pezzo surface} if it is non-singular, or
\emph{singular del Pezzo surface} if its singularities are rational
double points. A \emph{generalized del Pezzo surface} is a
non-singular projective surface whose anticanonical class is big and
nef. Generalized del Pezzo surfaces are precisely ordinary del Pezzo
surfaces and minimal desingularizations of singular del Pezzo
surfaces.

Universal torsors are closely related to Cox rings: Cox rings of generalized
del Pezzo surfaces are finitely generated, and a Zariski open subset of
$\Spec(\Cox(\tS))$ is a universal torsor over the generalized del Pezzo
surface $\tS$ \cite[Theorems~5.6, 5.7]{MR2498061}. Therefore, determining
generators of $\Cox(\tS)$ and all relations between them gives coordinates and
defining equations for the affine variety $\Spec(\Cox(\tS))$ containing a
universal torsor.

While there is no canonical choice of generators of $\Cox(\tS)$, we will see
(Lemma~\ref{lem:generators}) that $\Cox(\tS)$ can be described by a certain
minimal number of generators and relations that are homogeneous with respect
to the natural $\Pic(\tS)$-grading of $\Cox(\tS)$ and whose degrees with
respect to this grading are uniquely defined. Whenever we mention (numbers and
degrees of) generators and relations of Cox rings, we refer to such minimal
sets of homogeneous generators and relations.

As universal torsors of a generalized del Pezzo surface $\tS$ of
degree $d$ have dimension $\dim(\tS)+\rk(\Pic(\tS))=12-d$, we can
distinguish three cases:
\begin{enumerate}
\item\label{it:toric} Precisely in the toric case \cite{MR95i:14046},
  \cite[Corollary~2.10]{MR2001i:14059}, Cox rings are polynomial rings in
  $12-d$ variables (without relations), and their universal torsors can be
  presented as open subsets of affine space $\AAA^{12-d}$.
\item\label{it:one_rel} For Cox rings that are isomorphic to the quotients of
  polynomial rings in $13-d$ variables by one relation, the corresponding
  universal torsors are open subsets of affine hypersurfaces in $\AAA^{13-d}$.
\item\label{it:more_rel} If there are at least $14-d$ generators and more than
  one relation in the Cox ring, universal torsors are open subsets of
  (possibly non-complete) intersections of at least two hypersurfaces in
  affine space of dimension $\ge 14-d$.
\end{enumerate}

For applications such as Manin's conjecture, it is interesting to distinguish
between these three cases. In case~\eqref{it:toric}, Manin's conjecture is
known by the work of Batyrev and Tschinkel \cite{MR1620682}. It turns out that
it tends to be easier to prove Manin's conjecture for del Pezzo surfaces in
case~\eqref{it:one_rel} than in case~\eqref{it:more_rel}.

For an ordinary del Pezzo surface $\tS$, the Cox ring depends essentially on
the degree $d$. For $d \ge 6$, $\tS$ is toric, so we are in
case~\eqref{it:toric}. For $d \le 5$, the Cox ring has been determined in
\cite{MR2029863}, \cite{MR2358614}, \cite{MR2529093}, \cite{MR2579393}, and it
turns out that we are in case~\eqref{it:more_rel}. Here, $\Cox(\tS)$ has one
generator corresponding to each negative curve on $\tS$, plus two generators
of anticanonical degree in case $d=1$ \cite[Theorem~3.2]{MR2029863}.

For generalized del Pezzo surfaces, the shape of the Cox ring depends
not only on the degree. The classification of generalized del Pezzo
surfaces is closely related to the classification of singular del
Pezzo surfaces, which was done in degree $3$ by Schl\"afli
\cite{schlaefli} and Cayley \cite{cayley} and in degrees $1$ and $2$
by Du Val \cite{duval}. This leads to a finite number of types (see
Definition~\ref{def:type}), where each type in degree $d \le 7$
corresponds to a subsystem of the root system $R_d$ in
Table~\ref{tab:root_system}. Each type may be described by the degree,
the types of the $\ADE$-singularities and the number of the lines on
the corresponding singular del Pezzo surfaces (where the number of
lines may be ignored in most cases; see
Remark~\ref{rem:type_classification}).

\begin{table}[ht]
  \centering
  \begin{equation*}
    \begin{array}{cccccccc}
      \hline
      d & 7 & 6 & 5 & 4 & 3 & 2 & 1 \\
      \hline
      R_d & \Aone & \Atwo+\Aone & \Afour & \Dfive & \Esix & \Eseven & \Eeight\\
      \hline
    \end{array}
  \end{equation*}
  \smallskip
  \caption{Root systems $R_d$ of del Pezzo surfaces of degree $d$}
  \label{tab:root_system}
\end{table}

As the example of the $\Esix$ cubic surface \cite{MR2029868} shows,
$\Cox(\tS)$ may have generators whose degrees are neither classes of negative
curves nor the anticanonical class. However, the following result allows us to
determine the degrees of generators of $\Cox(\tS)$ recursively, starting with
del Pezzo surfaces of high degree $\ge 7$, which are always toric, so that
their Cox ring is known by \cite{MR95i:14046}. We denote the class of a
divisor $D$ by $\cl D$; in particular, $-\cl K_\tS$ is the class of an
arbitrary anticanonical divisor $-K_\tS$ on $\tS$.

\begin{theorem}\label{thm:generators1}
  Let $\tS$ be a generalized del Pezzo surface. Then the number and
  degrees of generators and relations of $\Cox(\tS)$ depend only on
  the type of $\tS$.

  If $\deg(\tS) \in \{2, \dots, 7\}$, then the $\Pic(\tS)$-degree of
  each generator of $\Cox(\tS)$ is the class of a negative curve or
  $-\cl K_\tS$ or $\rho^*\cl D$ for some contraction $\rho:\tS \to
  \tS'$ of a negative curve on $\tS$ and the degree $\cl D$ of a
  generator of $\Cox(\tS')$.
\end{theorem}

In Section~\ref{sec:strategy}, we present facts about del Pezzo surfaces and
their universal torsors and Cox rings. We prove Theorem~\ref{thm:generators1}
(in the slightly more precise version stated as Theorem~\ref{thm:generators})
and several auxiliary results. We show how to determine the shape of the Cox
rings for all types of generalized del Pezzo surfaces in every degree $d$,
given only their abstract classification via root systems.

We compute the degrees of generators and relations of the Cox ring for
each type in a purely ``combinatorial'' way, without using projective
models, in Steps~\ref{step:types} and
\ref{step:generator_degrees}. Whenever this leads to
case~\eqref{it:one_rel} of Cox rings with one relation, we show how to
determine an explicit description of the corresponding singular
del Pezzo surfaces $S$ and of $\Cox(\tS)$. It
seems interesting to note that we first construct $\Cox(\tS)$ by
choosing a relation in Step~\ref{step:relation}, and then extract an
explicit description of $S$ from this in Step~\ref{step:model},
together with extra information on the structure on $S$, its
desingularization $\tS$ and generators of $\Cox(\tS)$.

We construct $S$ and $\tS$ in a way that they are defined and
\emph{split} over $\QQ$ (i.e., $\tS$ is obtained from $\Ptwo$ by a
series of blow-ups, each defined over $\QQ$, and $S$ and its lines and
singularities are all defined over $\QQ$). Also $\Cox(\tS)$ is defined
over $\QQ$. Therefore, our results cover the case of split del Pezzo
surfaces over $\QQ$. The classification, geometry and arithmetic of
\emph{non-split} del Pezzo surfaces over $\QQ$ is more involved, see
for example \cite{MR89f:11083}.

Some (but not all) of our del Pezzo surfaces can be equipped with a
complexity one torus action, so that their Cox rings can also be
determined using the methods of Hausen and S\"u\ss{} \cite{MR2671185}
once a projective model is known.  For some singular del Pezzo
surfaces, Cox rings have been determined by this method; see
\cite[Theorems~3.23--3.26]{arXiv:1106.0854} and the references therein.

The fact that our construction starts with the abstract classification of del
Pezzo surfaces instead of an explicit description of each type (for example by
defining equations) seems both natural and useful in practice (lists of
defining equations do not seem to be readily available in degrees other than
$3$ and $4$).

This leads to the following classification (see also Table~\ref{tab:overview}):

\begin{theorem}\label{thm:classification}
  The Cox rings of generalized del Pezzo surfaces $\tS$ of the following types
  have a minimal set of $13-\deg(\tS)$ generators with one relation:
  \begin{itemize}
  \item degree $6$: types $\Aone$ (with three lines), $\Atwo$;
  \item degree $5$: types $\Aone$, $\Atwo$, $\Athree$, $\Afour$;
  \item degree $4$: types $3\Aone$, $\Atwo+\Aone$, $\Athree$ (with five lines),
    $\Athree+\Aone$, $\Afour$, $\Dfour$, $\Dfive$;
  \item degree $3$: types $\Dfour$, $\Athree+2\Aone$, $2\Atwo+\Aone$, $\Afour+\Aone$,
    $\Dfive$, $\Afive+\Aone$, $\Esix$;
  \item degree~$2$: types $\Dfive+\Aone$, $\Esix$, $2\Athree+\Aone$,
    $\Afive+\Atwo$, $\Dfour+3\Aone$, $\Dsix+\Aone$,~$\Eseven$;
  \item degree $1$: types $\Esix+\Atwo$, $\Eseven+\Aone$, $\Eeight$.
  \end{itemize}
  All other types are either toric or have Cox rings with at least two
  relations.
\end{theorem}

\begin{table}[ht]
  \centering
  \begin{tabular}{cccc}
    \hline
    degree & toric & $1$ relation & $\ge 2$ relations\\
    \hline
    9 & 1 type & -- & --\\
    8 & 3 types & -- & --\\
    7 & 2 types & -- & --\\
    6 & 4 types & 2 types & --\\
    5 & 2 types & 4 types & 1 type\\
    4 & 3 types & 7 types & 6 types\\
    3 & 1 type & 7 types & 13 types\\
    2 & -- & 7 types & 39 types\\
    1 & -- & 3 types & 71 types\\
    \hline
    $\sum$ & 16 types & 30 types & 130 types\\
    \hline
  \end{tabular}
  \smallskip
  \caption{Relations in Cox rings of generalized del Pezzo surfaces}
  \label{tab:overview}
\end{table}

In Section~\ref{sec:classification}, we give a detailed presentation
of the results. For Cox rings with one relation, our choice of data
listed there is lead on the one hand by the question which information
is crucial to determine the Cox rings and to describe them precisely,
and on the other hand by the question which information is needed for
applications such as Manin's conjecture (see \cite{MR2290499} and
\cite{MR2520770} for a systematic description of some aspects).

\begin{ack}
  This work was supported by grant DE 1646/2-1 of the Deutsche
  Forschungsgemeinschaft, by grant 200021\_124737/1 of the Schweizer
  Nationalfonds and by the Center for Advanced Studies of LMU
  M\"unchen. I thank V.~Batyrev, J.~Hausen, B.~Hassett, J.~Heinloth,
  A.~Laface and the referee for useful remarks.
\end{ack}

\section{Construction}\label{sec:strategy}

For the minimal desingularization $\tS$ of every singular del Pezzo surface
$S$, we determine whether its Cox ring $\Cox(\tS)$ has precisely one
relation. If this is the case, we give an explicit and natural description of
$S$ in (weighted) projective space with its lines and singularities, of $\tS$
as a blow-up of $\Ptwo$ with its negative curves, and of the generators and
relation in $\Cox(\tS)$.

In this section, we describe how this can be achieved starting with the
abstract classification of generalized del Pezzo surfaces via root
systems. When we determine the number of generators and relations of
$\Cox(\tS)$, it turns out that the information mentioned above can be
extracted along the way. In particular, we construct an anticanonical
embedding of $S$ via the Cox ring.

\subsection{Classification of generalized del Pezzo surfaces}\label{sec:del_pezzo}

Let $\KK$ be an algebraically closed field of characteristic $0$.  In
this section, we recall the classification of generalized del Pezzo
surfaces via root systems. For basic properties of generalized del
Pezzo surfaces, see also \cite[\S III]{MR579026}, \cite{MR89f:11083},
\cite{MR2227002}, \cite[\S 8]{MR2964027} (where they are called
\emph{rational surfaces of negative type}, \emph{generalized del Pezzo
  surfaces}, \emph{Gorenstein log del Pezzo surfaces}, \emph{weak del
  Pezzo surfaces}, respectively).

Let $\tS$ be a generalized del Pezzo surface. Its Picard group $\Pic(\tS)$ is
equipped with the intersection form $(\cdot , \cdot)$, which is a
non-degenerate bilinear form. For $n \in \ZZ$, a curve $C$ on $\tS$ is called
an \emph{$(n)$-curve} if $C$ is isomorphic to $\Pone$ and its class $\cl{C} \in
\Pic(S)$ has self intersection number $(\cl{C},\cl C)=n$. A \emph{negative curve}
is an $(n)$-curve with $n<0$. The degree $\deg(\tS)$ of a generalized del
Pezzo surface is the self intersection number of its anticanonical class
$-\cl K_\tS$. The Picard group $\Pic(\tS)$ of a generalized del Pezzo surface is
free of rank $10-\deg(\tS)$.

Any generalized del Pezzo surface $\tS$ is isomorphic to
$\Ptwo$ (of degree $9$), $\Pone \times \Pone$, the Hirzebruch surface $F_2$
(both of degree $8$) or a blow-up of $\Ptwo$ in $r \le 8$ points in
\emph{almost general position} (of degree $9-r$), i.e., we have a map $\rho:
\tS \to \Ptwo$ that is a composition
\begin{equation}\label{eq:blow_up_sequence}
  \tS = \tS_r \xrightarrow{\rho_r} \tS_{r-1} \to \dots \to \tS_1 \xrightarrow{\rho_1} \tS_0 = \Ptwo,
\end{equation}
where $\rho_i : \tS_i \to \tS_{i-1}$ is the blow-up of a closed point $p_i \in
\tS_{i-1}$ not lying on a $(-2)$-curve on $\tS_{i-1}$ (see
\cite[Proposition~8.1.16]{MR2964027}).

A generalized del Pezzo surface $\tS$ is ordinary if and only if it does not
contain $(-2)$-curves, i.e., if it is $\Ptwo$ or $\Pone \times \Pone$ or the
blow-up of $\Ptwo$ in $r$ points in \emph{general position}. In the notation
above, general position means that $p_1, \dots, p_r$ do not lie on
$(-1)$-curves, and additionally for $r=8$, the point $p_8$ is not a double
point on a $(2)$-curve. Equivalently, $p_1, \dots, p_8$ are distinct points in
$\Ptwo$, with no three on one line, no six on one conic and no eight on a
cubic curve in $\Ptwo$ with one of them a singular point of that curve.

If $\tS$ is a blow-up of $\Ptwo$ in $r$ points as
in~(\ref{eq:blow_up_sequence}), let $\ell_0 = \rho^*\OO_\Ptwo(1)$ and
$\ell_i$ be the class of the total transform on $\tS$ of the
exceptional divisor of $\rho_i$, for $i=1, \dots, r$. Then $\ell_0,
\ell_1, \dots, \ell_r$ form a basis of $\Pic(\tS)$. The intersection
form is defined on this basis by
\begin{equation*}
  (\ell_i,\ell_j) =
  \begin{cases}
    1, &i=j=0,\\
    -1, &i=j>0,\\
    0, &i \ne j.
  \end{cases}
\end{equation*}
In particular, excluding the cases $\Pone \times \Pone$ and $F_2$ in
degree $8$, the Picard group $\Pic(\tS)$ together with its intersection form
depends only on the degree of $\tS$.

For the class of an irreducible curve $C$ on $\tS$, three cases occur:
\begin{itemize}
\item $\cl C=\ell_i$ with $i \in \{1, \dots, r\}$ if and only if $C$ is the
  $(-1)$-curve that is the strict transform of the exceptional divisor of
  $\rho_i$ that is not blown up subsequently.
\item $\cl C=\ell_i-\ell_j$ with $i < j \in \{1, \dots, r\}$ if and only if
  $C$ is the $(-2)$-curve that is the strict transform of the exceptional
  divisor of $\rho_i$, with $p_j$ on its strict transform on $\tS_{j-1}$.
\item $\cl C=a_0\ell_0-a_1\ell_1-\dots-a_r\ell_r$ with $a_0 > 0$ and $a_1,
  \dots, a_r \ge 0$ if and only if $C$ is the strict transform of the curve
  $\rho(C)$ of degree $a_0$ in $\Ptwo$, with $p_i$ a point of multiplicity
  $a_i$ on its strict transform on $\tS_{i-1}$, for $i=1, \dots, r$. Then $C$
  has self intersection number $a_0^2-a_1^2-\dots-a_r^2$.
\end{itemize}

\begin{definition}\label{def:type}
  Two generalized del Pezzo surfaces $\tS$, $\tS'$ have the same
  \emph{type} if there is an isomorphism $\Pic(\tS) \cong \Pic(\tS')$
  preserving the intersection form that gives a bijection between
  their sets of classes of negative curves.
\end{definition}

For $\tS$ of degree $d \le 7$, the group of isomorphisms of $\Pic(\tS)$
preserving the intersection form is the Weyl group $W(R_d)$ of the root system
\begin{equation*}
  R_d = \{\cl D \in \Pic(\tS) \mid (\cl D,\cl D)=-2,\ (\cl D,-\cl K_\tS)=0\},
\end{equation*}
as in Table~\ref{tab:root_system}.

The types of generalized del Pezzo surfaces of degree $d$ are in
bijection to the subsystems of the root system $R_d$ (up to
automorphisms of $R_d$; labeled according to the
$\ADE$-classification), except the subsystems of types $7\Aone$ of
$R_2$ and of types $7\Aone$, $8\Aone$ and $\Dfour+4\Aone$ of $R_1$
(which occur only over fields of characteristic $2$); see
\cite{MR713283}, \cite{MR2227002}.  Here, the classes of the
$(-2)$-curves are the simple roots of the subsystem $R$. Hence the
Dynkin diagram of $R$ describes the configuration of $(-2)$-curves on
the corresponding generalized del Pezzo surface. The classes of the
$(-1)$-curves are precisely the elements $\cl D \in \Pic(\tS)$ with
$(\cl D,\cl D)=-1$, $(\cl D,-\cl K_\tS)=1$ and $(\cl D,\cl E) \ge 0$
for all classes $\cl E$ of $(-2)$-curves on~$\tS$.

\begin{remark}\label{rem:type_classification}
  We note that $R_d$ may contain root systems $R,R'$ that are abstractly
  isomorphic, without the existence of an automorphism of $R_d$ mapping $R$ to
  $R'$. In this case, $R, R'$ correspond to two different types. For example,
  there are two types of generalized del Pezzo surfaces of degree $4$ with an
  $\Athree$-configuration of $(-2)$-curves: one with four and one with five
  $(-1)$-curves.

  It turns out that the type of a generalized del Pezzo surface is uniquely
  determined by its degree, its configuration of $(-2)$-curves and its number
  of $(-1)$-curves. Indeed, this can be checked from the data in
  \cite{MR89f:11083} for degree $d \ge 4$ and from \cite{MR80f:14021} for
  $d=3$. The article \cite{MR713283} lists all cases with $d \in \{1,2\}$
  where the configuration of $(-2)$-curves together with $d$ does not
  determine the type uniquely. In each case, there are precisely two types,
  and it is straightforward to check that they differ by their number of
  $(-1)$-curves.

  For any type, we mention the number of $(-1)$-curves in addition to the
  degree and the $\ADE$-type of the root system only when this is
  necessary to identify the type uniquely.
\end{remark}

In particular, there is only a finite number of types. Our first step
is to find them together with the configuration of classes of negative
curves, encoded in the \emph{extended Dynkin diagram of negative curves}. This gives the
classification in Tables~\ref{tab:degree_789}--\ref{tab:degree_1},
which recovers the lists in \cite{MR89f:11083} for degree $d \ge 4$,
\cite{MR80f:14021} for $d=3$ and \cite{MR713283} for $d \le 2$. For
simplicity, for $d \le 2$, we only list the types with at most $13-d$
negative curves; all others have Cox rings with more than one relation
by Theorem~\ref{thm:generators} below.

\begin{step}[Find all types of generalized del Pezzo surfaces]\label{step:types}
  We apply the procedure described in \cite[Section~0.3]{MR2227002} to find
  all types of generalized del Pezzo surfaces. For each type, this gives an
  \emph{extended Dynkin diagram of negative curves} (where a vertex
  marked by a circle corresponds to the class of a $(-2)$-curve, a vertex
  marked by a square corresponds to the class of a $(-1)$-curve, with $(\cl
  E,\cl E')$ edges between the vertices of $\cl E$ and $\cl E'$). It is empty
  for $\Ptwo$ and $\Pone \times \Pone$. Otherwise, the extended Dynkin diagram of negative curves
  characterizes each type uniquely.

  For each type of degree $d \ge 7$, we determine a sequence of contractions
  $\rho_r, \dots, \rho_1$ of $(-1)$-curves corresponding to
  (\ref{eq:blow_up_sequence}) as follows: for $i=r, \dots, 1$, the effect of
  $\rho_i$ is to remove a class $\cl E$ with self intersection number $-1$
  from the extended Dynkin diagram of negative curves, to increase the self intersection number of
  the other classes $\cl E'$ by $(\cl E, \cl E')^2$ and to add $(\cl E, \cl
  E')\cdot (\cl E, \cl E'')$ edges between $\cl E', \cl E''$.

  This allows us to express the classes of negative curves in terms of the
  basis $\ell_0, \dots, \ell_r$ of $\Pic(\tS)$:
  \begin{itemize}
  \item $\cl E = \ell_i$ if $(\cl E, \cl E)=-1$ and (the strict transform of)
    $\cl E$ is contracted by $\rho_i$.
  \item $\cl E = \ell_i-\ell_j$ with $i<j$ if $(\cl E, \cl E)=-2$ such that
    $\rho_j$ increases the self intersection number to $-1$ and $\rho_i$
    contracts its strict transform.
  \item $\cl E = a_0\ell_0-a_1\ell_1-\dots-a_r\ell_r$ if $\cl E$ is not
    contracted by $\rho$, where $a_0>0$ is the self intersection number of its
    image in $\Ptwo$ and $a_i \ge 0$ is the intersection number of
    its strict transform on $\tS_i$ with $\ell_i$.
  \end{itemize}
\end{step}

\begin{result}\label{res:types}
  This step gives a list of all types. We identify each type by the
  degree, the $\ADE$-type of the configuration of $(-2)$-curves and
  the number of $(-1)$-curves. Furthermore, we give the extended
  Dynkin diagram of negative curves and list the classes $\cl E_i$ of
  the negative curves in terms of the basis $\ell_0, \dots, \ell_r$
  the Picard group.
\end{result}

\subsection{Degrees of generators of the Cox rings}\label{sec:generatoring_degrees}

We summarize basic facts about universal torsors and Cox rings of
generalized del Pezzo surfaces; see also \cite{MR89f:11082},
\cite{MR2029868}, \cite{arXiv:1003.4229}. Then we show how to find the
degrees in $\Pic(\tS)$ of generators of $\Cox(\tS)$ for $\tS$ of any
given type, based on the description of the classes of negative curves
on $\tS$ obtained in Step~\ref{step:types}. Several of the arguments
and intermediate results in this section have appeared (sometimes
implicitly) in previous work on Cox rings of rational surfaces (see
\cite{MR2029863}, \cite{MR2029868}, \cite{MR2115006},
\cite{MR2824848}, for example).

Let $\tS$ be a generalized del Pezzo surface of degree $d$. Let $D_0, \dots,
D_r$ be divisors on $\tS$ such that their classes form a basis of $\Pic(\tS)$,
where $r=9-d$. Then the \emph{Cox ring} of $\tS$ with respect to $D_0, \dots,
D_r$ is defined as
\begin{equation*}
  \Cox(\tS) = \Cox(\tS,D_0, \dots,D_r) = \bigoplus_{(\nu_0, \dots, \nu_r) \in
    \ZZ^{r+1}} H^0(\tS, \OO(\nu_0D_0+\dots+\nu_rD_r)),
\end{equation*}
where the multiplication is induced by the multiplication of global
sections. A different choice of divisors gives rise to a non-canonically
isomorphic Cox ring. Our classification results will be independent of this
choice, so that we use the notation $\Cox(\tS)$.

Since $\tS$ is rational and $-\cl K_\tS$ is big and nef, $\Cox(\tS)$ is a
finitely generated $\KK$-algebra by \cite[Theorem~5.7]{MR2498061}. It is
naturally graded by $\Pic(\tS)$. For $\cl D \in \Pic(\tS)$, let
$\Cox(\tS)_{\cl D}$ be its degree-$\cl D$-part, which is isomorphic to
$H^0(\tS,\OO(D))$. Whenever we talk about homogeneous elements $s \in
\Cox(\tS)$, this is meant with respect to this grading, and we write $\deg(s)
\in \Pic(\tS)$ for the degree. See \cite{arXiv:1003.4229} for more
information.

For $i=0, \dots, r$, let $\LL_i$ be the invertible sheaf $\OO(D_i)$, and let
$\LL_i^\circ$ be the $\Gm$-torsor over $\tS$ obtained by removing the
zero-section from $\LL_i$. Then
\begin{equation*}
  \TT = \LL_0^\circ \times_\tS \dots \times_\tS \LL_r^\circ
\end{equation*}
with the natural map to $\tS$ is a \emph{universal torsor} over $\tS$ in the
sense of \cite[\S 2]{MR89f:11082}. It is a $\Tns$-torsor over $\tS$, where $\Tns =
\Hom(\Pic(\tS), \Gm) \cong \Gm^{r+1}$ is the N\'eron-Severi torus with
character group $\Pic(\tS)$. Then $\TT$ is an open subset of the affine
variety $\Spec(\Cox(\tS))$ \cite[Theorem~5.6]{MR2498061}.

We say that $\e_1, \dots, \e_t$ are \emph{(a minimal system of) generators} of
$\Cox(\tS)$ if they are homogeneous, generate $\Cox(\tS)$ as a $\KK$-algebra,
and if no proper subset of $\e_1, \dots, \e_t$ is enough to generate
$\Cox(\tS)$.  Then the $\KK$-algebra homomorphism
\begin{equation}\label{eq:cox_polynomial}
  \Psi: \KK[X_1, \dots, X_t] \to \Cox(\tS)
\end{equation}
defined by $\Psi(X_i) = \e_i$ for $i=1, \dots, t$ is surjective.  By
Hilbert's basis theorem, the kernel of $\Psi$ is a finitely generated
ideal. If $R_1, \dots, R_m \in \KK[X_1, \dots, X_t]$ are homogeneous
(with respect to the $\Pic(\tS)$-grading induced from $\Cox(\tS)$ via
$\Psi$) and generate $\ker(\Psi)$, but no proper subset generates this
kernel, we say that $R_1(\e_1, \dots, \e_t)=\dots=R_m(\e_1, \dots,
\e_t)=0$ are \emph{(a minimal system of) relations} in $\Cox(\tS)$.

\begin{lemma}\label{lem:generators}
  Let $\tS$ be a generalized del Pezzo surface. For any minimal system
  of generators and any corresponding minimal system of relations of
  $\Cox(\tS)$, the numbers of their elements and the degrees of their
  elements depend only on $\tS$.
\end{lemma}

\begin{proof}
  The effective cone of $\tS$ is strictly convex. Therefore, we have a partial
  ordering on $\Pic(\tS)$, where $\cl D' \le \cl D$ if and only if $\cl D-\cl
  D'$ is the class of an effective divisor. As $\Cox(\tS)_{\cl D}$ is
  non-trivial only for classes of effective divisors $D$, this part is only
  affected by generators and relations of the finitely many effective degrees
  $\cl D' \le \cl D$.

  Therefore, for a minimal set of homogeneous generators of
  $\Cox(\tS)$, the number of generators of degree $\cl D$ is
  $h^0(\tS,\OO(D))-\dim H_{\cl D}$ where $H_{\cl D}$ is the
  intersection of $\Cox(\tS)_{\cl D}$ with the subring of $\Cox(\tS)$
  generated by all homogeneous $s \in \Cox(\tS)$ with $\deg(s)<\cl D$,
  which is independent of any choices.

  Similarly, the number of relations of degree $\cl D$ is $\dim
  \ker(\Psi)_{\cl D} - \dim J_{\cl D}$ where $\ker(\Psi)_{\cl D}$ is the
  degree-$\cl D$-part of $\ker(\Psi)$ and $J_{\cl D}$ is its intersection with
  the ideal generated by all homogeneous elements $R$ of $\ker(\Psi)$ with
  $\deg(R) < \cl D$.
\end{proof}

\begin{remark}\label{rem:toric}
  Toric generalized del Pezzo surfaces are easily classified: $\Ptwo$, $\Pone
  \times \Pone$ and $F_2$ are toric. If $\tS$ is a toric blow-up of $\Ptwo$ in
  $r$ points, it must be obtained by a sequence of blow-ups of points
  invariant under the torus action. These are the intersection points of the
  two torus-invariant curves corresponding to adjacent rays in the fan of the
  toric varieties. Such a point may be blown up if the self intersection
  numbers of these two curves are both $\ge -1$.

  It is straightforward to produce a list of all possibilities of such
  sequences of blow-ups. By \cite{MR95i:14046}, the generators of $\Cox(\tS)$
  correspond to the rays in the fan of $\tS$. In particular, this includes
  information on the configuration of negatives curves on $\tS$, so that one
  can determine all types (as in Definition~\ref{def:type}) of toric
  generalized del Pezzo surfaces, shown in Figure~\ref{fig:blow-ups} together
  with all blow-ups between them.

\begin{figure}[ht]
  \begin{equation*}
    \xymatrix    @R=0.2in @C=0.2in{
      d=9 & 8 & 7 & 6 & 5 & 4 &3\\
      &\dpbox{7}{\Pone \times \Pone} & & \dpbox{6}{\Bl_3\Ptwo} \ar@{->}[dl] &
       & \dpbox{4}{4\Aone} \ar@{->}[dl]\\
      \dpbox{9}{\Ptwo}& \dpbox{8}{\Bl_1\Ptwo} \ar@{->}[l] &
      \dpbox{7}{\Bl_2\Ptwo} \ar@{->}[ul] \ar@{->}[l] & 
      \dpbox{6}{\Aone \text{($4$ l.)}} \ar@{->}[l] \ar@{->}[dl] & \dpbox{5}{2\Aone}
      \ar@{->}[l] \ar@{->}[dl] \ar@{->}[ul] & \dpbox{4}{\Atwo+2\Aone}
      \ar@{->}[l] \ar@{->}[dl] & \dpbox{3}{3\Atwo} \ar@{->}[l]\\
      &\dpbox{8}{F_2} & \dpbox{7}{\Aone} \ar@{->}[l] \ar@{->}[ul]&
      \dpbox{6}{2\Aone} \ar@{->}[ul] \ar@{->}[l] & \dpbox{5}{\Atwo+\Aone}
      \ar@{->}[ul] \ar@{->}[l] \ar@{->}[dl] & \dpbox{4}{\Athree+2\Aone} \ar@{->}[l] \\
      & & & \dpbox{6}{\Atwo+\Aone} \ar@{->}[ul] 
    }
  \end{equation*}
  \caption{Toric generalized del Pezzo surfaces.}
  \label{fig:blow-ups}
\end{figure}
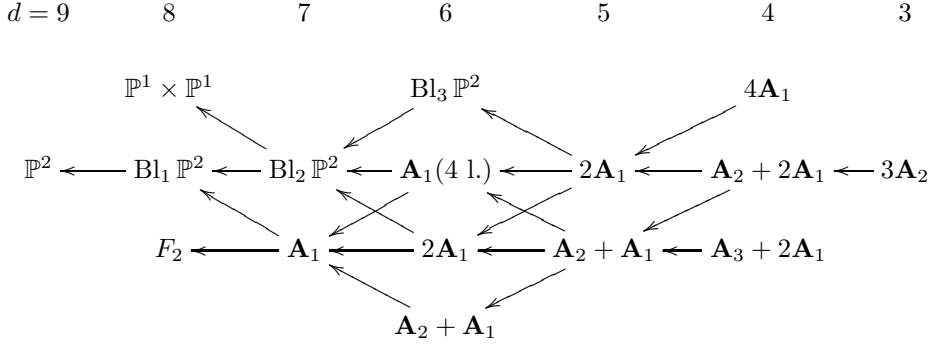\end{remark}

Our next goal is to determine the degrees of generators and relations of
$\Cox(\tS)$. Our strategy is similar to \cite[Theorem~3.2]{MR2029863}. The
following technical lemma is used in the proof of
Theorem~\ref{thm:generators}.

\begin{lemma}\label{lem:base_point_free}
  Let $\tS$ be a generalized del Pezzo surface of degree $d \ge 2$ that is a
  blow-up of $\Ptwo$. Let $D$ be a nef divisor on $\tS$. Then the linear
  system $|D|$ has no base points.
\end{lemma}

\begin{proof}
  For $d=9$, we have $\tS \cong \Ptwo$, with $D$ trivial or very ample, so the
  statement is clear.
  
  For $d < 9$, if there is a $(-1)$-curve $E$ on $\tS$ with $(\cl D,\cl E)=0$,
  we consider the map $\rho:\tS \to \tS'$ that contracts $E$. Then $\tS'$ is a
  generalized del Pezzo surface of degree $d+1$. We have $D = \rho^*(D')$ for
  some divisor $D'$ on $\tS'$. As the negative curves on $\tS'$ are images of
  negative curves of $\tS$ under $\rho$, the divisor $D'$ is nef.  By
  induction on $d$, the linear system $|D'|$ has no base points.
  Consequently, $|D|=|\rho^*D'|$ has no base points.
  
  If there is no $(-1)$-curve $E \subset \tS$ with $(D,E) = 0$, then let $m
  > 0$ be the minimum of $(\cl D,\cl E)$ over all $(-1)$-curves $E$.  Since $(\cl E,
  -\cl K_\tS) = 1$ for all $(-1)$-curves and $(\cl E, -\cl K_\tS) = 0$ for all
  $(-2)$-curves, $D'=D-m\cdot(-K_\tS)$ is nef, and $(\cl D',\cl E) = 0$ for some
  $(-1)$-curve $E$. As discussed before, $|D'|$ has no base points. Since $d
  \ge 2$, the system $|-K_\tS|$ is also base point free. Therefore, $|D|$
  has no base points.
\end{proof}

All generalized del Pezzo surfaces of degree $\ge 7$ are toric, so the degrees
of the generators of their Cox rings are easily determined (see
Remark~\ref{rem:toric}). For a generalized del Pezzo surface of lower degree,
the following result allows us to determine recursively a finite number of
potential degrees of generators of the Cox ring.

\begin{theorem}\label{thm:generators}
  Let $\tS$ be a generalized del Pezzo surface of degree $d \in \{2, \dots,
  7\}$. For each negative curve $E$ on $\tS$, there is precisely one generator
  of $\Cox(\tS)$ of degree $\cl E$. Furthermore, there are at most two
  generators of $\Cox(\tS)$ of degree $-\cl K_\tS$. The degree $\cl D$ of any
  other generator of $\Cox(\tS)$ is nef, and there is a contraction $\rho :
  \tS \to \tS'$ of a $(-1)$-curve on $\tS$ such that $\cl D = \rho^*(\cl D')$
  for some nef $\cl D'$ that is the degree of a generator of $\Cox(\tS')$.
\end{theorem}

\begin{proof}
  For $d \le 7$, the effective cone of $\tS$ is minimally generated by its
  negative curves (see \cite[Theorem~3.10]{MR2377367}, for
  example). Therefore, $h^0(\tS, \OO(E))=1$ for each negative curve $E$, and
  any minimal system of homogeneous generators of $\Cox(\tS)$ contains
  precisely one section corresponding to $E$.

  The degrees of all other generators of $\Cox(\tS)$ are nef. Indeed,
  this is a special case of \cite[Theorem~1]{MR2115006}, but we give a
  short direct argument here. If an effective divisor $D$
  corresponding to a generator $s$ of $\Cox(\tS)$ is not nef, then
  $(\cl D,\cl E)<0$ for some negative curve $E$, which implies that
  $E$ is a component of $D$. On the other hand, $D$ is a prime divisor
  because otherwise $s$ would be the product of sections corresponding
  to its components. Therefore, $D$ must coincide with $E$. See also
  \cite[Proposition~1]{MR2824848}.

  In the following, a \emph{distinguished section} is an element of
  $\Cox(\tS)$ that can be expressed as a polynomial in some sections
  whose degrees are as in the statement of the result.  For each nef
  $\cl D \in \Pic(\tS)$, we must show that each element of
  $H^0(\tS,\OO(D))$ is a distinguished section.  We use induction over
  the partial order on the effective cone in $\Pic(\tS)$. If $\cl D$
  is trivial, the claim is clear.

  If $(\cl D,\cl E) = 0$ for some $(-1)$-curve $E$, consider the contraction
  $\rho: \tS \to \tS'$ of $E$. Then $D = \rho^*(D')$ for some nef divisor $D'$
  on $\tS'$, and
  \begin{equation*}
    \rho^* : H^0(\tS', \OO(D')) \to H^0(\tS, \OO(D))
  \end{equation*}
  is an isomorphism. If $\cl D'$ is not the degree of a generator of
  $\Cox(\tS')$, then each $s' \in H^0(\tS', \OO(D'))$ is a polynomial in
  generators of $\Cox(\tS')$ of smaller degree, so the claim is true by
  induction over the degree after applying $\rho^*$. Otherwise, $\cl D'$ is
  the degree of a nef generator of $\Cox(\tS')$, so any section of degree $\cl
  D=\rho^*(\cl D')$ is distinguished by definition.

  If $(\cl D,\cl E) \ge 1$ for all $(-1)$-curves $E$ and $(\cl D,\cl E) \ge 0$
  for all $(-2)$-curves $E$, let $m$ be the minimum of $(\cl D,\cl E)$ over
  all $(-1)$-curves $E$. Let $F$ be a $(-1)$-curve with $(\cl D,\cl F) =
  m$, and let $s$ be a section corresponding to $F$.

  Consider the nef divisor $D'=D-m\cdot(-K_\tS)$. Since $|D'|$ is
  base point free by Lemma~\ref{lem:base_point_free}, we can choose a section
  $s' \in H^0(\tS, \OO(D'))$ whose support does not contain $F$. Then $s'$ is a
  distinguished section by induction.

  We consider the commutative diagram
  \begin{equation*}
    \xymatrix@R=0.3in @C=0.3in{
      && H^0(\tS,\OO(-mK_\tS)) \ar@{->}[r]\ar@{->}[d]^{\cdot s'} & H^0(F, \OO(-mK_\tS)|_{F})\ar@{->}[d]_\sim^f\\
      0\ar@{->}[r] & H^0(\tS, \OO(D-F))\ar@{->}[r]^{\cdot s} & H^0(\tS, \OO(D))\ar@{->}[r] & H^0(F,\OO(D)|_{F}).
    }
  \end{equation*}

  By induction, any section in $H^0(\tS, \OO(D-F))$ is
  distinguished. Therefore, it is enough show that $H^0(\tS, \OO(D))$ contains
  distinguished sections that restrict to generators of $H^0(F, \OO(D)|_{F})$,
  with $\OO(D)|_{F} \cong \OO_\Pone(m)$ because of $(\cl D,\cl F)=m$.
  
  We have distinguished sections $s_1, s_2 \in H^0(\tS, \OO(-K_\tS))$ that
  restrict to a basis of $H^0(F,\OO(-K_\tS)|_{F})$. Indeed, $|-K_\tS|$ is base
  point free by Lemma~\ref{lem:base_point_free} (using $d \ge 2$). This shows
  that $\Cox(\tS)$ has at most two generators of degree $-\cl K_\tS$.

  For any other $\cl D$, the distinguished sections $s's_1^is_2^{m-i}$ for
  $i=0, \dots, m$ restrict to a basis of $H^0(F,\OO(D)|_{F})$ since
  $s_1^is_2^{m-i}$ restrict to a basis of $H^0(F, \OO(-mK_\tS)|_{F})$ and
  since multiplication by $s'$ restricts to the isomorphism $f$ by $(\cl D', \cl
  F)=0$.
\end{proof}

\begin{lemma}\label{lem:dimension_nef}
  Let $\tS$ be a generalized del Pezzo surface whose Cox ring has $t$
  generators of degrees $\cl D_1, \dots, \cl D_t$ and one relation of degree
  $\cl D_0$.  For any $\cl D \in \Pic(\tS)$, let $d(\cl D)$ be the number of
  ways to express $\cl D$ as a non-negative integral linear combination of
  $\cl D_1, \dots, \cl D_t$. Then $\Cox(\tS)_{\cl D}$ has dimension
  $d(\cl D)-d(\cl D-\cl D_0)$. If $\cl D$ is nef, this is equal to the Euler
  characteristic $\chi(\OO(D)) = \frac 1 2 ((\cl D,\cl D)+(\cl D,-\cl K_\tS))+1$.
\end{lemma}

\begin{proof}
  Since there is only one relation, we have the exact sequence
  \begin{equation*}
    0 \to \KK[X_1, \dots, X_t] \to \KK[X_1, \dots, X_t] \to \Cox(\tS) \to 0,
  \end{equation*}
  where the second map is $\Psi$ as in (\ref{eq:cox_polynomial}) and
  the first map is multiplication by the relation $R$ generating
  $\ker(\Psi)$. If we consider the first $\KK[X_1, \dots, X_t]$ with
  the $\Pic(\tS)$-grading induced from $\Cox(\tS)$ shifted by the degree
  $\cl D_0$ of the relation, this is an exact sequence of graded
  vector spaces, where the degree-$\cl D$-part has dimension $d(\cl D)$
  in the middle and dimension $d(\cl D-\cl D_0)$ on the left hand side
  because of the shift. The first claim follows.

  The second claim holds by the vanishing theorem of Kawamata--Viehweg as in
  \cite[Corollary~1.10]{MR2029868} and the Riemann--Roch formula.
\end{proof}

The following result allows us to show that the Cox ring has more than one
relation even in some cases with $\le 13-d$ negative curves.

\begin{prop}\label{prop:more_generators}
  Let $\tS$ be a generalized del Pezzo surface. Let $\rho: \tS \to
  \tS'$ be the contraction of a $(-1)$-curve on $\tS$. If $\Cox(\tS)$
  has $t$ generators, then $\Cox(\tS')$ has at most $t-1$ generators.
\end{prop}

\begin{proof}
  Let $r = 9-d$. Via the natural embedding $\rho^*: \Pic(\tS') \inj
  \Pic(\tS)$, we may assume that $\Pic(\tS')$ has the basis $\ell_0,
  \dots, \ell_{r-1}$, and $\Pic(\tS)$ has the basis $\ell_0, \dots,
  \ell_r$, where $\ell_r$ is the class of the exceptional divisor of
  $\rho$. Let us assume that $\Cox(\tS)$ is minimally generated by
  $\e_1, \dots, \e_t$. Since this includes sections corresponding to
  all negative curves by Theorem~\ref{thm:generators}, we may assume
  that $\deg(\e_t)=\ell_r$. Furthermore, in a minimal set of
  generators, the curve corresponding to a generator of $\Cox(\tS)$
  must be irreducible. Therefore, in all cases except the exceptional
  divisor of $\rho$, it is the strict transform of an irreducible
  curve corresponding to a section $\e_i'$ on $\tS'$, so $\deg(\e_i) =
  \deg(\e_i')- k_i\ell_r$ and $\rho^*(\e_i') = \e_i\e_t^{k_i}$, where
  $k_i$ is a non-negative integer.

  We claim that $\e_1', \dots, \e_{t-1}'$ generate
  $\Cox(\tS')$. Indeed, for each $\cl D' \in \Pic(\tS')$, the map $\rho^*:
  H^0(\tS', \OO(D')) \to H^0(\tS, \rho^*\OO(D'))$ is an isomorphism. Therefore,
  any $\rho^*(s)$ with $\deg(s)=\cl D'$ can be expressed as a linear
  combination of monomials in $\e_1, \dots, \e_t$. We note that any
  $\e_1^{a_1}\cdots \e_{t-1}^{a_{t-1}}\e_t^{a_t}$ of degree $\rho^*\cl D'$
  is the pullback of $(\e_1')^{a_1}\cdots(\e_{t-1}')^{a_{t-1}}$ of
  degree $\cl D'$ via $\rho^*$ because of the relation between
  $\deg(\e_i), \deg(\e_i')$ mentioned above (where $a_t =
  \sum_{i=1}^{t-1} k_ia_i$ because of the degree). Therefore, $s$ can
  be expressed using only $\e_1', \dots, \e_{t-1}'$.
\end{proof}

\begin{step}[Determine the degrees of generators and of the relation
  of the Cox ring]\label{step:generator_degrees}
  For any type of degree $d$, each of the $s$ classes $\cl E_1, \dots,
  \cl E_s$ of negative curves is the degree of a generator $\e_i$ of
  $\Cox(\tS)$. If $s \ge 14-d$, then $\Cox(\tS)$ has more than one
  relation. If $s \le 13-d$, we check whether the contraction of some
  $(-1)$-curve leads to a generalized del Pezzo surface $\tS'$ of
  degree $d+1$ whose Cox ring has more than one relation and hence
  more than $13-(d+1)$ generators. If this is the case, the Cox ring of
  $\tS$ has more than $13-d$ generators by
  Proposition~\ref{prop:more_generators} and hence more than one
  relation.

  Otherwise, we determine finitely many nef degrees $\cl D_i$ that may
  be degrees of generators of $\Cox(\tS)$ recursively using
  Theorem~\ref{thm:generators}. We assume that their ordering $\cl
  D_1, \dots, \cl D_t$ is compatible with the partial ordering of
  $\Pic(\tS)$ (i.e., if $\cl D_i < \cl D_j$, then $i < j$).  For $i=1,
  \dots, t$, we apply Lemma~\ref{lem:dimension_nef} to $\cl D_i$ as
  follows: with the knowledge of the number of generators in each
  degree $\cl D_1, \dots, \cl D_{i-1}$ and whether a relation of some
  degree $\cl D_0<\cl D_i$ is known, we compute $d(\cl D_i)$ or $d(\cl
  D_i)-d(\cl D_i-\cl D_0)$ under the assumption that there is no
  generator and no relation of degree $\cl D_i$. If this number is
  smaller than $\chi(\cl D_i)$, there must be generators of degree
  $\cl D_i$. If it is larger, there must be a relation; in this case,
  we look for the smallest of the finitely many degrees $\cl D_0 \le
  \cl D_i$ where this relation occurs.

  If we end up finding $13-d$ generators, but no relation yet, we look
  for it in multiples of $-\cl K_\tS$ (where it must show up because
  $-\cl K_\tS$ is in the interior of the effective cone), and then in all
  smaller degrees.
\end{step}

\begin{result}
  We list the degrees $\cl E_i$ of the nef generators of the Cox ring
  and the degree of the relation in terms of $\ell_0, \dots, \ell_r$
  and compute their self intersection numbers. We add these degrees to
  the extended Dynkin diagram of negative curves to obtain the
  \emph{extended Dynkin diagram of Cox ring generators} (with a vertex
  for the degree of each generator of the Cox ring, with edges
  corresponding to intersection numbers, and with circles and squares
  marking classes of $(-2)$- and $(-1)$-curves as in
  Step~\ref{step:types}).  In case of toric generalized del Pezzo
  surfaces, the extended Dynkin diagram of Cox ring generators is
  ``circular'' by Remark~\ref{rem:toric}.
\end{result}

\begin{remark}
  In degree $d=1$, Lemma~\ref{lem:base_point_free} turns out to be false in
  some cases, and we do not know whether Theorem~\ref{thm:generators} is
  always true. However, working in the following
  Step~\ref{step:relation} under the assumption that
  Theorem~\ref{thm:generators} still holds gives the correct result, as we can
  confirm at the end using the techniques of \cite{MR2029868} or
  \cite{MR2671185}. The latter applies only to $\tS$ with a suitable
  $\KK^*$-action. For every type of generalized del Pezzo surface where we
  find $13-d$ generators, there is one isomorphy class with such a
  structure. By the proof of Lemma~\ref{lem:generators} and
  Lemma~\ref{lem:dimension_nef}, other isomorphy classes of the same type have
  the same number and degrees of generators.
\end{remark}

\subsection{Construction of generalized del Pezzo surfaces}\label{sec:construction}

The type does not determine the surface up to isomorphism. For
example, the type of an ordinary del Pezzo surface depends only on the
degree $d$ (except for $d=8$, where $\Pone \times \Pone$ and a blow-up
of $\Ptwo$ in one point have different types), but it is unique up to
isomorphism only if $d \ge 5$. For each type where the Cox rings have
one relation, it turns out that there are only one or two isomorphism
classes. Two isomorphism classes occur only for type $\Dfour$ in
degree $3$, for types $\Dfive+\Aone$ and $\Esix$ in degree $2$ and for
types $\Esix+\Atwo$, $\Eseven+\Aone$ and $\Eeight$ in degree $1$ (see
\cite{MR80f:14021} and \cite{MR1933881}).

So far, everything is purely combinatorial and depends only on the
type, but not on the isomorphy class of $\tS$. A good choice of the
relation in $\Cox(\tS)$ will determine $\tS$.

\begin{lemma}\label{lem:relation}
  Let $\tS$ be a generalized del Pezzo surface of degree $d \le 7$
  whose Cox ring has $13-d$ generators $\e_1, \dots, \e_{13-d}$ with
  corresponding divisors $E_1, \dots, E_{13-d}$, and one relation
  $R(\e_1, \dots, \e_{13-d})=0$. Let $\rho: \tS \to \Ptwo$ be the
  contraction of (the strict transforms of) $E_i$ for $i \in I$ with
  $\#I=9-d$.

  Then we may choose coordinates $y_0,y_1,y_2$ in $\Ptwo$ such that
  $\rho(E_{i_j})=\{y_j=0\}$ for $j=0,1,2$ and such that $\rho(E_{i_3})
  = \{R'(y_0,y_1,y_2)=0\}$ for some homogeneous polynomial $R'$ of
  degree $(\cl E_{i_3},\ell_0)$, with $I \cup \{i_0, \dots, i_3\} = \{1,
  \dots, 13-d\}$.

  Furthermore, we may normalize $y_0,y_1,y_2$ and $R'$ such that the
  substitution of the variables $\e_{i_j}$ by $y_j$, for
  $j=0,\dots,3$, and $\e_i$ by $1$, for $i \in I$, in the relation
  $R(\e_1, \dots, \e_{13-d})=0$ results in $y_3-R'(y_0,y_1,y_2)=0$.
\end{lemma}

\begin{proof}
  As described in Section~\ref{sec:del_pezzo}, any such $\tS$ is
  obtained as a composition of $9-d$ blow-ups of $\Ptwo$ as
  in~(\ref{eq:blow_up_sequence}), and the strict transforms of the
  exceptional divisors are negative curves on $\tS$ (corresponding to
  generators $\e_i$ of $\Cox(\tS)$ satisfying
  $(\deg(\e_i),\ell_0)=0$, for $i \in I$ of cardinality $9-d$).

  For any choice of independent sections $y_0,y_1,y_2 \in
  H^0(\Ptwo,\OO(1))$, the sections $\rho^*(y_j)$ have degree $\ell_0$
  and are linearly independent. Looking at the three possibilities for
  $\cl E_i$ as a linear combination of $\ell_0, \dots, \ell_r$ listed in
  Step~\ref{step:types} (where the third option also holds
  for the nef degrees of the additional generators), we see that the
  coefficient $a_{i,0}=(\deg(\e_i),\ell_0)$ of $\ell_0$ must be $1$
  for at least three indices $i \in \{i_0,i_1,i_2\}$. Then we may
  choose $y_0,y_1,y_2$ such that
  $\rho(E_{i_j})=\{y_j=0\}$. Furthermore, for the remaining index
  $i_3$, the curve $\rho(E_{i_3})$ has degree
  $a_{i_3,0}=(\deg(\e_{i_3}),\ell_0)$, so it is described by a
  homogeneous polynomial $R'(y_0,y_1,y_2)$ of that degree.

  Let $y_3 \in H^0(\Ptwo,\OO(a_{i_3,0}))$ be a section vanishing on
  $\rho(E_{i_3})$. For $j=0, \dots, 3$, we normalize $y_j$ such that
  $\rho^*(y_j)$ is a product of $\e_{i_j}$ and a monic monomial in
  $\{\e_i \mid i \in I\}$, and we normalize $R'$ such that
  $y_3-R'(y_0,y_1,y_2)=0$ is a relation in $H^0(\Ptwo,\OO(a_{i_3,0}))$.

  Therefore, its pullback via $\rho^*$ is a relation of degree
  $a_{i_3,0}\ell_0$ in $\Cox(\tS)$. But by assumption, any relation in
  $\Cox(\tS)$ is a product of $R(\e_1, \dots, \e_{13-d})=0$ and some
  polynomial expression $R''(\e_1, \dots \e_{13-d})$. Looking at the
  term $\rho^*(y_3)$ of $\rho^*(y_3-R'(y_0,y_1,y_2))$ and the
  $\ell_0$-component of the degrees, we see that $R''$ can only
  involve $\{\e_i \mid i \in I\}$. The final claim follows by our
  choice of normalizations.
\end{proof}

If the Cox ring of some type of generalized del Pezzo surfaces of
degree $d \le 7$ has $13-d$ generators and one relation, fixing an
explicit polynomial expression for the relation determines $\tS$
uniquely. Indeed, this is abstractly clear, because $\tS$ is the GIT
quotient of $\Spec(\Cox(\tS))$ by the action of $\Tns$ determined by
the degrees of the generators of $\Cox(\tS)$.

But we can also obtain an explicit description of $\tS$ as a blow-up $\rho :
\tS \to \Ptwo$ of the projective plane from the relation $R$ and from the
degrees of the generators of $\Cox(\tS)$ found in
Step~\ref{step:generator_degrees}. We choose coordinates in $\Ptwo$ as in
Lemma~\ref{lem:relation} and determine the images of the divisors $E_i$
corresponding to the generators $\e_i$ of $\Cox(\tS)$ that are not contracted
by $\rho$.  The multiplicity of $p_j$ on the strict transform of $E_i$ is
$a_{i,j}=(\cl E_i,\ell_j)$, which is easily read off the expressions for all
$\cl E_i$ in terms of $\ell_0, \dots, \ell_r$ found in
Step~\ref{step:generator_degrees}. It turns out that $p_1, \dots, p_{9-d}$ are
always determined uniquely by these multiplicities and the images $\rho(E_i)$
of the non-contracted $E_i$.

We may also determine $\rho(E_i)$ for the contracted $E_i$ ($i \in I$)
as follows. As described in the proof of Lemma~\ref{lem:relation},
$\rho^*(y_j)$ is the product of $\e_{i_j}$ and a monic monomial in all
$\e_i$ with $i \in I$. This monomial can be determined just from the
knowledge of the degrees of all $\e_i$ and $\rho^*(y_j)$. If $\e_i$
occurs in the resulting expression for $\rho^*(y_j)$, then $\rho(E_i)$
must lie on $\rho(E_{i_j})$. It turns out that this always determines
$\rho(E_i)$ uniquely.

Therefore, the relation $R$ also determines the generators $\e_1, \dots,
\e_{13-d}$ of $\Cox(\tS)$ up to constant factors: for $i \in I$, the section
$\e_i$ corresponds to the strict transform of the exceptional divisor of the
corresponding blow-up; for $i \in \{i_0, \dots, i_3\}$, the section $\e_i$
corresponds to the strict transform of $\rho(E_i)$.

\begin{step}[Determine the relation in the Cox ring]\label{step:relation}
  As there is a relation in degree $\cl D_0$, there are $h^0(\tS,\OO(D_0))+1$
  ways to write $\cl D_0$ as a non-negative integral linear combination of
  $\cl E_1, \dots, \cl E_{13-d}$. Therefore, the relation $R$ is a linear
  combination of the corresponding $h^0(\tS,\OO(D_0))+1$ monomials in the
  generators $\e_1, \dots, \e_{13-d}$ of $\Cox(\tS)$.

  To construct an explicit description of $\tS$ and $\Cox(\tS)$, we
  choose a candidate for the relation (see also
  Example~\ref{ex:relation_choice}; confirming the correctness of the
  choice of relation is part of Step~\ref{step:model} below) for each
  isomorphism class of $\tS$ of each type and determine $\rho(E_i)$
  for all $i=1, \dots, 13-d$ from it.
\end{step}

\begin{result}
  We list the relation $R(\e_1, \dots, \e_{13-d})$ and the images of the
  curves $E_i$ corresponding to the generators $\e_i$ of $\Cox(\tS)$ unter
  $\rho : \tS \to \Ptwo$. Here, $\rho(E_i)$ is a point for $i \in I$,
  $\rho(E_{i_j}) = \{y_j=0\}$ for $j=0,1,2$, and $\rho(E_{i_3}) =
  \{R'(y_0,y_1,y_2)=0\}$.
\end{result}

\begin{example}\label{ex:relation_choice}
  Different choices of a relation $R$ may lead to isomorphic or
  non-isomorphic $\tS$ (as mentioned in
  Section~\ref{sec:del_pezzo}). Furthermore, it is not possible to
  choose an arbitrary relation between the monomials in $\e_1, \dots,
  \e_{13-d}$ of degree $\cl D_0$.

  For example, in case of type $\Dfour$ in degree $d=3$, the monomials
  of the degree $\ell_0$ of the relation are $\e_2\e_5^2\e_8,
  \e_3\e_6^2\e_9, \e_4\e_7^2\e_{10}, \e_1\cdots\e_7$. The curves
  contracted by $\rho$ are $E_2, \dots, E_7$, and $\rho(E_1),
  \rho(E_8), \rho(E_9), \rho(E_{10})$ are lines in $\Ptwo$. The two
  non-isomorphic choices that can be found in
  Section~\ref{sec:degree_3} correspond to the configurations where
  $\rho(E_1)$ meets the other lines in distinct points, but these may
  or may not meet each other in the same point. However, one cannot
  choose the relation $\e_2\e_5^2\e_8+\e_3\e_6^2\e_9+\e_1\cdots\e_7=0$
  because this would lead to $\rho(E_1), \rho(E_8), \rho(E_9)$ meeting
  in one point, which is impossible because the first blow-up of $p_1$
  must be on $\rho(E_1)\cap \rho(E_8)$ but not on $\rho(E_9)$.
\end{example}

\subsection{Singular del Pezzo surfaces}\label{sec:singular}

Singular del Pezzo surfaces of degree $3$ were classified by
Schl\"afli \cite{schlaefli} and Cayley \cite{cayley}; see
\cite{MR80f:14021} for a more modern treatment. In degree $4$, see
\cite{MR13:972c}. All types of degree $\ge 4$ can also be found in
\cite{MR89f:11083}. A classification of del Pezzo surfaces of degrees
$2$ and $1$ was done by Du Val \cite{duval}, see also \cite{MR713283}.
A modern presentation of these results can be found in \cite{MR767407}
and \cite{MR2227002}.

For a generalized del Pezzo surface $\tS$ of degree $d$, the linear system
$|-K_\tS|$ is $d$-dimensional. Let
\begin{equation*}
  \PP[d]=
  \begin{cases}
    \PP^d, & d \ge 3,\\
    \PP(2,1,1,1), & d=2,\\
    \PP(3,2,1,1), & d=1.
  \end{cases}
\end{equation*}
For any $d$, we have a morphism $\pi: \tS \to \PP[d]$, defined as follows:
\begin{itemize}
\item For $d \ge 3$, the morphism $\pi$ is defined by the choice of a
  basis $x_0, \dots, x_d \in H^0(\tS,\OO(-K_\tS))$.
\item For $d = 2$, it is defined by a basis $x_1,x_2,x_3 \in H^0(\tS,
  \OO(-K_\tS))$ and some $x_0 \in H^0(\tS, -2K_\tS)$ linearly
  independent of $x_1^2, x_2^2, x_3^2, x_1x_2,x_1x_3,x_2x_3$.
\item For $d = 1$, it is defined by suitable $x_0 \in H^0(\tS,
  -3K_\tS)$, $x_1 \in H^0(\tS, -2K_\tS)$ and a basis $x_2,x_3 \in
  H^0(\tS, \OO(-K_\tS))$.
\end{itemize}

Let $S \subset \PP[d]$ be the image of $\tS$ under $\pi$. If $\tS$ is
ordinary, $\pi$ is an embedding, so $S \cong \tS$. Otherwise, $\pi$
contracts the $(-2)$-curves on $\tS$ to rational double points (i.e.,
$\ADE$-singularities) on $S$ and is an embedding in their
complement. More precisely, the preimage under $\pi$ of a singularity
$p$ of $\ADE$-type on $S$ is a connected tree of $(-2)$-curves on
$\tS$ such that the corresponding connected component of $(-2)$-curves
in the extended Dynkin diagram of negative curves has the same $\ADE$-type as $p$. The
image of a $(-1)$-curve on $\tS$ under $\pi$ is a line on $S$ if
$d \ge 3$. For simplicity, we use the term \emph{line on $S$} also for
the image of a $(-1)$-curve on $\tS$ under $\pi$ for $d=1,2$.

Therefore, $S$ is a singular del Pezzo surface with minimal desingularization
$\tS$. We have $-\cl K_\tS = \pi^*(-\cl K_S)$.

The \emph{type} of a singular del Pezzo surface $S$ is the type of its minimal
desingularization $\tS$. As in Remark~\ref{rem:type_classification}, the type
of $S$ can be described by its degree, the types of its $\ADE$-singularities
and the number of its lines.

A singular del Pezzo surface of degree $d$ can be described explicitly as
follows:
\begin{itemize}
\item For $d \ge 5$, the singular del Pezzo surface $S$ is the
  non-complete intersection of $d(d-3)/2$ quadrics in $\PP^d$
  \cite[Theorem~4.4]{MR646042}.
\item For $d=4$, it is the complete intersection of two quadrics in $\PP^4$.
\item For $d=3$, it is a cubic surface in $\PP^3$.
\item For $d=2$, it is a surface in $\PP(2,1,1,1)$ of weighted degree
  $4$.
\item For $d=1$, it is a surface in $\PP(3,2,1,1)$ of weighted degree
  $6$.
\end{itemize}

The composition of a rational inverse of $\pi : \tS \to S$ and $\rho :
\tS \to \Ptwo$ is a birational map $\phi : S \rto \Ptwo$, resulting in
the following diagram:
\begin{equation*}
  \xymatrix{\tS \ar@{->}^\rho[dr]\ar@{->}^\pi[d] & \\
    \PP[d] \supset S\ar@{-->}^\phi[r] & \Ptwo}
\end{equation*}

\begin{step}[Determine the singular del Pezzo surface]\label{step:model}
  For this, we determine all monomials in $\e_1, \dots, \e_{13-d}$ of
  degree $-\cl K_\tS$. Keeping the relation $R$ in mind, we choose
  monomials $M_0, \dots, M_d$ that form a basis of
  $H^0(\tS,\OO(-K_\tS))$. For $d \ge 3$, we define $\pi: \tS \to \PP^d$ via
  $\pi^*(x_i) = M_i$, for $i=0, \dots, d$. We choose additionally
  suitable independent monomials of degree $-2K_\tS$ if $d=1,2$ and
  $-3K_\tS$ if $d=1$ and get $\pi: \tS \to \PP[d]$ by defining
  $\pi^*(x_i)$ to be a chosen monomial in $H^0(\tS,-d_iK_\tS)$ for
  each coordinate $x_i$ of weighted degree $d_i$ in
  $\PP[d]=\PP(d_0,d_1,d_2,d_3)$.

  Using $R$, we express any other monomial $M$ of degree $-\cl K_\tS$ as a
  linear combination $M=\sum_{i=0}^d \alpha_iM_i$ of $M_0, \dots,
  M_d$. Then $M=\pi^*(\sum_{i=0}^d \alpha_ix_i)$.

  Next, we want to find equations defining $S$ in $\PP[d]$. If $F(x_0,
  \dots, x_d)=0$ on $S$ for a (weighted) homogeneous polynomial of
  (weighted) degree $d_0$, then also $F(\pi^*(x_0), \dots,
  \pi^*(x_d))=0$, giving a relation of degree $-d_0K_\tS$ in
  $\Cox(\tS)$. As the ideal of relations is generated by $R$ by
  assumption, this must be a multiple of $R$.

  As described above, we know how many independent equations of which
  degree define $S$. Therefore, we can write down all monomials in
  $x_0, \dots, x_d$ of that degree, apply $\pi^*$ and find relations
  between the resulting monomials in $\e_1, \dots, \e_{13-d}$ using
  $R$ until the corresponding homogeneous equations in $x_0, \dots,
  x_d$ are enough to define~$S$.

  To confirm that we have chosen a suitable relation in $\Cox(\tS)$ in
  Step~\ref{step:relation} that gives rise to a singular del Pezzo
  surface $S$ of the expected type, we compute the minimal
  desingularization of $S$ by finding and blowing up the
  singularities, and we determine the configuration of the
  $(-2)$-curves obtained from the singularities and the $(-1)$-curves
  obtained as the strict transforms of the lines on $S$. We compare
  the resulting extended Dynkin diagram of negative curves with the
  one computed in Step~\ref{step:types}.  If we have two isomorphism
  classes of the same type, we find a geometric property that shows
  that our construction gives two non-isomorphic surfaces of this
  type.  If this reveals that we have made an unsuitable choice of
  relation, we go back to Step~\ref{step:relation}.

  Next, we determine $\pi(E_i)$ for $i=1, \dots, 13-d$. Here, we use that
  $\e_i$ occurs in the monomial $\pi^*(x_j)$ if and only if $x_j$ vanishes on
  $\pi(E_i)$, and similarly for $\pi^*$ of the linear combinations of $x_0,
  \dots, x_d$ mentioned above. 

  This gives another description of the divisors corresponding to the
  generators of $\Cox(\tS)$. If one computes $\tS$ from $S$ by
  resolving the singularities, one can determine which exceptional
  divisor is which $(-2)$-curve using the extended Dynkin diagram of
  Cox ring generators.

  Finally, to find an explicit expression for $\phi: S \rto \Ptwo$, we
  first constructs its inverse $\psi: \Ptwo \rto S$, which is defined
  by $(y_0:y_1:y_2)\mapsto \xx \in S \subset \PP[d]$. Here, the $i$-th
  component $x_i$ of $\xx$ is obtained from the expression
  $\pi^*(x_i)$ in $\e_1, \dots, \e_{13-d}$ by substituting $\e_i$ by
  $1$ for $i \in I$, and $\e_{i_j}$ by $y_j$ for $j=0, \dots, 3$
  (where we rewrite $y_3$ in terms of $y_0,y_1,y_2$ using
  $y_3-R'(y_0,y_1,y_2)=0$, which is obtained from $R$ by the same
  substitutions; see Lemma~\ref{lem:relation}). If $x_i$ has weighted
  degree $d_i$, then the expression in terms of $y_0,y_1,y_2$ has
  degree $3d_i$.

  Then it is straightforward to compute the inverse $\phi$. For $d \ge
  3$, we choose the sequence of contractions in the Dynkin diagram of
  negative curves in Step~\ref{step:types} such that $\phi$ is a
  projection to three coordinates $x_i$. For $d=2$, $\phi$ is given by
  monomials of weighted degree $2$. For $d=1$, it is given by
  monomials of weighted degree $3$.
\end{step}

\begin{result}\label{res:model}
  We list the equations defining $S \subset \PP[d]$. We describe the
  anticanonical map $\pi: \tS \to \PP[d]$ by expressing
  $\pi^*(x_i)$ in terms of $\e_1, \dots, \e_{13-d}$. We also list all
  other $\pi^*(\sum \alpha_i x_i)$ that are monomials in $\e_1, \dots
  ,\e_{13-d}$.

  We list $\pi(E_1), \dots, \pi(E_{13-d})$: for a $(-2)$-curve, this is a
  point, namely the singularity it arises from (for which we give the
  $\ADE$-type); for a $(-1)$-curve, it is a line (which we describe as an
  intersection of $d-1$ hyperplanes if $d \ge 3$); for a nef curve, it is a
  curve of higher degree (which we describe as an intersection of
  hypersurfaces).

  Whenever more than two $E_i$ intersect pairwise on $\tS$ (which
  results in a triangle in the extended Dynkin diagram of Cox ring
  generators in case of three pairwise intersecting curves), we also
  determine whether they meet in one common point (sometimes, this
  helps to distinguish different isomorphy classes of $\tS$; it is
  also interesting for applications). If so, we compute its image
  under $\pi$.

  We list the projection $\phi : S \rto \Ptwo$.  We list its inverse $\psi :
  \Ptwo \rto S$ only for $d=1,2$, as it is easily recovered from $\phi$ for $d
  \ge 3$.
\end{result}

\section{Classification}\label{sec:classification}

In this section, we present our classification of Cox rings of
generalized del Pezzo surface, as summarized in
Theorem~\ref{thm:classification}. It is obtained using the results of
Section~\ref{sec:strategy}, following the strategy of
Steps~\ref{step:types}--\ref{step:model}. The presentation of del
Pezzo surfaces where the Cox ring has one relation is as described in
Results~\ref{res:types}--\ref{res:model} in
Section~\ref{sec:strategy}.

\subsection{Degree 7, 8, 9}\label{sec:degree_789}

Ordinary del Pezzo surfaces are $\Ptwo$ of degree $9$, the blow-up
$\Bl_1\Ptwo$ of $\Ptwo$ in one point and $\Pone \times \Pone$ of
degree $8$, and the blow-up $\Bl_2\Ptwo$ of $\Ptwo$ in two distinct
points of degree $7$. The additional generalized del Pezzo surfaces
are the Hirzebruch surface $F_2$ of degree $8$ and the blow-up of
$\Bl_1\Ptwo$ in a point on its exceptional divisor. 

It is easy to check that all of them are toric (see
Remark~\ref{rem:toric}).  We list them together with some basic
properties and their cyclic Dynkin diagrams of Cox ring generators in
Table~\ref{tab:degree_789}.

\begin{table}[ht]
  \centering
\[\begin{array}{cccccc}
  \hline
  \text{degree} & \text{type} & \text{singularities} 
  & \text{$(-1)$-curves} & \text{type} & \text{Dynkin diagram} \\
  \hline
  9 & \Ptwo & - & 0 & \toric & (1,1,1) \\
  8 & \Bl_1\Ptwo & - & 1 & \toric & (1,0,-1,0) \\
  8 & \Pone \times \Pone & - & 0 & \toric & (0,0,0,0) \\
  8 & F_2 & \Aone & 0 & \toric & (2,0,-2,0) \\
  7 & \Bl_2\Ptwo & - & 3 & \toric & (0,-1,-1,-1,0) \\
  7 & & \Aone & 2 & \toric & (1,0,-2,-1,-1) \\
  \hline
\end{array}\]
  \smallskip
  \caption{Del Pezzo surfaces of degree $\ge 7$}
  \label{tab:degree_789}
\end{table}

\subsection{Degree 6}\label{sec:degree_6}

The classification of generalized del Pezzo surfaces $\tS$ of degree $6$ can
be found in \cite[Proposition~8.3]{MR89f:11083}, see
Table~\ref{tab:degree_6}. Ordinary sextic del Pezzo surfaces and types $\Aone$
with four lines, $2\Aone$ and $\Atwo+\Aone$ are toric, with cyclic Dynkin
diagrams of Cox ring generators
\begin{equation*}
  \begin{split}
    &(-1,-1,-1,-1,-1,-1), \qquad (0,-1,-1,-2,-1,-1),\\
    &(0,-2,-1,-2,-1,0),\qquad (1,0,-2,-2,-1,-2),
  \end{split}
\end{equation*}
For every type, considering the configuration of blown-up points shows that
the corresponding del Pezzo surface is unique up to isomorphism. For the
remaining types, we determine $\Cox(\tS)$ using the procedure described in
Section~\ref{sec:strategy}.

\begin{table}[ht]
  \centering
  \begin{equation*}
    \begin{array}{cccc}
      \hline
      \text{singularities} 
      & \text{lines} & \text{type} & \text{reference}\\
      \hline
      - & 6 & \toric & \text{\cite{MR1620682}}\\
      \Aone & 4 & \toric & \text{\cite{MR1620682}}\\
      \Aone & 3 & \onerel & \text{\cite{MR1906155}, \cite{MR2559866}}\\
      2\Aone & 2 & \toric & \text{\cite{MR1620682}}\\
      \Atwo & 2 & \onerel & \text{\cite{MR1906155}, \cite{MR2769338}}\\
      \Atwo+\Aone & 1 & \toric & \text{\cite{MR1620682}}\\
      \hline
    \end{array}
  \end{equation*}
  \smallskip
  \caption{Del Pezzo surfaces of degree 6}
  \label{tab:degree_6}
\end{table}

\coxcase{ii}{$\Aone$}\coxembed{\Psix}{
    {}&x_0x_1-x_3x_6 = x_0x_2-x_4x_6 = x_0x_5-x_3x_4 = x_1x_4-x_2x_3\\
    = {}&x_2x_3-x_5x_6 = x_0x_1+x_0x_2+x_1x_2 = x_0x_3+x_0x_4+x_2x_3\\
    = {}&x_1x_5+x_3^2+x_3x_4 = x_2x_5+x_3x_4+x_4^2 = 0}
\coxbirat 6 {x_3:x_4:x_5}  
\coxtwoone{\Aone}{(0:0:0:0:0:0:1)}{E_1}{
  \coximod{1}{\ell_0-\ell_1-\ell_2-\ell_3}{(0:0:0:0:0:0:1)}{y_2}
}
\coxonem{E_2, E_3, E_4}{
  \coximp{2}{\ell_1}{x_0=x_1=x_3=x_4=x_5=0}{(0:1:0)}
  \coximp{3}{\ell_2}{x_0=x_2=x_3=x_4=x_5=0}{(1:0:0)}
  \coximpd{4}{\ell_3}{x_1=x_2=x_3=x_4=x_5=0}{(1:-1:0)}
}
\coxother{$(0)$-curves $E_5, E_6, E_7$}{
  \coxims{5}{\ell_0-\ell_1}{x_0=x_1=x_3=x_6=x_2x_5+x_4^2=0}{y_0}
  \coxims{6}{\ell_0-\ell_2}{x_0=x_2=x_4=x_6=x_1x_5+x_3^2=0}{y_1}
  \coximds{7}{\ell_0-\ell_3}{x_1=x_2=x_3+x_4=x_6=x_0x_5-x_3x_4=0}{-y_0-y_1} 
}
\coxdynkinp{E_5 \cap E_6 \cap E_7}{(0:0:0:0:0:1:0)}{E_5
  \ar@{-}[rr]\ar@{-}[dd]\ar@{-}[dr]& & \li{2} \ar@{-}[dr] \\ & E_6\ar@{-}[r] &
  \li{3} \ar@{-}[r] & \ex{1}\\E_7\ar@{-}[rr]\ar@{-}[ur] & & \li{4}\ar@{-}[ur]}
\coxonerel{\e_1, \dots, \e_7}{\e_2\e_5+\e_3\e_6+\e_4\e_7}{\ell_0} 
\coxantie{(\e_1\e_2\e_3\e_5\e_6, \e_1\e_2\e_4\e_5\e_7, \e_1\e_3\e_4\e_6\e_7,
  \e_1^2\e_2^2\e_3\e_4\e_5,\e_1^2\e_2\e_3^2\e_4\e_6, \e_1^3\e_2^2\e_3^2\e_4^2,
  \e_5\e_6\e_7)}{ \pi^*(-x_0-x_1) = \e_1\e_2^2\e_5^2,
  \ \pi^*(-x_0-x_2) = \e_1\e_3^2\e_6^2,\\
  \pi^*(-x_1-x_2) = \e_1\e_4^2\e_7^2, \ \pi^*(-x_3-x_4) =
  \e_1^2\e_2\e_3\e_4^2\e_7}

\coxcase{iv}{$\Atwo$}
\coxembed{\Psix}{{}&x_0x_5-x_3x_4 = x_0x_6-x_1x_4 = x_0x_6-x_2x_3 = x_3x_6-x_1x_5\\
    = {}&x_4x_6-x_2x_5 = x_1x_6+x_3^2+x_3x_4 = x_2x_6+x_3x_4+x_4^2\\
    = {}&x_6^2+x_3x_5+x_4x_5 = x_1x_2+x_0x_3+x_0x_4 = 0}
\coxbirat 6 {x_4:x_5:x_6}
\coxtwoone{\Atwo}{(1:0:0:0:0:0:0)}{E_1, E_2}{
  \coximop{1}{\ell_2-\ell_3}{(1:0:0:0:0:0:0)}{(1:0:0)}
  \coximopd{2}{\ell_1-\ell_2}{(1:0:0:0:0:0:0)}{(1:0:0)}
}
\coxonem{E_3, E_4}{
  \coximp{3}{\ell_3}{x_1=x_3=x_4=x_5=x_6=0}{(1:0:0)}
  \coximd{4}{\ell_0-\ell_1-\ell_2}{x_2=x_3=x_4=x_5=x_6=0}{y_1}
}
\coxother{$(0)$-curve $E_5$ and $(1)$-curves $E_6, E_7$}{
  \coxims{5}{\ell_0-\ell_1}{x_0x_5-x_3^2=x_1=x_2=x_3+x_4=x_6=0}{y_2}
  \coxims{6}{2\ell_0-\ell_1-\ell_2-\ell_3}{x_0=x_1=x_3=x_2x_5-x_4x_6=x_2x_6+x_4^2=x_4x_5+x_6^2=0}{-y_0y_1-y_2^2}
  \coximds{7}{\ell_0}{x_0=x_2=x_4=x_3x_6-x_1x_5=x_1x_6+x_3^2=x_3x_5+x_6^2=0}{y_0}
}
\coxdynkinp{E_5 \cap E_6 \cap E_7}
{(0:0:0:0:0:1:0)}{E_6\ar@{-}[rr]\ar@{=}[dd]\ar@{-}[dr] & &
  \li{3}\ar@{-}[dr]\\ &
  E_5\ar@{-}[r] & \ex{2} \ar@{-}[r] & \ex{1}\\
  E_7\ar@{-}[rr]\ar@{-}[ur] & & \li{4}\ar@{-}[ur]} 
\coxonerel{\e_1, \dots, \e_7}{\e_2\e_5^2 +\e_3\e_6 +
  \e_4\e_7}{2\ell_0-\ell_1-\ell_2} 
\coxantie{ (\e_6\e_7, \e_1\e_2\e_3\e_5\e_6,
  \e_1\e_2\e_4\e_5\e_7,\e_1^2\e_2\e_3^2\e_4\e_6,
  \e_1^2\e_2\e_3\e_4^2\e_7, \e_1^4\e_2^2\e_3^3\e_4^3,
  \e_1^3\e_2^2\e_3^2\e_4^2\e_5)}{\pi^*(-x_1-x_2) = \e_1\e_2^2\e_5^3,\
  \pi^*(-x_3-x_4) = \e_1^2\e_2^2\e_3\e_4\e_5^2}

\subsection{Degree 5}\label{sec:degree_5}

See \cite[Proposition~8.4]{MR89f:11083} for the classification of generalized
quintic del Pezzo surface. The generators and relations in the Cox ring of
ordinary quintic del Pezzo surfaces have been determined in \cite{MR2029863},
see also \cite{MR1260765}; because of the ten lines, there are at least two
(in fact five) relations. Additionally, we have six types of generalized del
Pezzo surfaces as in Table~\ref{tab:degree_5}.
  
\begin{table}[ht]
  \centering
  \begin{equation*}
    \begin{array}{cccc}
      \hline
      \text{singularities} 
      & \text{lines} & \text{type} & \text{reference}\\
      \hline
      - & 10 & \tworel & \text{\cite{MR1909606}, \cite{MR2099200}}\\
      \Aone & 7 & \onerel & \text{\cite{bauer_diplom}}\\
      2\Aone & 5 & \toric & \text{\cite{MR1620682}} \\
      \Atwo & 4 & \onerel & \text{\cite{arXiv:0710.1583}} \\
      \Atwo+\Aone & 3 & \toric & \text{\cite{MR1620682}} \\
      \Athree & 2 & \onerel & \text{\cite{MR1906155}} \\
      \Afour & 1 & \onerel & \text{\cite{MR1906155}, \cite{moll_diplom}} \\
      \hline
    \end{array}
  \end{equation*}
  \smallskip
  \caption{Del Pezzo surfaces of degree 5}
  \label{tab:degree_5}
\end{table}

Types $2\Aone$ and $\Atwo+\Aone$ are toric, with cyclic Dynkin diagrams of Cox ring generators
\begin{equation*}
  (-1,-1,-1,-1,-2,-1,-2), \qquad (0,-1,-1,-2,-2,-1,-2).
\end{equation*}
The remaining types have Cox rings with precisely one relation and are
described in detail below. For each type, there is a quintic del Pezzo surface
that is unique up to isomorphism, as one can see from the
configurations of blown-up points.

\coxcase{i}{$\Aone$}
\coxembed{\Pfive}{
    {}&x_0x_4-x_1x_2 = x_0x_5-x_1x_3 = x_2x_5-x_3x_4\\ =
    {}&x_1x_2+x_1x_3+x_2x_3 = x_1x_4+x_1x_5+x_2x_5 = 0}
\coxbirat 5 {x_1:x_2:x_4} 
\coxtwoone{\Aone}{(1:0:0:0:0:0)}{E_1}{
  \coximod{1}{\ell_0-\ell_1-\ell_2-\ell_3}{(1:0:0:0:0:0)}{y_2}
}
\coxonem{E_2, \dots, E_8}{
  \coximp{2}{\ell_1}{x_1=x_2=x_4=x_5=0}{(1:-1:0)}
  \coximp{3}{\ell_2}{x_1=x_3=x_4=x_5=0}{(0:1:0)}
  \coximp{4}{\ell_3}{x_2=x_3=x_4=x_5=0}{(1:0:0)}
  \coximp{5}{\ell_4}{x_0=x_1=x_2=x_3=0}{(0:0:1)}
  \coxim{6}{\ell_0-\ell_1-\ell_4}{x_0=x_1=x_2=x_4=0}{-y_0-y_1}
  \coxim{7}{\ell_0-\ell_2-\ell_4}{x_0=x_1=x_3=x_5=0}{y_0}
  \coximd{8}{\ell_0-\ell_3-\ell_4}{x_0=x_2=x_3=x_4+x_5=0}{y_1}
}
\coxdynkin{ & \li 2 \ar@{-}[r]& \li 6\ar@{-}[dr]\\
  \ex{1} \ar@{-}[r]\ar@{-}[ur] 
  \ar@{-}[dr]& \li 3\ar@{-}[r] & \li 7\ar@{-}[r] & \li 5\\
  & \li 4\ar@{-}[r] & \li 8\ar@{-}[ur]}
\coxonerel{\e_1, \dots,
  \e_8}{\e_2\e_6+\e_3\e_7+\e_4\e_8}{\ell_0-\ell_4}
\coxantie{(\e_5^2\e_6\e_7\e_8,
  \e_1\e_2\e_3\e_5\e_6\e_7, \e_1\e_2\e_4\e_5\e_6\e_8,\e_1\e_3\e_4\e_5\e_7\e_8,
  \e_1^2\e_2^2\e_3\e_4\e_6, \e_1^2\e_2\e_3^2\e_4\e_7)}
{ \pi^*(-x_1-x_2) = \e_1\e_2^2\e_5\e_6^2,\ \pi^*(-x_1-x_3) = \e_1\e_3^2\e_5\e_7^2,\\
  \pi^*(-x_2-x_3) = \e_1\e_4^2\e_5\e_8^2,\ \pi^*(-x_4-x_5) = \e_1^2\e_2\e_3\e_4^2\e_8}

\coxcase{iii}{$\Atwo$} 
\coxembed{\Pfive}{ {}&x_0x_2-x_1x_5 = x_0x_2-x_3x_4 =
  x_0x_3+x_1^2+x_1x_4\\ = {}&x_0x_5+x_1x_4+x_4^2 = x_3x_5+x_1x_2+x_2x_4 = 0}
\coxbirat 5 {x_0:x_1:x_4}
\coxtwoone{\Atwo}{(0:0:1:0:0:0)}{E_1, E_2}{
  \coximo{1}{\ell_0-\ell_3-\ell_4}{(0:0:1:0:0:0)}{y_0}
  \coximopd{2}{\ell_1-\ell_2}{(0:0:1:0:0:0)}{(0:1:-1)}
}
\coxonem{E_3, \dots, E_6}{
  \coximp{3}{\ell_3}{x_0=x_1=x_3=x_4=0}{(0:0:1)}
  \coximp{4}{\ell_4}{x_0=x_1=x_4=x_5=0}{(0:1:0)}
  \coximp{5}{\ell_2}{x_0=x_3=x_5=x_1+x_4=0}{(0:1:-1)}
  \coximds{6}{\ell_0-\ell_1-\ell_2}{x_2=x_3=x_5=x_1+x_4=0}{-y_1-y_2}
}
\coxother{$(0)$-curves $E_7, E_8$}{
  \coxim{7}{\ell_0-\ell_3}{x_1=x_2=x_3=x_0x_5+x_4^2=0}{y_1}
  \coximd{8}{\ell_0-\ell_4}{x_2=x_4=x_5=x_0x_3+x_1^2=0}{y_2}
}
\coxdynkinp{E_6 \cap E_7 \cap E_8}{(1:0:0:0:0:0)}{E_7 \ar@{-}[rrr]\ar@{-}[dr]\ar@{-}[dd]
  &&& \li{3}\ar@{-}[dr]\\
  & \li{6} \ar@{-}[r] & \li{5}\ar@{-}[r] & \ex{2}\ar@{-}[r] & \ex{1}\\
  E_8\ar@{-}[rrr]\ar@{-}[ur] &&& \li{4}\ar@{-}[ur]}
\coxonerel{\e_1, \dots, \e_8}{\e_2\e_5^2\e_6 + \e_3\e_7 + \e_4\e_8}{\ell_0}
\coxantie{(\e_1^3\e_2^2\e_3^2\e_4^2\e_5,
  \e_1^2\e_2\e_3^2\e_4\e_7, \e_6\e_7\e_8,\e_1\e_2\e_3\e_5\e_6\e_7, \e_1^2\e_2\e_3\e_4^2\e_8, \e_1\e_2\e_4\e_5\e_6\e_8)}
{\pi^*(-x_1-x_4) = \e_1^2\e_2^2\e_3\e_4\e_5^2\e_6,\ \pi^*(-x_3-x_5) =
  \e_1\e_2^2\e_5^3\e_6^2}

\coxcase{v}{$\Athree$}
\coxembed{\Pfive}{{}&x_0x_2-x_1^2 = x_0x_3-x_1x_4 = x_2x_4-x_1x_3\\ =
        {}&x_2x_4+x_4^2+x_0x_5 = x_2x_3+x_3x_4+x_1x_5 = 0}
\coxbirat 5 {x_0:x_1:x_4}
\coxtwoone{\Athree}{(0:0:0:0:0:1)}{E_1, E_2, E_3}{
  \coximop{1}{\ell_1-\ell_2}{(0:0:0:0:0:1)}{(0:0:1)}
  \coximop{2}{\ell_2-\ell_3}{(0:0:0:0:0:1)}{(0:0:1)}
  \coximod{3}{\ell_0-\ell_1-\ell_2-\ell_4}{(0:0:0:0:0:1)}{y_0}
}
\coxonem{E_4, E_5}{
  \coximp{4}{\ell_4}{x_0=x_1=x_3=x_4=0}{(0:1:0)}
  \coximpd{5}{\ell_3}{x_0=x_1=x_2=x_4=0}{(0:0:1)}
}
\coxother{$(0)$-curves $E_6, E_7$ and $(1)$-curve $E_8$}{
  \coxim{6}{\ell_0-\ell_1}{x_1=x_2=x_3=x_0x_5+x_4^2=0}{y_1}
  \coxim{7}{\ell_0-\ell_4}{x_3=x_4=x_5=x_0x_2-x_1^2=0}{y_2}
  \coximds{8}{2\ell_0-\ell_1-\ell_2-\ell_3}{x_5=x_2+x_4=x_0x_3+x_1x_2=x_0x_4+x_1^2=x_1x_3+x_2^2=0}{-y_0y_2-y_1^2}
}
\coxdynkinp{E_6 \cap E_7 \cap E_8}{(1:0:0:0:0:0)}{E_6\ar@{-}[rr]\ar@{-}[dr]\ar@{-}[dd] & & \ex{1}\ar@{-}[dr]\\
  & E_8\ar@{-}[r] & \li{5} \ar@{-}[r] & \ex{2}\\
  E_7\ar@{-}[r]\ar@{=}[ur] & \li{4}\ar@{-}[r] & \ex{3}\ar@{-}[ur]} 
\coxonerel{\e_1, \dots, \e_8}{\e_1\e_6^2 +
    \e_3\e_4^2\e_7 + \e_5\e_8}{2\ell_0-\ell_1-\ell_2}
  \coxantie{(\e_1^2\e_2^4\e_3^3\e_4^2\e_5^3, \e_1^2\e_2^3\e_3^2\e_4\e_5^2\e_6,
    \e_1^2\e_2^2\e_3\e_5\e_6^2,\e_1\e_2\e_3\e_4\e_6\e_7,
    \e_1\e_2^2\e_3^2\e_4^2\e_5\e_7, \e_7\e_8)}{\pi^*(-x_2-x_4) =
    \e_1\e_2^2\e_3\e_5^2\e_8}

\coxcase{vi}{$\Afour$} 
\coxembed{\Pfive}{ {}&x_0x_2-x_1^2 = x_0x_3-x_1x_4 =
  x_2x_4-x_1x_3\\ = {}&x_1x_2+x_4^2+x_0x_5 = x_2^2+x_3x_4+x_1x_5 = 0}
\coxbirat 5 {x_0:x_1:x_4}
\coxtwoone{\Afour}{(0:0:0:0:0:1)}{E_1, \dots, E_4}{
  \coximop{1}{\ell_1-\ell_2}{(0:0:0:0:0:1)}{(0:0:1)}
  \coximop{2}{\ell_2-\ell_3}{(0:0:0:0:0:1)}{(0:0:1)}
  \coximop{3}{\ell_3-\ell_4}{(0:0:0:0:0:1)}{(0:0:1)}
  \coximod{4}{\ell_0-\ell_1-\ell_2-\ell_3}{(0:0:0:0:0:1)}{y_0}
}
\coxone{E_5}{
  \coximpd{5}{\ell_4}{x_0=x_1=x_2=x_4=0}{(0:0:1)}
}
\coxother{$(0)$-curve $E_6$, $(1)$-curve $E_7$ and $(5)$-curve $E_8$}{
  \coxim{6}{\ell_0-\ell_1}{x_1=x_2=x_3=x_0x_5+x_4^2=0}{y_1}
  \coxims{7}{\ell_0}{x_3=x_4=x_0x_5+x_1x_2=x_0x_2-x_1^2=x_1x_5+x_2^2=0}{y_2}
  \coximds{8}{3\ell_0-\ell_1-\ell_2-\ell_3-\ell_4}{x_5=x_0x_2-x_1^2=x_0x_3-x_1x_4=x_2x_4-x_1x_3=x_1x_2+x_4^2=x_2^2+x_3x_4=0}{-y_1^3-y_0y_2^2}
}
\coxdynkinp{E_6 \cap E_7 \cap E_8}{(1:0:0:0:0:0)}
{E_7\ar@{-}[rr]\ar@3{-}[dr]\ar@{-}[dd] & & \ex{4}\ar@{-}[dr]\\
  & E_8\ar@{-}[r] & \li{5} \ar@{-}[r] & \ex{3}\\
  E_6\ar@{=}[ur]\ar@{-}[r] & \ex{1}\ar@{-}[r] & \ex{2}\ar@{-}[ur]} 
\coxonerel{\e_1, \dots, \e_8}{\e_1^2\e_2\e_6^3 +
  \e_4\e_7^2 + \e_5\e_8}{3\ell_0-\ell_1-\ell_2-\ell_3}
\coxanti{(\e_1^2\e_2^4\e_3^6\e_4^3\e_5^5,
  \e_1^2\e_2^3\e_3^4\e_4^2\e_5^3\e_6, \e_1^2\e_2^2\e_3^2\e_4\e_5\e_6^2,
  \e_1\e_2\e_3\e_4\e_6\e_7, \e_1\e_2^2\e_3^3\e_4^2\e_5^2\e_7, \e_8)}

\subsection{Degree 4}\label{sec:degree_4}

By classical results, which can be found in \cite[Book IV, \S
XIII.11]{MR13:972c}, every del Pezzo surface of degree $4$ is the
intersection of two quadrics in $\Pfour$, given by symmetric $5 \times
5$ matrices $A,B$, where $A$ can be assumed to be non-singular.
Besides the smooth quartic del Pezzo surface, there are $15$ singular
types, which can be distinguished by the \emph{Segre symbol}
(describing the structure of the Jordan form) of $A^{-1}B$. The
extended Dynkin diagrams of negative curves can
be found in \cite[Proposition~6.1]{MR89f:11083}.  In total, this leads
to the $16$ types listed in Table~\ref{tab:degree_4}.

\begin{table}[ht]
  \centering
  \begin{equation*}
    \begin{array}{ccccc}
      \hline
      \text{singularities} & \text{lines} & \text{type} & \text{reference} \\
      \hline
      - & 16 & \tworel & \text{\cite{MR2838351}}\\
      \Aone & 12 & \tworel & \\
      2\Aone & 9 & \tworel & \\
      2\Aone & 8 & \tworel & \text{\cite{MR2874644}, \cite{MR2980925}} \\
      \Atwo & 8 & \tworel & \\
      3\Aone & 6 & \onerel & \text{\cite{MR2853047}} \\
      \Atwo+\Aone & 6 & \onerel & \text{\cite{MR2853047}} \\
      \Athree & 5 & \onerel & \text{\cite{arXiv:1302.6151}} \\
      \Athree & 4 & \tworel & \text{\cite{MR2961294}} \\
      4\Aone & 4 & \toric & \text{\cite{MR1620682}}\\
      \Atwo+2\Aone & 4 & \toric & \text{\cite{MR1620682}}\\
      \Athree+\Aone & 3 & \onerel & \text{\cite{MR2520770}, \cite{arXiv:1304.3352}} \\
      \Afour & 3 & \onerel & \text{\cite{MR2543667}, \cite{arXiv:1304.3352}} \\
      \Dfour & 2 & \onerel & \text{\cite{MR2290499}, \cite{arXiv:1304.3352}} \\
      \Athree + 2\Aone & 2 & \toric & \text{\cite{MR1620682}} \\
      \Dfive & 1 & \onerel & \text{\cite{MR1906155}, \cite{MR2320172}} \\
      \hline
    \end{array}
  \end{equation*}
  \smallskip
  \caption{Del Pezzo surfaces of degree 4}
  \label{tab:degree_4}
\end{table}

It turns out that types $4\Aone$, $\Atwo+2\Aone$ and $\Athree+2\Aone$ are
toric, with cyclic Dynkin diagrams of Cox ring generators
\begin{equation*}
  (1,2,1,2,1,2,1,2),\qquad (2,1,2,1,1,2,1,2),\qquad(0,2,1,2,2,2,1,2).
\end{equation*}
Ordinary quartic del Pezzo surfaces and types $\Aone$, $2\Aone$ with eight or
nine lines and $\Atwo$ have Cox rings with more than one relation by
Theorem~\ref{thm:generators} because their number of negative curves is larger
than $13-d=9$.

For type $\Athree$ with four lines, this is not as
obvious. We have only seven negative curves in the following
configuration:
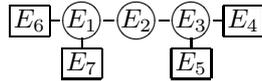
\begin{equation*}
  \dynkin{\li 6 \ar@{-}[r] & \ex{1} \ar@{-}[r] & \ex{2} \ar@{-}[r] & \ex{3}
    \ar@{-}[d]\ar@{-}[r] & \li 4\\
    & \li 7 \ar@{-}[u] & & \li 5}
\end{equation*}
Step~\ref{step:generator_degrees} shows that $\Cox(\tS)$ must have at least
four additional nef generators, one in each of the degrees $\cl E_1+\cl
E_2+\cl E_3+\cl E_i+\cl E_j$ with $i=4,5$ and $j=6,7$, whose space of global
sections is $2$-dimensional, but only a $1$-dimensional subspace is generated
by sections $\e_1, \dots, \e_7$ corresponding to $E_1, \dots, E_7$.

For the remaining seven types, the generalized quartic del Pezzo surface is
unique up to isomorphism because it can be obtained by contraction of a
$(-1)$-curve from a cubic surface that is unique up to isomorphism by
\cite{MR80f:14021}. We compute the Cox rings as follows.

\coxcase{v}{$3\Aone$}
\coxembed{\Pfour}{x_0x_2-x_3x_4 = x_0x_1+x_0x_2+x_1x_2 = 0}
\coxbiratproj 4 {x_1:x_2:x_4}{E_7}
\coxtwomulti{\Aone, \Aone, \Aone}{(0:1:0:0:0), (0:0:0:1:0), (0:0:0:0:1)}{$E_1, E_2, E_3$}{
  \coximo{1}{\ell_0-\ell_1-\ell_2-\ell_3}{(0:1:0:0:0)}{-y_0-y_1}
  \coximo{2}{\ell_0-\ell_3-\ell_4-\ell_5}{(0:0:0:1:0)}{y_2}
  \coximopd{3}{\ell_1-\ell_2}{(0:0:0:0:1)}{(0:0:1)} 
}
\coxonem{E_4, \dots, E_9}{ 
  \coxim{4}{\ell_0-\ell_1-\ell_4}{x_0=x_1=x_3=0}{y_0}
  \coximp{5}{\ell_4}{x_0=x_1=x_4=0}{(0:1:0)} 
  \coximp{6}{\ell_2}{x_0=x_2=x_3=0}{(0:0:1)}
  \coximp{7}{\ell_3}{x_0=x_2=x_4=0}{(1:-1:0)}
  \coxim{8}{\ell_0-\ell_1-\ell_5}{x_1=x_2=x_3=0}{y_1}
  \coximpd{9}{\ell_5}{x_1=x_2=x_4=0}{(1:0:0)} }
\coxdynkin{&\li{6} \ar@{-}[r]\ar@{-}[dl] & \ex{1} \ar@{-}[r]& \li{7}\ar@{-}[dr]\\
  \ex{3}\ar@{-}[r]\ar@{-}[dr] & \li{8}\ar@{-}[rr] & & \li{9}\ar@{-}[r] & \ex{2}\\
  & \li{4}\ar@{-}[rr] & & \li{5}\ar@{-}[ur]}
\coxonerel{\e_1, \dots, \e_9}{\e_4\e_5 + \e_1\e_6\e_7 +
  \e_8\e_9}{\ell_0-\ell_1}
\coxantie{(\e_1\e_2\e_3\e_4\e_5\e_6\e_7,
   \e_2\e_3\e_4\e_5\e_8\e_9, \e_1\e_2\e_3\e_6\e_7\e_8\e_9,
  \e_1\e_3^2\e_4\e_6^2\e_8,
  \e_1\e_2^2\e_5\e_7^2\e_9)}
{\pi^*(-x_0-x_1) =\e_2\e_3\e_4^2\e_5^2,
  \pi^*(-x_0-x_2) = \e_1^2\e_2\e_3\e_6^2\e_7^2,
  \pi^*(-x_1-x_2) = \e_2\e_3\e_8^2\e_9^2}


\coxcase{vi}{$\Atwo+\Aone$}
\coxembed{\Pfour}{x_0x_1-x_2x_3 = x_1x_2+x_2x_4+x_3x_4 = 0}
\coxbiratproj 4 {x_1:x_2:x_3}{E_9}
\coxtwotwo{\Atwo}{(1:0:0:0:0)}{E_1, E_2}{\Aone}{(0:0:0:0:1)}{E_3}{
  \coximop{1}{\ell_1-\ell_2}{(1:0:0:0:0)}{(0:1:-1)}
  \coximo{2}{\ell_0-\ell_1-\ell_4-\ell_5}{(1:0:0:0:0)}{y_0}
  \coximod{3}{\ell_0-\ell_1-\ell_2-\ell_3}{(0:0:0:0:1)}{-y_1-y_2}
}
\coxonem{E_4, \dots, E_9}{
  \coximp{4}{\ell_3}{x_0=x_2=x_3=0}{(1:0:0)}
  \coxim{5}{\ell_0-\ell_3-\ell_5}{x_0=x_2=x_4=0}{y_1}
  \coximp{6}{\ell_4}{x_1=x_3=x_4=0}{(0:1:0)}
  \coximp{7}{\ell_5}{x_1=x_2=x_4=0}{(0:0:1)}
  \coxim{8}{\ell_0-\ell_3-\ell_4}{x_0=x_3=x_1+x_4=0}{y_2}
  \coximpd{9}{\ell_2}{x_1=x_2=x_3=0}{(0:1:-1)}
}
\coxdynkin{& \ex{1} \ar@{-}[r]\ar@{-}[dl] & \li{9} \ar@{-}[r]& \ex{3}\ar@{-}[dr]\\
  \ex{2}\ar@{-}[r]\ar@{-}[dr] & \li{6}\ar@{-}[rr] & & \li{8}\ar@{-}[r] & \li{4}\\
  & \li{7}\ar@{-}[rr] & & \li{5}\ar@{-}[ur]}
\coxonerel{\e_1, \dots, \e_9}{\e_5\e_7 + \e_1\e_3\e_9^2 + \e_6\e_8}{\ell_0-\ell_3}
\coxantie{(\e_3\e_4^2\e_5\e_8,
  \e_1^2\e_2^2\e_3\e_6\e_7\e_9^2,
  \e_1\e_2\e_3\e_4\e_5\e_7\e_9,\e_1\e_2\e_3\e_4\e_6\e_8\e_9,
  \e_1\e_2^2\e_5\e_6\e_7^2)}{\pi^*(-x_1-x_4) =
  \e_1\e_2^2\e_6^2\e_7\e_8,\ \pi^*(-x_2-x_3) =
  \e_1^2\e_2\e_3^2\e_4\e_9^3}


\coxcase{vii}{$\Athree$}
\coxembed{\Pfour}{x_0x_1-x_2x_3 = x_2x_4+x_0x_3+x_1x_3 = 0}
\coxbiratproj 4 {x_0:x_1:x_2}{E_5}
\coxtwoone{\Athree}{(0:0:0:0:1)}{E_1, E_2, E_3}{
  \coximop{1}{\ell_1-\ell_4}{(0:0:0:0:1)}{(0:1:0)}
  \coximo{2}{\ell_0-\ell_1-\ell_2-\ell_3}{(0:0:0:0:1)}{y_2}
  \coximopd{3}{\ell_2-\ell_5}{(0:0:0:0:1)}{(1:0:0)}
}
\coxonem{E_4, \dots, E_8}{
  \coximp{4}{\ell_4}{x_0=x_2=x_3=0}{(0:1:0)}
  \coximp{5}{\ell_3}{x_0=x_1=x_2=0}{(1:-1:0)}
  \coximp{6}{\ell_5}{x_1=x_2=x_3=0}{(1:0:0)}
  \coxim{7}{\ell_0-\ell_1-\ell_4}{x_0=x_3=x_4=0}{y_0}
  \coximd{8}{\ell_0-\ell_2-\ell_5}{x_1=x_3=x_4=0}{y_1}
}
\coxother{$(0)$-curve $E_9$}{
  \coximd{9}{\ell_0-\ell_3}{x_4=x_0+x_1=x_0^2+x_2x_3}{-y_0-y_1}
}
\coxdynkinp{E_7 \cap E_8 \cap E_9}{(0:0:1:0:0)}{\li{7} \ar@{-}[r] \ar@{-}[dd] \ar@{-}[dr] & \li{4} \ar@{-}[r] & \ex{1} \ar@{-}[dr]\\
  & E_9 \ar@{-}[r] & \li{5} \ar@{-}[r] & \ex{2}\\
  \li{8} \ar@{-}[r] \ar@{-}[ur] & \li{6} \ar@{-}[r] & \ex{3} \ar@{-}[ur]}
\coxonerel{\e_1, \dots, \e_9}{\e_5\e_9 + \e_1\e_4^2\e_7 + \e_3\e_6^2\e_8}{\ell_0}
\coxantie{(\e_1^2\e_2^2\e_3\e_4^2\e_5\e_7,
  \e_1\e_2^2\e_3^2\e_5\e_6^2\e_8, \e_1^2\e_2^3\e_3^2\e_4\e_5^2\e_6,
  \e_1\e_2\e_3\e_4\e_6\e_7\e_8,\e_7\e_8\e_9)}{\pi^*(-x_0-x_1) =
  \e_1\e_2^2\e_3\e_5^2\e_9}


\coxcase{xi}{$\Athree+\Aone$}
\coxembed{\Pfour}{x_0x_3-x_2x_4 = x_0x_1+x_1x_3+x_2^2 = 0}
\coxbiratproj 4 {x_0:x_2:x_3}{E_5}
\coxtwotwo{\Athree}{(0:0:0:0:1)}{E_2, E_3, E_4}{\Aone}{(0:1:0:0:0)}{E_1}{
  \coximo{1}{\ell_0-\ell_1-\ell_2-\ell_3}{(0:1:0:0:0)}{-y_0-y_2}
  \coximo{2}{\ell_0-\ell_1-\ell_4-\ell_5}{(0:0:0:0:1)}{y_1}
  \coximop{3}{\ell_1-\ell_2}{(0:0:0:0:1)}{(1:0:-1)}
  \coximopd{4}{\ell_2-\ell_3}{(0:0:0:0:1)}{(1:0:-1)}
}
\coxonem{E_5, E_6, E_7}{
  \coximp{5}{\ell_3}{x_0=x_2=x_3=0}{(1:0:-1)}
  \coximp{6}{\ell_4}{x_0=x_1=x_2=0}{(0:0:1)}
  \coximpd{7}{\ell_5}{x_1=x_2=x_3=0}{(1:0:0)}
}
\coxother{$(0)$-curves $E_8, E_9$}{
  \coxim{8}{\ell_0-\ell_5}{x_3=x_4=x_0x_1+x_2^2=0}{y_2}
  \coximd{9}{\ell_0-\ell_4}{x_0=x_4=x_1x_3+x_2^2=0}{y_0}
}
\coxdynkinp{E_1 \cap E_8 \cap E_9}{(0:1:0:0:0)}{E_8 \ar@{-}[dr] \ar@{-}[dd] \ar@{-}[rrrr]& & & & \li{7} \ar@{-}[dr]\\
  & \ex{1} \ar@{-}[r] & \li{5} \ar@{-}[r] & \ex{4} \ar@{-}[r] & \ex{3} \ar@{-}[r] & \ex{2}\\
  E_9 \ar@{-}[ur] \ar@{-}[rrrr] & & & & \li{6} \ar@{-}[ur]}
\coxonerel{\e_1, \dots, \e_9}{\e_6\e_9+\e_7\e_8+\e_1\e_3\e_4^2\e_5^3}{\ell_0}
\coxantie{(\e_1\e_2\e_3\e_4\e_5\e_6\e_9,
  \e_2^3\e_3^2\e_4\e_6^2\e_7^2,
  \e_1\e_2^2\e_3^2\e_4^2\e_5^2\e_6\e_7,\e_1\e_2\e_3\e_4\e_5\e_7\e_8,
  \e_1\e_8\e_9)}{\pi^*(-x_0-x_3) = \e_1^2\e_2\e_3^2\e_4^3\e_5^4}


\coxcase{xii}{$\Afour$}
\coxembed{\Pfour}{x_0x_1-x_2x_3 = x_0x_4+x_1x_2+x_3^2 = 0}
\coxbiratproj 4 {x_0:x_2:x_3}{E_5}
\coxtwoone{\Afour}{(0:0:0:0:1)}{E_1, \dots, E_4}{
  \coximop{1}{\ell_3-\ell_4}{(0:0:0:0:1)}{(0:0:1)}
  \coximop{2}{\ell_4-\ell_5}{(0:0:0:0:1)}{(0:0:1)}
  \coximo{3}{\ell_0-\ell_1-\ell_3-\ell_4}{(0:0:0:0:1)}{y_0}
  \coximopd{4}{\ell_1-\ell_2}{(0:0:0:0:1)}{(0:1:0)}
}
\coxonem{E_5, E_6, E_7}{
  \coximp{5}{\ell_5}{x_0=x_2=x_3=0}{(0:0:1)}
  \coximp{6}{\ell_2}{x_0=x_1=x_3=0}{(0:1:0)}
  \coximd{7}{\ell_0-\ell_1-\ell_2}{x_1=x_3=x_4=0}{y_2}
}
\coxother{$(1)$-curve $E_8$ and $(0)$-curve $E_9$}{
  \coxims{8}{2\ell_0-\ell_3-\ell_4-\ell_5}{x_4=x_0x_1-x_2x_3=x_1x_2+x_3^2=x_0x_3+x_2^2=0}{-y_0y_2-y_1^2}
  \coximd{9}{\ell_0-\ell_3}{x_1=x_2=x_0x_4+x_3^2=0}{y_1} 
}
\coxdynkinp{E_7 \cap E_8 \cap E_9}{(1:0:0:0:0)}{E_8 \ar@{-}[rrrr] \ar@{-}[dd] \ar@{=}[dr] & & & & \li{5} \ar@{-}[dr]\\
  & \li{7} \ar@{-}[r] & \li{6} \ar@{-}[r] & \ex{4} \ar@{-}[r] & \ex{3} \ar@{-}[r] & \ex{2} \ar@{-}[dl]\\
  E_9 \ar@{-}[rrrr] \ar@{-}[ur] & & & & \ex{1}}
\coxonerel{\e_1, \dots, \e_9}{\e_5\e_8 +
  \e_1\e_9^2 + \e_3\e_4^2\e_6^3\e_7}{2\ell_0-\ell_3-\ell_4}
\coxanti{(\e_1^2\e_2^4\e_3^3\e_4^2\e_5^3\e_6,
  \e_1\e_2\e_3\e_4\e_6\e_7\e_9, \e_1^2\e_2^3\e_3^2\e_4\e_5^2\e_9, \e_1\e_2^2\e_3^2\e_4^2\e_5\e_6^2\e_7, \e_7\e_8)}


\coxcase{xiii}{$\Dfour$}
\coxembed{\Pfour}{x_0x_3-x_1x_4 = x_0x_1+x_1x_3+x_2^2 = 0}
\coxbiratproj 4 {x_0:x_1:x_2}{E_5}
\coxtwoone{\Dfour}{(0:0:0:0:1)}{E_1, \dots, E_4}{
  \coximop{1}{\ell_2-\ell_3}{(0:0:0:0:1)}{(1:0:0)}
  \coximop{2}{\ell_1-\ell_2}{(0:0:0:0:1)}{(1:0:0)}
  \coximo{3}{\ell_0-\ell_1-\ell_2-\ell_5}{(0:0:0:0:1)}{y_1}
  \coximopd{4}{\ell_3-\ell_4}{(0:0:0:0:1)}{(1:0:0)}
}
\coxonem{E_5, E_6}{
  \coximp{5}{\ell_5}{x_0=x_1=x_2=0}{(0:0:1)}
  \coximpd{6}{\ell_4}{x_1=x_2=x_3=0}{(1:0:0)}
}
\coxother{$(0)$-curves $E_7, E_8, E_9$}{
  \coxim{7}{\ell_0-\ell_1}{x_2=x_0+x_3=x_0^2+x_1x_4=0}{y_2}
  \coxim{8}{\ell_0-\ell_5}{x_0=x_4=x_1x_3+x_2^2=0}{y_0}
  \coximds{9}{2\ell_0-\ell_1-\ell_2-\ell_3-\ell_4}{x_3=x_4=x_0x_1+x_2^2=0}{-y_0y_1-y_2^2}
}
\coxdynkinp{E_7 \cap E_8 \cap E_9}{(0:1:0:0:0)}{E_8  \ar@{-}[r]\ar@{-}[dr] \ar@{=}[dd]& \li{5} \ar@{-}[r] & \ex{3} \ar@{-}[dr]\\
  & E_7 \ar@{-}[r] & \ex{2} \ar@{-}[r] & \ex{1} \\
  E_9 \ar@{-}[r] \ar@{-}[ur] & \li{6} \ar@{-}[r] & \ex{4} \ar@{-}[ur]}
\coxonerel{\e_1, \dots, \e_9}{\e_3\e_5^2\e_8 + \e_4\e_6^2\e_9 +
  \e_2\e_7^2}{2\ell_0-\ell_1-\ell_2}
\coxantie{(\e_1^2\e_2\e_3^2\e_4\e_5^2\e_8,
  \e_1^4\e_2^2\e_3^3\e_4^3\e_5^2\e_6^2,
  \e_1^3\e_2^2\e_3^2\e_4^2\e_5\e_6\e_7,
  \e_1^2\e_2\e_3\e_4^2\e_6^2\e_9, \e_8\e_9)}{\pi^*(-x_0-x_3) =
  \e_1^2\e_2^2\e_3\e_4\e_7^2}


\coxcase{xv}{$\Dfive$}
\coxembed{\Pfour}{x_0x_1-x_2^2 = x_3^2+x_0x_4+x_1x_2 = 0}
\coxbiratproj 4 {x_0:x_2:x_3}{E_6}
\coxtwoone{\Dfive}{(0:0:0:0:1)}{E_1, \dots, E_5}{
  \coximop{1}{\ell_3-\ell_4}{(0:0:0:0:1)}{(0:0:1)}
  \coximop{2}{\ell_4-\ell_5}{(0:0:0:0:1)}{(0:0:1)}
  \coximo{3}{\ell_0-\ell_1-\ell_2-\ell_3}{(0:0:0:0:1)}{y_0}
  \coximop{4}{\ell_2-\ell_3}{(0:0:0:0:1)}{(0:0:1)}
  \coximopd{5}{\ell_1-\ell_2}{(0:0:0:0:1)}{(0:0:1)}
}
\coxonem{E_6}{
  \coximpd{6}{\ell_5}{x_0=x_2=x_3 = 0}{(0:0:1)}
}
\coxother{$(1)$-curve $E_7$, $(0)$-curve $E_8$ and $(4)$-curve $E_9$}{
  \coxims{7}{\ell_0}{x_3=x_0x_1-x_2^2=x_0x_4+x_1x_2=x_1^2+x_2x_4=0}{y_2}
  \coxim{8}{\ell_0-\ell_1}{x_1=x_2=x_0x_4+x_3^2=0}{y_1}
  \coximds{9}{3\ell_0-\ell_1-\ell_2-\ell_3-\ell_4-\ell_5}{x_4=x_0x_1-x_2^2=x_3^2+x_1x_2=0}{-y_0y_2^2-y_1^3} 
}
\coxdynkinp{E_7 \cap E_8 \cap E_9}{(1:0:0:0:0)}{E_8 \ar@{-}[r]\ar@{=}[dd] \ar@{-}[dr]  & \ex{5} \ar@{-}[r] & \ex{4}\ar@{-}[dr] \\
      & E_7\ar@{-}[r]  & \ex{3} \ar@{-}[r] & \ex{1}\\
      E_9 \ar@{-}[r]\ar@3{-}[ur] & \li{6} \ar@{-}[r] & \ex{2}\ar@{-}[ur]}
\coxonerel{\e_1, \dots, \e_9}{\e_3\e_7^2 +
    \e_2\e_6^2\e_9 + \e_4\e_5^2\e_8^3}{3\ell_0-\ell_1-\ell_2-\ell_3}
  \coxanti{(&\e_1^6\e_2^5\e_3^3\e_4^4\e_5^2\e_6^4,
    \e_1^2\e_2\e_3\e_4^2\e_5^2\e_8^2,
    \e_1^4\e_2^3\e_3^2\e_4^3\e_5^2\e_6^2\e_8,\e_1^3\e_2^2\e_3^2\e_4^2\e_5\e_6\e_7, \e_9)}

\subsection{Degree 3}\label{sec:degree_3}

The classification of cubic del Pezzo surfaces is classical and goes
back to Schl\"afli \cite{schlaefli}. Together with their numbers of lines, the
list in Table~\ref{tab:degree_3} can be found in \cite{MR80f:14021}.

\begin{table}[ht]
  \centering
  \begin{equation*}
    \begin{array}{cccc}      \hline
      \text{singularities} 
      & \text{lines} & \text{type} & \text{reference}\\
      \hline
      - & 27 & \tworel & \\
      \Aone & 21 & \tworel & \\
      2\Aone & 16 & \tworel & \\
      \Atwo & 15 & \tworel & \\
      3\Aone & 12 & \tworel & \\
      \Atwo+\Aone & 11 & \tworel & \\
      \Athree & 10 & \tworel & \\
      4\Aone & 9 & \tworel & \text{\cite{MR2075628}} \\
      \Atwo + 2\Aone & 8 & \tworel & \\
      \Athree + \Aone & 7 & \tworel & \\
      2\Atwo & 7 & \tworel & \\
      \Afour & 6 & \tworel & \\
      \Dfour & 6 & \onerel & \text{\cite{MR2250046}, \cite{arXiv:1207.2685}} \\
      \Athree + 2\Aone & 5 & \onerel & \\
      2\Atwo+\Aone & 5 & \onerel & \text{\cite{MR2990624}}\\
      \Afour+\Aone & 4 & \onerel & \\
      \Afive & 3 & \tworel & \\
      \Dfive & 3 & \onerel & \text{\cite{MR2520769}}\\
      3\Atwo & 3 & \toric & \text{\cite{MR1620682}, \dots, \cite{arXiv:1204.0383}}\\
      \Afive+\Aone & 2 & \onerel & \text{\cite{arXiv:1205.0373}}\\
      \Esix & 1 & \onerel & \text{\cite{MR2332351}}\\
      \hline
    \end{array}
  \end{equation*}
  \smallskip
  \caption{Cubic surfaces}
  \label{tab:degree_3}
\end{table}

The only toric type is $3\Atwo$, with cyclic Dynkin diagram of Cox ring generators
\begin{equation*}
  (-2,-2,-1,-2,-2,-1,-2,-2,-1).
\end{equation*}
Ordinary cubic surfaces and types $\Aone$, $2\Aone$, $\Atwo$, $3\Aone$,
$\Atwo+\Aone$, $\Athree$, $4\Aone$, $\Atwo+2\Aone$, $\Athree+\Aone$ and
$2\Atwo$ have two or more relations in their Cox rings because their number of
negative curves is more than $13-d=10$.

For type $\Afour$, the extended Dynkin diagram of negative curves is:
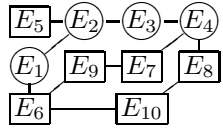
\begin{equation*}
  \dynkin{\li 5 \ar@{-}[r]& \ex{2} \ar@{-}[r] & \ex{3} \ar@{-}[r] & \ex{4}\\
    \ex{1}\ar@{-}[ur] & \li 9\ar@{-}[r] & \li 7\ar@{-}[ur] & \li 8\ar@{-}[u]\\
    \li 6\ar@{-}[u]\ar@{-}[ur]\ar@{-}[rr] & & \li{10} \ar@{-}[ur]}
\end{equation*}
Contracting $E_6$ gives type $\Athree$ with four lines of
degree $4$, whose Cox ring has $\ge 2$ relations, so the same is true
for type $\Afour$ of degree $3$ by Proposition~\ref{prop:more_generators}.

For type $\Afive$, the extended Dynkin diagram of negative curves is:
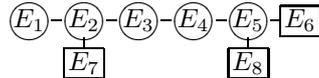
\begin{equation*}
  \dynkin{\ex{1} \ar@{-}[r] & \ex{2} \ar@{-}[r] \ar@{-}[d] & \ex{3}
    \ar@{-}[r] & \ex{4} \ar@{-}[r] & \ex{5} \ar@{-}[r] \ar@{-}[d] & \li 6\\
    & \li 7 & & & \li 8}
\end{equation*}
Step~\ref{step:generator_degrees} shows that we have at least three additional
generators, namely in degrees $\cl E_2+\cl E_3+\cl E_4+\cl E_5+\cl E_7+\cl
E_i$ for $i=6,8$ and in degree $\cl E_1+2\cl E_2+2\cl E_3+2\cl E_4+2\cl
E_5+2\cl E_6+\cl E_7$.

By \cite{MR80f:14021}, the type determines the surface up to isomorphism
except for ordinary cubic surfaces and types $\Aone, 2\Aone, \Atwo, 3\Aone,
\Atwo+\Aone, \Athree, 2\Atwo$ (where we have infinite families with one or
more parameters; the Cox rings have more than one relation) and type $\Dfour$
(with precisely two isomorphy classes; see also \cite[Remark~4.1]{MR2029868}).

\coxcase{xii}{$\Dfour$, first isomorphy class}
\coxembed{\Pthree}{x_0(x_1+x_2+x_3)^2-x_1x_2x_3 = 0}
\coxbiratproj 3 {x_1:x_2:x_3}{(1:0:0:0)}
\coxtwoone{\Dfour}{(1:0:0:0)}{E_1, \dots, E_4}{
  \coximo{1}{\ell_0-\ell_1-\ell_2-\ell_3}{(1:0:0:0)}{y_0+y_1+y_2}
  \coximop{2}{\ell_1-\ell_4}{(1:0:0:0)}{(0:1:-1)}
  \coximop{3}{\ell_2-\ell_5}{(1:0:0:0)}{(1:0:-1)}
  \coximopd{4}{\ell_3-\ell_6}{(1:0:0:0)}{(1:-1:0)}
}
\coxonem{E_5, \dots, E_{10}}{
  \coximp{5}{\ell_4}{x_1=x_2+x_3 = 0}{(0:1:-1)}
  \coximp{6}{\ell_5}{x_2=x_1+x_3 = 0}{(1:0:-1)}
  \coximp{7}{\ell_6}{x_3 = x_1+x_2 = 0}{(1:-1:0)}
  \coxim{8}{\ell_0-\ell_1-\ell_4}{x_0 = x_1 = 0}{y_0}
  \coxim{9}{\ell_0-\ell_2-\ell_5}{x_0 = x_2 = 0}{y_1}
  \coximd{10}{\ell_0-\ell_3-\ell_6}{x_0 = x_3 = 0}{y_2}
}
\coxdynkinpp{$E_8 \cap E_9 \cap E_{10}=\emptyset$}{\li{8} \ar@{-}[rr] \ar@{-}[dr] \ar@{-}[dd] & & \li{5} \ar@{-}[r]& \ex{2} \ar@{-}[dr]\\
& \li{9} \ar@{-}[r] & \li{6} \ar@{-}[r] & \ex{3} \ar@{-}[r] & \ex{1}\\
\li{{10}} \ar@{-}[rr] \ar@{-}[ur] & & \li{7} \ar@{-}[r] & \ex{4} \ar@{-}[ur]}
\coxonerel{\e_1, \dots, \e_{10}}{\e_2\e_5^2\e_8 +
  \e_3\e_6^2\e_9 + \e_4\e_7^2\e_{10} -
  \e_1\e_2\e_3\e_4\e_5\e_6\e_7}{\ell_0}
\coxantie{(\e_8\e_9\e_{10},
  \e_1^2\e_2^2\e_3\e_4\e_5^2\e_8, \e_1^2\e_2\e_3^2\e_4\e_6^2\e_9,
  \e_1^2\e_2\e_3\e_4^2\e_7^2\e_{10})}{\pi^*(x_1+x_2+x_3) =
  \e_1^3\e_2^2\e_3^2\e_4^2\e_5\e_6\e_7}
\coxcase{xii}{$\Dfour$, second isomorphy class}
It differs from the first one as follows: Anticanonical model $S \subset \Pthree$:
\begin{equation*}
x_0(x_1+x_2+x_3)^2 + x_1x_2(x_1+x_2)=0.
\end{equation*}
The $(-1)$-curve
$E_{10}$ meets $E_8, E_9$ in a point $p$ with $\pi(p)=(0:0:0:1)$:
\begin{align*}
  \coximd{10}{\ell_0-\ell_3-\ell_6}{x_0 = x_1+x_2 = 0}{-y_0-y_1}
\end{align*}
Cox ring: relation
\begin{equation*}
  \e_2\e_5^2\e_8 + \e_3\e_6^2\e_9 + \e_4\e_7^2\e_{10}=0.
\end{equation*}
\coxantie{(\e_8\e_9\e_{10},
  \e_1^2\e_2^2\e_3\e_4\e_5^2\e_8, \e_1^2\e_2\e_3^2\e_4\e_6^2\e_9,
  \e_1^3\e_2^2\e_3^2\e_4^2\e_5\e_6\e_7 +
  \e_1^2\e_2\e_3\e_4^2\e_7^2\e_{10})}{\pi^*(-x_1-x_2) =
  \e_1^2\e_2\e_3\e_4^2\e_7^2\e_{10},\ \pi^*(x_1+x_2+x_3) =
  \e_1^3\e_2^2\e_3^2\e_4^2\e_5\e_6\e_7}


\coxcase{xiii}{$\Athree+2\Aone$}
\coxembed{\Pthree}{x_3^2(x_1+x_2)+x_0x_1x_2 = 0}
\coxbiratproj 3 {x_1:x_2:x_3}{(1:0:0:0)}
\coxtwomulti{\Athree,\Aone,\Aone}{(1:0:0:0), (0:1:0:0), (0:0:1:0)}{$E_2,E_3,E_4$
  and $E_1$ and $E_5$}{
  \coximop{1}{\ell_3-\ell_4}{(0:1:0:0)}{(1:0:0)}
  \coximo{2}{\ell_0-\ell_3-\ell_4-\ell_5}{(1:0:0:0)}{y_1}
  \coximop{3}{\ell_5-\ell_6}{(1:0:0:0)}{(0:0:1)}
  \coximo{4}{\ell_0-\ell_1-\ell_2-\ell_5}{(1:0:0:0)}{y_0}
  \coximopd{5}{\ell_1-\ell_2}{(0:0:1:0)}{(0:1:0)}}
\coxonem{E_6, \dots, E_{10}}{
  \coximp{6}{\ell_2}{x_1=x_3=0}{(0:1:0)}
  \coximp{7}{\ell_4}{x_2=x_3=0}{(1:0:0)}
  \coximp{8}{\ell_6}{x_1=x_2=0}{(0:0:1)}
  \coxim{9}{\ell_0-\ell_5-\ell_6}{x_0=x_1+x_2=0}{-y_0-y_1}
  \coximd{10}{\ell_0-\ell_1-\ell_3}{x_0=x_3=0}{y_2}
}
\coxdynkin{\ex{2} \ar@{-}[r]\ar@{-}[d] & \li{7} \ar@{-}[r]& \ex{1}\ar@{-}[dr]\\
  \ex{3}\ar@{-}[r]\ar@{-}[d] & \li{8}\ar@{-}[r] & \li{9}\ar@{-}[r] & \li{10}\\
  \ex{4}\ar@{-}[r] & \li{6}\ar@{-}[r] & \ex{5}\ar@{-}[ur]}
\coxonerel{\e_1, \dots, \e_{10}}{\e_4\e_5\e_6^2 + \e_1\e_2\e_7^2 + \e_8\e_9}{\ell_0-\ell_5}
\coxantie{(\e_1\e_5\e_9\e_{10}^2,
  \e_2\e_3^2\e_4^2\e_5\e_6^2\e_8,
  \e_1\e_2^2\e_3^2\e_4\e_7^2\e_8,
  \e_1\e_2\e_3\e_4\e_5\e_6\e_7\e_{10})}{\pi^*(-x_1-x_2) =
  \e_2\e_3^2\e_4\e_8^2\e_9}


\coxcase{xiv}{$2\Atwo+\Aone$}
\coxembed{\Pthree}{x_3^2(x_1+x_3)+x_0x_1x_2 = 0}
\coxbiratproj 3 {x_0:x_2:x_3}{(0:1:0:0)}
\coxtwomulti{\Atwo,\Atwo,\Aone}{(1:0:0:0),(0:0:1:0),(0:1:0:0)}{$E_2,E_3$
  and $E_4,E_5$ and $E_1$}{
  \coximo{1}{2\ell_0-\ell_1-\ell_2-\ell_3-\ell_4-\ell_5-\ell_6}{(0:1:0:0)}{-y_0y_1-y_2^2}
  \coximop{2}{\ell_2-\ell_3}{(1:0:0:0)}{(1:0:0)}
  \coximop{3}{\ell_1-\ell_2}{(1:0:0:0)}{(1:0:0)}
  \coximop{4}{\ell_5-\ell_6}{(0:0:1:0)}{(0:1:0)}
  \coximopd{5}{\ell_4-\ell_5}{(0:0:1:0)}{(0:1:0)} }
\coxonem{E_6, \dots, E_{10}}{
  \coximp{6}{\ell_6}{x_0=x_3=0}{(0:1:0)}
  \coxim{7}{\ell_0-\ell_1-\ell_4}{x_1=x_3=0}{y_2}
  \coximp{8}{\ell_3}{x_2=x_3=0}{(1:0:0)}
  \coxim{9}{\ell_0-\ell_4-\ell_5,}{x_0=x_1+x_3=0}{y_0}
  \coximd{10}{\ell_0-\ell_1-\ell_2}{x_2=x_1+x_3=0}{y_1}
}
\coxdynkin{ & \li{8} \ar@{-}[r] & \ex{1} \ar@{-}[r] & \li{6}\ar@{-}[dr]  \\
      \ex{2} \ar@{-}[ur]\ar@{-}[r]\ar@{-}[dr]   & \ex{3}\ar@{-}[r]  & \li{7} \ar@{-}[r] & \ex{5} \ar@{-}[r] & \ex{4}\\
      & \li{{10}} \ar@{-}[rr]& & \li{9}\ar@{-}[ur]}
\coxonerel{\e_1, \dots, \e_{10}}{\e_3\e_5\e_7^2 + \e_1\e_6\e_8 + \e_9\e_{10}}{2\ell_0-\ell_1-\ell_2-\ell_4-\ell_5}
\coxantie{(\e_1\e_4^2\e_5\e_6^2\e_9, \e_2\e_3^2\e_4\e_5^2\e_7^3, \e_1\e_2^2\e_3\e_8^2\e_{10}, \e_1\e_2\e_3\e_4\e_5\e_6\e_7\e_8)}{\pi^*(-x_1-x_3) = \e_2\e_3\e_4\e_5\e_7\e_9\e_{10}}


\coxcase{xv}{$\Afour+\Aone$}
\coxembed{\Pthree}{x_2x_3^2+x_1^2x_3+x_0x_1x_2 = 0}
\coxbiratproj 3 {x_1:x_2:x_3}{(1:0:0:0)}
\coxtwotwo{\Afour}{(1:0:0:0)}{E_1, \dots, E_4}{\Aone}{(0:0:1:0)}{E_5}{
  \coximo{1}{\ell_0-\ell_1-\ell_2-\ell_6}{(1:0:0:0)}{y_1}
  \coximop{2}{\ell_2-\ell_3}{(1:0:0:0)}{(0:0:1)}
  \coximop{3}{\ell_1-\ell_2}{(1:0:0:0)}{(0:0:1)}
  \coximo{4}{\ell_0-\ell_1-\ell_4-\ell_5}{(1:0:0:0)}{y_0}
  \coximopd{5}{\ell_4-\ell_5}{(0:0:1:0)}{(0:1:0)}
}
\coxonem{E_6, \dots, E_9}{
  \coximp{6}{\ell_5}{x_1=x_3=0}{(0:1:0)}
  \coximp{7}{\ell_3}{x_1=x_2=0}{(0:0:1)}
  \coximp{8}{\ell_6}{x_2=x_3=0}{(1:0:0)}
  \coximd{9}{\ell_0-\ell_4-\ell_6}{x_0=x_3=0}{y_2}
}
\coxother{$(0)$-curves $E_{10}$}{
  \coximds{10}{2\ell_0-\ell_1-\ell_2-\ell_3-\ell_4}{x_0=x_1^2+x_2x_3=0}{-y_0^2-y_1y_2}
}
\coxdynkinp{E_5 \cap E_9 \cap E_{10}}{(0:0:1:0)}{\ex{5} \ar@{-}[r]\ar@{-}[dr]\ar@{-}[dd]    & \li{6} \ar@{-}[r] & \ex{4} \ar@{-}[r] & \ex{3}\ar@{-}[dr]  \\
  & \li{9}\ar@{-}[r] & \li{8} \ar@{-}[r] & \ex{1} \ar@{-}[r] & \ex{2}\\
  E_{10} \ar@{-}[rrr] \ar@{-}[ur] & & & \li{7}\ar@{-}[ur]}
\coxonerel{\e_1, \dots, \e_{10}}{\e_1\e_8^2\e_9 + \e_3\e_4^2\e_5\e_6^3 + \e_7\e_{10}}{2\ell_0-\ell_1-\ell_2-\ell_4}
\coxanti{
    (\e_5\e_9\e_{10}, \e_1\e_2^2\e_3^2\e_4^2\e_5\e_6^2\e_7,
    \e_1^2\e_2^3\e_3^2\e_4\e_7^2\e_8, \e_1\e_2\e_3\e_4\e_5\e_6\e_8\e_9)}


\coxcase{xvii}{$\Dfive$}
\coxembed{\Pthree}{x_3x_0^2+x_0x_2^2+x_1^2x_2 = 0}
\coxbiratproj 3 {x_0:x_1:x_2}{(0:0:0:1)}
\coxtwoone{\Dfive}{(0:0:0:1)}{E_1, \dots, E_5}{
  \coximop{1}{\ell_2-\ell_3}{(0:0:0:1)}{(0:0:1)}
  \coximop{2}{\ell_3-\ell_4}{(0:0:0:1)}{(0:0:1)}
  \coximop{3}{\ell_1-\ell_2}{(0:0:0:1)}{(0:0:1)}
  \coximo{4}{\ell_0-\ell_1-\ell_2-\ell_5}{(0:0:0:1)}{y_0}
  \coximopd{5}{\ell_5-\ell_6}{(0:0:0:1)}{(0:1:0)}
}
\coxonem{E_6, E_7, E_8}{
  \coximp{6}{\ell_4}{x_0=x_1=0}{(0:0:1)}
  \coximp{7}{\ell_6}{x_0=x_2=0}{(0:1:0)}
  \coximd{8}{\ell_0-\ell_5-\ell_6}{x_2=x_3=0}{y_2}
}
\coxother{$(0)$-curves $E_9, E_{10}$}{
  \coxim{9}{\ell_0-\ell_1}{x_1=x_0x_3+x_2^2=0}{y_1}
  \coximds{10}{2\ell_0-\ell_1-\ell_2-\ell_3-\ell_4}{x_3=x_0x_2+x_1^2=0}{-y_0y_2-y_1^2}
}
\coxdynkinp{E_8 \cap E_9 \cap E_{10}}{(1:0:0:0)}{E_{10} \ar@{-}[rr] \ar@{=}[dr] \ar@{-}[dd] & & \li{6}  \ar@{-}[rr]& & \ex{2} \ar@{-}[dr]\\
      & \li{8} \ar@{-}[r] & \li{7} \ar@{-}[r] & \ex{5} \ar@{-}[r] & \ex{4} \ar@{-}[r] & \ex{1}\\
      E_9 \ar@{-}[rrrr] \ar@{-}[ur] & & & & \ex{3} \ar@{-}[ur]}
\coxonerel{\e_1, \dots, \e_{10}}{\e_2\e_6^2\e_{10} +
    \e_4\e_5^2\e_7^3\e_8 + \e_3\e_9^2}{2\ell_0-\ell_1-\ell_2}
  \coxanti{(\e_1^4\e_2^3\e_3^2\e_4^3\e_5^2\e_6^2\e_7,
    \e_1^3\e_2^2\e_3^2\e_4^2\e_5\e_6\e_9,
    \e_1^2\e_2\e_3\e_4^2\e_5^2\e_7^2\e_8, \e_8\e_{10})}


\coxcase{xix}{$\Afive+\Aone$}
\coxembed{\Pthree}{x_1^3+x_2x_3^2+x_0x_1x_2 = 0}
\coxbiratproj 3 {x_1:x_2:x_3}{(1:0:0:0)}
\coxtwotwo{\Afive}{(1:0:0:0)}{E_1, \dots, E_5}{\Aone}{(0:0:1:0)}{E_6}{
  \coximo{1}{\ell_0-\ell_1-\ell_5-\ell_6}{(1:0:0:0)}{y_0}
  \coximop{2}{\ell_1-\ell_2}{(1:0:0:0)}{(0:0:1)}
  \coximop{3}{\ell_2-\ell_3}{(1:0:0:0)}{(0:0:1)}
  \coximop{4}{\ell_3-\ell_4}{(1:0:0:0)}{(0:0:1)}
  \coximo{5}{\ell_0-\ell_1-\ell_2-\ell_3}{(1:0:0:0)}{y_1}
  \coximopd{6}{\ell_5-\ell_6}{(0:0:1:0)}{(0:1:0)}
}
\coxonem{E_7, E_8}{
  \coximp{7}{\ell_6}{x_1=x_3=0}{(0:1:0)}
  \coximpd{8}{\ell_4}{x_1=x_2=0}{(0:0:1)}
}
\coxother{$(0)$-curves $E_9$ and $(1)$-curve $E_{10}$}{
  \coxim{9}{\ell_0-\ell_5}{x_3=x_0x_2+x_1^2=0}{y_2}
  \coximds{10}{3\ell_0-\ell_1-\ell_2-\ell_3-\ell_4-2\ell_5}{x_0=x_1^3+x_2x_3^2=0}{-y_0^3-y_1y_2^2}
}
\coxdynkinp{E_6 \cap E_9 \cap E_{10}}{(0:0:1:0)}{E_9 \ar@{-}[rrrrr] \ar@{-}[dr] \ar@{-}[dd]& & & & & \ex{5} \ar@{-}[dr]\\
      & \ex{6} \ar@{-}[r] & \li{7} \ar@{-}[r] & \ex{1} \ar@{-}[r] & \ex{2} \ar@{-}[r] & \ex{3} \ar@{-}[r] & \ex{4}\\
      E_{10} \ar@{-}[rrrrr] \ar@{=}[ur] & & & & & \li{8} \ar@{-}[ur]}
\coxonerel{\e_1, \dots,\e_{10}}{\e_1^3\e_2^2\e_3\e_6\e_7^4 +
  \e_5\e_9^2 + \e_8\e_{10}}{3\ell_0-\ell_1-\ell_2-\ell_3-2\ell_5}
\coxanti{(\e_6\e_{10}, \e_1^2\e_2^2\e_3^2\e_4^2\e_5\e_6\e_7^2\e_8, \e_1\e_2^2\e_3^3\e_4^4\e_5^2\e_8^3, \e_1\e_2\e_3\e_4\e_5\e_6\e_7\e_9)}


\coxcase{xx}{$\Esix$}
\coxembed{\Pthree}{x_1x_2^2+x_2x_0^2+x_3^3=0}
\coxbiratproj 3 {x_0:x_2:x_3}{(0:1:0:0)}
\coxtwoone{\Esix}{(0:1:0:0)}{E_1, \dots, E_6}{
  \coximop{1}{\ell_1-\ell_2}{((0:1:0:0)}{(1:0:0)}
  \coximo{2}{\ell_0-\ell_1-\ell_2-\ell_3}{(0:1:0:0)}{y_1}
  \coximop{3}{\ell_2-\ell_3}{(0:1:0:0)}{(1:0:0)}
  \coximop{4}{\ell_5-\ell_6}{(0:1:0:0)}{(1:0:0)}
  \coximop{5}{\ell_4-\ell_5}{(0:1:0:0)}{(1:0:0)}
  \coximopd{6}{\ell_3-\ell_4}{(0:1:0:0)}{(1:0:0)}
}
\coxonem{E_7}{
  \coximpd{7}{\ell_6}{x_2=x_3=0}{(1:0:0)}
}
\coxother{$(0)$-curve $E_8$, $(1)$-curve $E_9$ and $(3)$-curve $E_{10}$}{
  \coxim{8}{\ell_0-\ell_1}{x_3=x_0^2+x_1x_2=0}{y_2}
  \coxim{9}{\ell_0}{x_0=x_1x_2^2+x_3^3=0}{y_0}
  \coximds{10}{3\ell_0-\ell_1-\ell_2-\ell_3-\ell_4-\ell_5-\ell_6}{x_1=x_0^2x_2+x_3^3=0}{-y_0^2y_1-y_2^3}
}
\coxdynkinp{E_8 \cap E_9 \cap E_{10}}{(0:0:1:0)}{E_{10} \ar@{-}[r] \ar@2{-}[dr] \ar@3{-}[dd] & \li{7} \ar@{-}[r] & \ex{4} \ar@{-}[r] & \ex{5} \ar@{-}[dr]\\
      & E_8 \ar@{-}[r] & \ex{1} \ar@{-}[r] & \ex{3} \ar@{-}[r] & \ex{6}\\
      E_9 \ar@{-}[rrr] \ar@{-}[ur] & & & \ex{2} \ar@{-}[ur]}
\coxonerel{\e_1, \dots, \e_{10}}
{\e_4^2\e_5\e_7^3\e_{10} + \e_2\e_9^2 + \e_1^2\e_3\e_8^3}
{3\ell_0-\ell_1-\ell_2-\ell_3}
\coxanti{(\e_1\e_2^2\e_3^2\e_4\e_5^2\e_6^3\e_9,
  \e_{10}, \e_1^2\e_2^3\e_3^4\e_4^4\e_5^5\e_6^6\e_7^3,
  \e_1^2\e_2^2\e_3^3\e_4^2\e_5^3\e_6^4\e_7\e_8)}

\subsection{Degree 2}\label{sec:degree_2}

In degree 2, the eleven types listed in Table~\ref{tab:degree_2} have
at most $13-d = 11$ negative curves. For simplicity, we do not list
the additional $35$ types with $12$ or more negative curves, whose Cox
rings have at least two relations by Theorem~\ref{thm:generators}.

\begin{table}[ht]
  \centering
\begin{equation*}
\begin{array}{cccc}
  \hline
  \text{singularities} & \text{$(-1)$-curves} &
  \text{type} & \text{reference}\\
  \hline
  \vdots & \vdots & \vdots & \\
  \Afive + \Aone & 5 & \tworel & \\
  \Dfive + \Aone & 5 & \onerel & \\
  \Esix & 4 & \onerel & \\
  \Asix & 4 & \tworel & \\
  \Dsix & 3 & \tworel & \\
  2\Athree+\Aone & 4 & \onerel & \\
  \Afive + \Atwo & 3 & \onerel & \\
  \Dfour + 3\Aone & 3 & \onerel & \\
  \Aseven & 2 & \tworel & \\
  \Dsix + \Aone & 2 & \onerel & \\
  \Eseven & 1 & \onerel & \text{\cite{arXiv:1011.3434}} \\
  \hline
\end{array}
\end{equation*}
  \smallskip
  \caption{Del Pezzo surfaces of degree $2$}
  \label{tab:degree_2}
\end{table}

For types $\Dsix, \Asix, \Afive+\Aone$, we apply
Proposition~\ref{prop:more_generators} as follows.  For type $\Dsix$, the extended
Dynkin diagram of negative curves is:
\begin{equation*}
  \dynkin{& \li 9 \ar@{-}[d] & & & \ex 6 \ar@{-}[d]\\
    \li 7 \ar@{-}[r] & \ex 1 \ar@{-}[r] & \ex 2 \ar@{-}[r] & \ex 3 \ar@{-}[r] & \ex 4 \ar@{-}[r] & \ex 5 \ar@{-}[r] & \li 8}
\end{equation*}
Contracting $E_8$ gives type $\Afive$ of degree $3$, whose Cox ring
has $\ge 2$ relations.
  
For type $\Asix$, the extended Dynkin diagram of negative curves is:
\begin{equation*}
  \dynkin{& \ex 1 \ar@{-}[r] \ar@{-}[d] & \li 8 \ar@{-}[r] & \li 7 \ar@{-}[r] & \ex 6 \ar@{-}[d]\\
    \li{10} \ar@{-}[r] & \ex 2 \ar@{-}[r] & \ex 3 \ar@{-}[r] & \ex 4 \ar@{-}[r] & \ex 5 \ar@{-}[r] & \li 9}
\end{equation*}
Contracting $E_8$ gives type $\Afive$ of degree $3$, whose Cox ring has $\ge
2$ relations.

For type $\Afive+\Aone$, the extended Dynkin diagram of negative curves is:
\begin{equation*}
  \dynkin{\li 7 \ar@{-}[r] & \ex 3 \ar@{-}[r] & \ex 4 \ar@{-}[r] & \ex 5 \ar@{-}[d]\\
    \ex 2 \ar@{-}[ur] \ar@{-}[d] & \li 8 \ar@{-}[r] & \li 9 \ar@{-}[ur] & \li{10}\\
    \ex 1 \ar@{-}[r] \ar@{-}[ur] & \li{11} \ar@{-}[r] & \ex 6 \ar@{-}[ur]}
\end{equation*}
Contracting $E_{11}$ gives type $\Afour$ of degree $3$, whose Cox ring has
$\ge 2$ relations.

For type $\Aseven$, the extended Dynkin diagram of negative curves is:
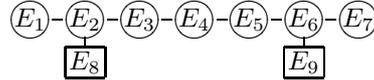
\begin{equation*}
  \dynkin{\ex 1 \ar@{-}[r] & \ex 2 \ar@{-}[r] \ar@{-}[d] & \ex 3 \ar@{-}[r] 
    & \ex 4 \ar@{-}[r] & \ex 5 \ar@{-}[r] & \ex 6 \ar@{-}[r] \ar@{-}[d] & \ex 7\\
    & \li 8 & & & & \li 9}
\end{equation*}
In Step~\ref{step:generator_degrees}, we discover at least three extra
generators in degrees $\cl E_2+\cl E_3+\cl E_4+\cl E_5+\cl E_6+\cl E_8+\cl E_9$,
$\cl E_1+2\cl E_2+2\cl E_3+2\cl E_4+2\cl E_5+2\cl E_6+\cl E_8+2\cl E_9$ and
$2\cl E_2+2\cl E_3+2\cl E_4+2\cl E_5+2\cl E_6+\cl E_7+2\cl E_8+\cl E_9$.

For the remaining seven types, we compute the Cox rings with precisely
one relation below.

By \cite[Theorems 1.2, 1.5]{MR1933881}, del Pezzo surfaces of types $\Eseven,
\Afive+\Atwo, 2\Athree+\Aone,\Dsix+\Aone,\Dfour+3\Aone$ are unique up to
isomorphism.
For types $\Esix$ and $\Dfive+\Aone$, we have two isomorphism
classes. Indeed, we cannot have more than two isomorphism classes
because blowing up $E_7 \cap E_9$ on type $\Esix$ gives a surface of
degree $1$, type $\Eseven+\Aone$, and blowing up $E_7 \cap E_{10}$ on
type $\Dfive+\Aone$ gives a surface of degree $1$, type
$\Esix+\Atwo$. Since types $\Eseven+\Aone, \Esix+\Atwo$ of degree $1$
have precisely two isomorphism classes, types $\Esix,\Dfive+\Aone$ of
degree $2$ can have at most two isomorphism classes. For each type, we
will give two surfaces below (depending on a parameter $\lambda \in
\{0,1\}$). The automorphism group will be infinite precisely for
$\lambda=0$, but we will also mention how their configurations of
negative curves differ geometrically, so that one clearly sees that
they cannot be isomorphic.


\coxcase{iv}{$\Dfive+\Aone$}
\coxembed{\PP(2,1,1,1)}{x_0^2+x_0x_1^2+x_1x_2^2x_3+\lambda x_0x_1x_2=0}
\coxbiratinv{x_0:x_1^2:x_1x_2}
{y_0y_1y_2^4:y_1y_2^2:y_2^3:-y_0^2y_1-y_0y_1^2-\lambda y_0y_1y_2}
\coxtwotwo{\Dfive}{(0:0:0:1)}{E_1, \dots, E_5}{\Aone}{(0:0:1:0)}{E_6}{
  \coximop{1}{\ell_6-\ell_7}{(0:0:0:1)}{(1:0:0)}
  \coximop{2}{\ell_5-\ell_6}{(0:0:0:1)}{(1:0:0)}
  \coximo{3}{\ell_0-\ell_1-\ell_3-\ell_5}{(0:0:0:1)}{y_2}
  \coximop{4}{\ell_1-\ell_2}{(0:0:0:1)}{(0:1:0)}
  \coximop{6}{\ell_3-\ell_4}{(0:0:0:1)}{(1:-1:0)}
  \coximod{5}{\ell_0-\ell_5-\ell_6-\ell_7}{(0:0:1:0)}{y_1}
}
\coxonem{E_7, \dots, E_{11}}{
  \coxims{7}{\ell_0-\ell_3-\ell_4}{x_3=x_0+x_1^2+\lambda x_1x_2=0}{-y_0-y_1-\lambda y_2}
  \coximp{8}{\ell_2}{x_0=x_2=0}{(0:1:0)}
  \coximp{9}{\ell_7}{x_0=x_1=0}{(1:0:0)}
  \coximp{10}{\ell_4}{x_2=x_0+x_1^2=0}{(1:-1:0)}
  \coximd{11}{\ell_0-\ell_1-\ell_2}{x_0=x_3=0}{y_0}
}
\coxdynkinpp{$E_6,E_7,E_{11}$ meet in one point $p$ with
  $\pi(p)=(0:0:1:0)$ if and only if $\lambda=0$}
{\li{7} \ar@{-}[rrr] \ar@{-}[dd] \ar@{-}[dr] & & & \li{{10}} \ar@{-}[r] & \ex{5} \ar@{-}[dr]\\
  & \ex{6} \ar@{-}[r] & \li{9} \ar@{-}[r] & \ex{1} \ar@{-}[r] & \ex{2} \ar@{-}[r] & \ex{3}\\
  \li{{11}} \ar@{-}[rrr] \ar@{-}[ur] & & & \li{8} \ar@{-}[r] & \ex{4}
  \ar@{-}[ur]}
\coxonerel{\e_1, \dots, \e_{11}} {\e_1^2\e_2\e_6\e_9^3 +
  \e_5\e_7\e_{10}^2 + \e_4\e_8^2\e_{11} + \lambda
  \e_1\e_2\e_3\e_4\e_5\e_8\e_9\e_{10}} {\ell_0}
\coxantie{
  (\e_1^2\e_2^3\e_3^4\e_4^3\e_5^2\e_6\e_8^2\e_9\e_{11},
  \e_1^2\e_2^2\e_3^2\e_4\e_5\e_6\e_9^2,
  \e_1\e_2^2\e_3^3\e_4^2\e_5^2\e_8\e_{10}, \e_6\e_7\e_{11})}
{\pi^*(-x_0-x_1^2-\lambda x_1x_2) = \e_1^2\e_2^3\e_3^4\e_4^2\e_5^3\e_6\e_7\e_9\e_{10}^2}


\coxcase{ii}{$\Esix$}
\coxembed{\PP(2,1,1,1)}{x_0^2+x_0x_3^2+x_1^3x_2+\lambda x_0x_1x_3=0}
\coxbiratinv{x_0:x_1^2:x_1x_3}
{y_0y_1^5:y_1^3:-y_0^2y_1-y_0y_2^2-\lambda y_0y_1y_2:y_1^2y_2}
\coxtwoone{\Esix}{(0:0:1:0)}{E_1, \dots, E_6}{
  \coximop{1}{\ell_6-\ell_7}{(0:0:1:0)}{(0:0:1)}
  \coximo{2}{\ell_0-\ell_1-\ell_2-\ell_6}{(0:0:1:0)}{y_1}
  \coximop{3}{\ell_2-\ell_3}{(0:0:1:0)}{(1:0:0)}
  \coximop{4}{\ell_3-\ell_4}{(0:0:1:0)}{(1:0:0)}
  \coximop{5}{\ell_4-\ell_5}{(0:0:1:0)}{(1:0:0)}
  \coximopd{6}{\ell_1-\ell_2}{(0:0:1:0)}{(1:0:0)}
}
\coxonem{E_7, \dots, E_{10}}{
  \coxim{7}{\ell_0-\ell_6-\ell_7}{x_0=x_2=0}{y_0}
  \coximp{8}{\ell_5}{x_1=x_0+x_3^2=0}{(1:0:0)}
  \coximp{9}{\ell_7}{x_0=x_1=0}{(0:0:1)}
  \coximds{10}{2\ell_0-\ell_1-\ell_2-\ell_3-\ell_4-\ell_5}{x_2=x_0+x_3^2+\lambda x_1x_3=0}{-y_0y_1-y_2^2-\lambda y_1y_2}
}
\coxother{$(0)$-curve $E_{11}$}{
  \coximd{11}{\ell_0-\ell_1}{x_3=x_0^2+x_1^3x_2=0}{y_2} }
\coxdynkinpp{with $E_7,E_{10}$ touching in $p$ and meeting $E_{11}$
  transversally also in $p$ with $\pi(p)=(0:1:0:0)$ if $\lambda=0$, and with
  $E_7,E_{10},E_{11}$ meeting in $p$ with $\pi(p)=(0:1:0:0)$ transversally and
  $E_7,E_{10}$ meeting also in $p'$ with $\pi(p')=(0:1:0:-1)$ if $\lambda=1$}{E_{11} \ar@{-}[rrrr] \ar@{-}[dd] \ar@{-}[dr] & & & & \ex{6} \ar@{-}[dr]\\
    & \li{{10}} \ar@{-}[r] & \li{8} \ar@{-}[r] & \ex{5} \ar@{-}[r] & \ex{4} \ar@{-}[r] & \ex{3}\\
    \li{7} \ar@{=}[ur] \ar@{-}[rr] & & \li{9} \ar@{-}[r] & \ex{1} \ar@{-}[r] & \ex{2}
    \ar@{-}[ur]}
\coxonerel{\e_1, \dots, \e_{11}}
{\e_1^2\e_2\e_7\e_9^3+\e_4\e_5^2\e_8^3\e_{10} + \e_6\e_{11}^2 + \lambda \e_1\e_2\e_3\e_4\e_5\e_6\e_8\e_9\e_{11}}
{2\ell_0-\ell_1-\ell_2}
\coxanti{
  (\e_1^4\e_2^5\e_3^6\e_4^4\e_5^2\e_6^3\e_7\e_9^3,
  \e_1^2\e_2^3\e_3^4\e_4^3\e_5^2\e_6^2\e_8\e_9, \e_7\e_{10},
  \e_1\e_2^2\e_3^3\e_4^2\e_5\e_6^2\e_{11})}


\coxcase{ix}{$2\Athree+\Aone$}
\coxembed{\PP(2,1,1,1)}{x_0^2+x_0x_1x_3+x_1x_2^2x_3=0}
\coxbiratinv{x_0:x_1^2:x_1x_2}
{y_0y_1(y_0y_1+y_2^2)^2:y_1(y_0y_1+y_2^2):y_2(y_0y_1+y_2^2):-y_0^2y_1}
\coxtwomulti{\Athree,\Athree,\Aone}{(0:0:0:1),(0:1:0:0),(0:0:1:0)}{$E_1,E_2,E_3$ and $E_4,E_5,E_6$ and $E_7$}{
  \coximo{1}{2\ell_0-\ell_1-\ell_2-\ell_3-\ell_4-\ell_5-\ell_6}{(0:0:0:1)}{-y_0y_1-y_2^2}
  \coximop{2}{\ell_6-\ell_7}{(0:0:0:1)}{(1:0:0)}
  \coximop{3}{\ell_5-\ell_6}{(0:0:0:1)}{(1:0:0)}
  \coximop{4}{\ell_3-\ell_4}{(0:1:0:0)}{(0:1:0)}
  \coximop{5}{\ell_2-\ell_3}{(0:1:0:0)}{(0:1:0)}
  \coximop{6}{\ell_1-\ell_2}{(0:1:0:0)}{(0:1:0)}
  \coximod{7}{\ell_0-\ell_5-\ell_6-\ell_7}{(0:0:1:0)}{y_1}
}
\coxonem{E_8, \dots, E_{11}}{
  \coximp{8}{\ell_7}{x_0=x_1=0}{(1:0:0)}
  \coximp{9}{\ell_4}{x_0=x_2=0}{(0:1:0)}
  \coxim{10}{\ell_0-\ell_1-\ell_5}{x_2=x_0+x_1x_3=0}{y_2}
  \coximd{11}{\ell_0-\ell_1-\ell_2}{x_0=x_3=0}{y_0}
}
\coxdynkin{&\ex 1 \ar@{-}[r] \ar@{-}[dl] & \li 9
  \ar@{-}[r] & \ex 4 \ar@{-}[dr]\\
    \ex 2 \ar@{-}[r] \ar@{-}[dr] & \li 8 \ar@{-}[r] & \ex 7 \ar@{-}[r] &
    \li{11} \ar@{-}[r] & \ex 5 \ar@{-}[dl]\\
    & \ex 3 \ar@{-}[r] & \li{10} \ar@{-}[r] & \ex 6}
\coxonerel{\e_1, \dots, \e_{11}}
{\e_1\e_4\e_9^2 + \e_3\e_6\e_{10}^2 + \e_7\e_8\e_{11}}
{2\ell_0-\ell_1-\ell_2-\ell_5-\ell_6}
\coxanti{
  (&\e_1^2\e_2^2\e_3\e_4^2\e_5^2\e_6\e_7\e_8\e_9^2\e_{11},
  \e_1\e_2^2\e_3\e_7\e_8^2,
  \e_1\e_2\e_3\e_4\e_5\e_6\e_9\e_{10}, \e_4\e_5^2\e_6\e_7\e_{11}^2)}


\coxcase{viii}{$\Afive+\Atwo$}
\coxembed{\PP(2,1,1,1)}{x_0^2+x_0x_2x_3+x_1^3x_3=0}
\coxbiratinv{x_0:x_1^2:x_1x_3}
{y_0^3y_1y_2^2:y_0y_1y_2:-y_0^2y_1-y_1^2y_2:y_0y_2^2}
\coxtwotwo{\Afive}{(0:0:1:0)}{E_1, \dots, E_5}{\Atwo}{(0:0:0:1)}{E_6,E_7}{
  \coximo{1}{\ell_0-\ell_3-\ell_4-\ell_5}{(0:0:1:0)}{y_0}
  \coximop{2}{\ell_5-\ell_6}{(0:0:1:0)}{(0:1:0)}
  \coximop{3}{\ell_6-\ell_7}{(0:0:1:0)}{(0:1:0)}
  \coximo{4}{\ell_0-\ell_1-\ell_5-\ell_6}{(0:0:1:0)}{y_2}
  \coximop{5}{\ell_1-\ell_2}{(0:0:1:0)}{(1:0:0)}
  \coximop{6}{\ell_3-\ell_4}{(0:0:0:1)}{(0:0:1)}
  \coximod{7}{\ell_0-\ell_1-\ell_2-\ell_3}{(0:0:0:1)}{y_1}
}
\coxonem{E_8, E_9, E_{10}}{
  \coximp{8}{\ell_4}{x_0=x_1=0}{(0:0:1)}
  \coximp{9}{\ell_2}{x_1=x_0+x_2x_3=0}{(1:0:0)}
  \coximd{10}{\ell_7}{x_0=x_3=0}{(0:1:0)}
}
\coxother{$(0)$-curve $E_{11}$}{
  \coximds{11}{2\ell_0-\ell_3-\ell_5-\ell_6-\ell_7}{x_2=x_0^2+x_1^3x_3=0}{-(y_0^2+y_1y_2)} }
\coxdynkinp{E_6 \cap E_7 \cap E_{11}}{(0:0:0:1)}{\ex{6} \ar@{-}[r] \ar@{-}[dr] \ar@{-}[dd] & \li 8 \ar@{-}[r]
  & \ex 1 \ar@{-}[r] & \ex 2 \ar@{-}[dr]\\
    & E_{11} \ar@{-}[r] & \li{10} \ar@{-}[rr] & & \ex 3 \ar@{-}[dl]\\
    \ex 7 \ar@{-}[r] \ar@{-}[ur] & \li 9 \ar@{-}[r] & \ex 5 \ar@{-}[r] & \ex 4}
  \coxonerel{\e_1, \dots, \e_{11}} {\e_1^2\e_2\e_6\e_8^3 +
    \e_4\e_5^2\e_7\e_9^3 + \e_{10}\e_{11}}
  {2\ell_0-\ell_3-\ell_5-\ell_6}
\coxanti{
  (\e_1^3\e_2^3\e_3^3\e_4^2\e_5\e_6^2\e_7\e_8^3\e_{10},
  \e_1\e_2\e_3\e_4\e_5\e_6\e_7\e_8\e_9,
  \e_6\e_7\e_{11}, \e_1\e_2^2\e_3^3\e_4^2\e_5\e_{10}^2)}


\coxcase{xi}{$\Dfour+3\Aone$}
\coxembed{\PP(2,1,1,1)}{x_0^2+x_1^2x_2x_3+x_1x_2^2x_3=0}
\coxbiratinv{x_0:x_1^2:x_1x_2}
{y_0y_1y_2^2(y_1+y_2)^2:y_1y_2(y_1+y_2):y_2^2(y_1+y_2):-y_0^2y_1}
\coxtwomulti{\Dfour,\Aone,\Aone,\Aone}{(0:0:0:1),(0:0:1:0),(0:1:0:0),(0:1:-1:0)}{$E_1, \dots, E_4$ and $E_5$ and $E_6$ and $E_7$}{
  \coximop{1}{\ell_1-\ell_2}{(0:0:0:1)}{(1:0:0)}
  \coximop{2}{\ell_2-\ell_3}{(0:0:0:1)}{(1:0:0)}
  \coximo{3}{\ell_0-\ell_1-\ell_4-\ell_5}{(0:0:0:1)}{y_2}
  \coximo{4}{\ell_0-\ell_1-\ell_6-\ell_7}{(0:0:0:1)}{-y_1-y_2}
  \coximo{5}{\ell_0-\ell_1-\ell_2-\ell_3}{(0:0:1:0)}{y_1}
  \coximop{6}{\ell_4-\ell_5}{(0:1:0:0)}{(0:1:0)}
  \coximopd{7}{\ell_6-\ell_7}{(0:1:-1:0)}{(0:1:-1)}
}
\coxonem{E_8, \dots, E_{11}}{
  \coximp{8}{\ell_5}{x_0=x_2=0}{(0:1:0)} 
  \coxim{9}{\ell_0-\ell_4-\ell_6}{x_0=x_3=0}{y_0}
  \coximp{10}{\ell_3}{x_0=x_1=0}{(1:0:0)}
  \coximd{11}{\ell_7}{x_0=x_1+x_2=0}{(0:1:-1)} 
}
\coxdynkin{& \ex 5 \ar@{-}[r] \ar@{-}[dl] & \li{10} \ar@{-}[r] & \ex 2 \ar@{-}[dr]\\
    \li 9 \ar@{-}[r] \ar@{-}[dr] & \ex 6 \ar@{-}[r] & \li 8 \ar@{-}[r] & \ex 3
    \ar@{-}[r] & \ex 1\\
    & \ex 7 \ar@{-}[r] & \li{11} \ar@{-}[r] & \ex 4 \ar@{-}[ur]}
\coxonerel{\e_1, \dots, \e_{11}}
{\e_2\e_5\e_{10}^2+\e_3\e_6\e_8^2+\e_4\e_7\e_{11}^2}
{\ell_0-\ell_1}
\coxantie{
  (\e_1^3\e_2^2\e_3^2\e_4^2\e_5\e_6\e_7\e_8\e_9\e_{10}\e_{11},
  \e_1^2\e_2^2\e_3\e_4\e_5\e_{10}^2,
  \e_1^2\e_2\e_3^2\e_4\e_6\e_8^2, \e_5\e_6\e_7\e_9^2)}
{\pi^*(-x_1-x_2) = \e_1^2\e_2\e_3\e_4^2\e_7\e_{11}^2}


\coxcase{x}{$\Dsix+\Aone$}
\coxembed{\PP(2,1,1,1)}{x_0^2+x_1x_2^3+x_1^2x_2x_3=0}
\coxbiratinv{x_0:x_1^2:x_1x_2}
{y_0y_1(y_0y_1+y_2^2)^2:y_1(y_0y_1+y_2^2):y_2(y_0y_1+y_2^2):-y_0^2y_1}
\coxtwotwo{\Dsix}{(0:0:0:1)}{E_1, \dots, E_6}{\Aone}{(0:1:0:0)}{E_7}{
  \coximo{1}{\ell_0-\ell_1-\ell_6-\ell_7}{(0:0:0:1)}{y_2}
  \coximop{2}{\ell_1-\ell_2}{(0:0:0:1)}{(1:0:0)}
  \coximop{3}{\ell_2-\ell_3}{(0:0:0:1)}{(1:0:0)}
  \coximop{4}{\ell_3-\ell_4}{(0:0:0:1)}{(1:0:0)}
  \coximop{5}{\ell_4-\ell_5}{(0:0:0:1)}{(1:0:0)}
  \coximo{6}{\ell_0-\ell_1-\ell_2-\ell_3}{(0:0:0:1)}{y_1}
  \coximopd{7}{\ell_6-\ell_7}{(0:1:0:0)}{(0:1:0)}
}
\coxonem{E_8, E_9}{
  \coximp{8}{\ell_5}{x_0=x_1=0}{(1:0:0)}
  \coximpd{9}{\ell_7}{x_0=x_2=0}{(0:1:0)}
}
\coxother{$(0)$-curves $E_{10}, E_{11}$}{
  \coxim{10}{\ell_0-\ell_6}{x_0=x_2^2+x_1x_3=0}{y_0} 
  \coximds{11}{3\ell_0-\ell_1-\ell_2-\ell_3-\ell_4-\ell_5-2\ell_6}{x_3=x_0^2+x_1x_2^3=0}{-y_0^2y_1-y_2^3} 
}
\coxdynkinp{E_7 \cap E_{10} \cap E_{11}}{(0:1:0:0)}{E_{10} \ar@{-}[rrrrr] \ar@{-}[dd] \ar@{-}[dr] & & & & & \ex 6 \ar@{-}[dr]\\
    & \ex 7 \ar@{-}[r] & \li 9 \ar@{-}[r] & \ex 1 \ar@{-}[r] & \ex 2
    \ar@{-}[r] & \ex 3 \ar@{-}[r] & \ex 4\\
    E_{11} \ar@{-}[rrrr] \ar@{=}[ur] & & & & \li 8 \ar@{-}[r] & \ex 5
    \ar@{-}[ur]}
\coxonerel{\e_1, \dots, \e_{11}}
{\e_1^3\e_2^2\e_3\e_7\e_9^4 + \e_6\e_{10}^2 +
    \e_5\e_8^2\e_{11}}
{3\ell_0-\ell_1-\ell_2-\ell_3-2\ell_6}
\coxanti{
  (\e_1^2\e_2^3\e_3^4\e_4^5\e_5^3\e_6^3\e_7\e_8\e_9\e_{10},
  \e_1\e_2^2\e_3^3\e_4^4\e_5^3\e_6^2\e_8^2,
  \e_1^2\e_2^2\e_3^2\e_4^2\e_5\e_6\e_7\e_9^2, \e_7\e_{11})}


\coxcase{i}{$\Eseven$}
\coxembed{\PP(2,1,1,1)}{x_0^2+x_1x_2^3+x_1^3x_3=0}
\coxbiratinv {x_0:x_1^2:x_1x_2}{y_0y_1^5:y_1^3:y_1^2y_2:-y_0^2y_1-y_2^3}
\coxtwoone{\Eseven}{(0:0:0:1)}{E_1, \dots, E_7}{
  \coximop{1}{\ell_1-\ell_2}{(0:0:0:1)}{(1:0:0)}
  \coximop{2}{\ell_2-\ell_3}{(0:0:0:1)}{(1:0:0)}
  \coximop{3}{\ell_3-\ell_4}{(0:0:0:1)}{(1:0:0)}
  \coximop{4}{\ell_4-\ell_5}{(0:0:0:1)}{(1:0:0)}
  \coximop{5}{\ell_5-\ell_6}{(0:0:0:1)}{(1:0:0)}
  \coximop{6}{\ell_6-\ell_7}{(0:0:0:1)}{(1:0:0)}
  \coximod{7}{\ell_0-\ell_1-\ell_2-\ell_3}{(0:0:0:1)}{y_1}
}
\coxonem{E_8}{
  \coximpd{8}{\ell_7}{x_0=x_1=0}{(1:0:0)}
}
\coxother{$(0)$-curve $E_9$, $(1)$-curve $E_{10}$ and $(2)$-curve $E_{11}$}{
  \coxim{9}{\ell_0-\ell_1}{x_2=x_0^2+x_1^3x_3=0}{y_2}
  \coxim{10}{\ell_0}{x_0=x_2^3+x_1^2x_3=0}{y_0}
  \coximds{11}{3\ell_0-\ell_1-\ell_2-\ell_3-\ell_4-\ell_5-\ell_6-\ell_7}{x_3=x_0^2+x_1x_2^3=0}{-y_0^2y_1-y_2^3}
}
\coxdynkinp{E_9 \cap E_{10} \cap E_{11}}{(0:1:0:0)}{E_9 \ar@{-}[rrrr] \ar@{=}[dr] \ar@{-}[dd] & & & & \ex{1} \ar@{-}[r] & \ex{2} \ar@{-}[dr]\\
    & E_{11} \ar@{-}[r] & \li{8} \ar@{-}[r] & \ex{6} \ar@{-}[r] & \ex{5} \ar@{-}[r] & \ex{4} \ar@{-}[r] & \ex{3} \\
    E_{10} \ar@{-}[rrrrr] \ar@3{-}[ur] & & & & & \ex{7} \ar@{-}[ur]}
\coxonerel{\e_1, \dots, \e_{11}}
{\e_1^2\e_2\e_9^3 + \e_7\e_{10}^2 +
    \e_4\e_5^2\e_6^3\e_8^4\e_{11}}
{3\ell_0-\ell_1-\ell_2-\ell_3}
\coxanti{
  (\e_1^3\e_2^6\e_3^9\e_4^7\e_5^5\e_6^3\e_7^5\e_8\e_{10},
  \e_1^2\e_2^4\e_3^6\e_4^5\e_5^4\e_6^3\e_7^3\e_8^2,
  \e_1^2\e_2^3\e_3^4\e_4^3\e_5^2\e_6\e_7^2\e_9, \e_{11})}


\subsection{Degree 1}\label{sec:degree_1}

We list the seven types of generalized del Pezzo surfaces of degree
$1$ with at most $13-d=12$ negative curves in
Table~\ref{tab:degree_1}. Additionally, there are $67$ types with more
than $12$ negative curves.

\begin{table}[ht]
  \centering
  \begin{equation*}
    \begin{array}{ccc}
      \hline
      \text{singularities} & \text{$(-1)$-curves} & \text{type}\\
      \hline
      \vdots & \vdots & \vdots \\
      \Dseven & 5 & \tworel\\
      \Eseven & 5 & \tworel\\
      \Esix + \Atwo & 4 & \onerel\\
      \Aeight & 3 & \tworel\\
      \Eseven + \Aone & 3 & \onerel\\
      \Deight & 2 & \tworel\\
      \Eeight & 1 & \onerel\\
      \hline
    \end{array}
  \end{equation*}
  \smallskip
  \caption{Del Pezzo surfaces of degree $1$}
  \label{tab:degree_1}
\end{table}

For types $\Dseven, \Eseven, \Aeight, \Deight$, we apply
Proposition~\ref{prop:more_generators} to show that the Cox rings have $\ge 2$
relations.  For type $\Dseven$, the extended Dynkin diagram of negative curves is:
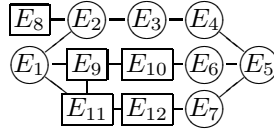
\begin{equation*}
  \dynkin{\li 8 \ar@{-}[r] & \ex 2 \ar@{-}[r] & \ex 3 \ar@{-}[r] & \ex 4 \ar@{-}[dr]\\
    \ex 1 \ar@{-}[ur] \ar@{-}[r] \ar@{-}[dr] & \li 9 \ar@{-}[r] \ar@{-}[d] & \li{10} \ar@{-}[r] & \ex 6 \ar@{-}[r] & \ex 5\\
    & \li{11} \ar@{-}[r] & \li{12} \ar@{-}[r] & \ex 7 \ar@{-}[ur]}
\end{equation*}
Contracting $E_9$ gives type $\Dsix$ of degree $2$, whose Cox ring has $\ge
2$ relations.
  
For type $\Eseven$, the extended Dynkin diagram of negative curves is:
  \begin{equation*}
    \dynkin{\li{12} \ar@{=}[r] \ar@{-}[dd] \ar@3{-}[dr] & \li 8 \ar@{-}[r] &
      \ex 7 \ar@{-}[r] \ar@{-}[d] & \ex 6 \ar@{-}[d] \\
      & \li 9 \ar@{=}[r] & \li{11} & \ex 5 \ar@{-}[d]\\
      \li{10} \ar@{-}[r] \ar@{-}[ur] & \ex 2 \ar@{-}[r] & \ex 3 \ar@{-}[r] & \ex 4 \ar@{-}[r] & \ex 1}
  \end{equation*}
  Contracting $E_{10}$ leads to type $\Dsix$ of degree $2$, whose Cox
  ring has $\ge 2$ relations.

  For type $\Aeight$, the extended Dynkin diagram of negative curves is:
  \begin{equation*}
    \dynkin{& \ex 2 \ar@{-}[r] \ar@{-}[d] & \ex 1 \ar@{-}[r] & \li 9 \ar@{-}[r] & \ex 8 \ar@{-}[r] & \ex 7 \ar@{-}[d]\\
    \li{10} \ar@{-}[r] & \ex 3 \ar@{-}[r] & \ex 4 \ar@{-}[rr] & & \ex 5 \ar@{-}[r] & \ex 6 \ar@{-}[r] & \li{11}}
  \end{equation*}
  Contracting $E_9$ gives type $\Asix$ of degree $2$, whose Cox ring has $\ge
  2$ relations.
  
  For type $\Deight$, the extended Dynkin diagram of negative curves is:
  \begin{equation*}
    \dynkin{\ex 1 \ar@{-}[r] & \ex 2 \ar@{-}[r] \ar@{-}[d] & \ex 3 \ar@{-}[r] & \ex 4 \ar@{-}[r] & \ex 5 \ar@{-}[r] & \ex 6 \ar@{-}[r] \ar@{-}[d] & \ex 7 \ar@{-}[r] & \li{10}\\
    & \li 9 & & & & \ex 8}
  \end{equation*}
  Contracting $E_{10}$ gives type $\Aseven$ of degree $2$, whose Cox ring has
  $\ge 2$ relations.

  According to \cite[Theorem 1.2]{MR1933881}, for each of the types
  $\Esix+\Atwo, \Eseven+\Aone, \Eeight$, there exist two isomorphy
  classes. They will be given below using a parameter $\lambda \in
  \{0,1\}$. Geometrically, the two isomorphism classes of each type differ as
  follows: for $\lambda=1$, their automorphism group is finite; for
  $\lambda=0$, it is infinite.


\coxcase{vi}{$\Esix+\Atwo$}
\coxembed{\PP(3,2,1,1)}{x_0^2+x_0x_2^2x_3+x_1^3+\lambda x_0x_1x_2=0}
\coxbiratinv{x_0:x_1x_2:x_2^3}
{y_0y_2^8:y_1y_2^5:y_2^3:-y_0^2y_2-y_1^3-\lambda y_0y_1y_2}
\coxtwotwo{\Esix}{(0:0:0:1)}{E_1, \dots, E_6}{\Atwo}{(0:0:1:0)}{E_7,E_8}{
  \coximo{1}{\ell_0-\ell_6-\ell_7-\ell_8}{(0:0:0:1)}{y_0}
  \coximo{2}{\ell_0-\ell_1-\ell_2-\ell_3}{(0:0:0:1)}{y_2}
  \coximop{3}{\ell_3-\ell_4}{(0:0:0:1)}{(1:0:0)}
  \coximop{4}{\ell_2-\ell_3}{(0:0:0:1)}{(1:0:0)}
  \coximop{5}{\ell_1-\ell_2}{(0:0:0:1)}{(1:0:0)}
  \coximop{6}{\ell_4-\ell_5}{(0:0:0:1)}{(1:0:0)}
  \coximop{7}{\ell_7-\ell_8}{(0:0:1:0)}{(0:0:1)}
  \coximopd{8}{\ell_6-\ell_7}{(0:0:1:0)}{(0:0:1)}
}
\coxonem{E_9, \dots, E_{12}}{
  \coximp{9}{\ell_5}{x_2=x_0^2+x_1^3=0}{(1:0:0)}
  \coximp{10}{\ell_8}{x_0=x_1=0}{(0:0:1)}
  \coxim{11}{\ell_0-\ell_1-\ell_6}{x_1=x_0+x_2^2x_3=0}{y_1}
  \coximdss{12}{3\ell_0-\ell_1-\ell_2-\ell_3-\ell_4-\ell_5-2\ell_6-\ell_7}{x_3=x_0^2+x_1^3+\lambda x_0x_1x_2=0}{-y_0^2y_2-y_1^3-\lambda y_0y_1y_2}
}
\coxdynkinpp{for $\lambda=0$, $\pi(E_7 \cap E_8 \cap E_{12}) = (0:0:1:0)$; for
  $\lambda=1$, $E_7,E_8,E_{12}$ meet pairwise in three different points that
  all project to $(0:0:1:0)$}{\ex 7 \ar@{-}[r] \ar@{-}[dd] \ar@{-}[dr] & \li{10} \ar@{-}[r] & \ex 1
  \ar@{-}[r] & \ex 2 \ar@{-}[dr]\\
  & \li{12} \ar@{-}[r] & \li 9 \ar@{-}[r] & \ex 6 \ar@{-}[r] & \ex 3 \ar@{-}[dl]\\
  \ex 8 \ar@{-}[r] \ar@{-}[ur] & \li{11} \ar@{-}[r] & \ex 5 \ar@{-}[r] & \ex
  4}
\coxonerel{\e_1, \dots, \e_{12}} {\e_1^2\e_2\e_7\e_{10}^3 +
  \e_4\e_5^2\e_8\e_{11}^3 + \e_6\e_9^2\e_{12}+ \lambda
  \e_1\e_2\e_3\e_4\e_5\e_6\e_9\e_{10}\e_{11}}
{3\ell_0-\ell_1-\ell_2-\ell_3-2\ell_6-\ell_7}
\coxanti{
  (\e_1^4\e_2^5\e_3^6\e_4^4\e_5^2\e_6^3\e_7^2\e_8\e_{10}^3,
  \e_1^2\e_2^3\e_3^4\e_4^3\e_5^2\e_6^2\e_7\e_8\e_{10}\e_{11},
  \e_1\e_2^2\e_3^3\e_4^2\e_5\e_6^2\e_9, \e_7\e_8\e_{12})}


\coxcase{vii}{$\Eseven+\Aone$}
\coxembed{\PP(3,2,1,1)}{x_0^2+x_1^3+x_1x_2^3x_3+\lambda x_0x_1x_2=0}
\coxbiratinv{x_0:x_1x_2:x_2^3}
{y_0y_2^8:y_1y_2^5:y_2^3:-y_0^2y_2-y_1^3-\lambda y_0y_1y_2}
\coxtwotwo{\Eseven}{(0:0:0:1)}{E_1, \dots, E_7}{\Aone}{(0:0:1:0)}{E_8}{
  \coximop{1}{\ell_5-\ell_6}{(0:0:0:1)}{(1:0:0)}
  \coximop{2}{\ell_4-\ell_5}{(0:0:0:1)}{(1:0:0)}
  \coximop{3}{\ell_3-\ell_4}{(0:0:0:1)}{(1:0:0)}
  \coximop{4}{\ell_2-\ell_3}{(0:0:0:1)}{(1:0:0)}
  \coximop{5}{\ell_1-\ell_2}{(0:0:0:1)}{(1:0:0)}
  \coximo{6}{\ell_0-\ell_1-\ell_7-\ell_8}{(0:0:0:1)}{y_1}
  \coximo{7}{\ell_0-\ell_1-\ell_2-\ell_3}{(0:0:0:1)}{y_2}
  \coximopd{8}{\ell_7-\ell_8}{(0:0:1:0)}{(0:0:1)}
}
\coxonem{E_9, E_{10}, E_{11}}{
  \coximp{9}{\ell_8}{x_0=x_1=0}{(0:0:1)} 
  \coximp{10}{\ell_6}{x_2=x_0^2+x_1^3=0}{(1:0:0)} 
  \coximdss{11}{3\ell_0-\ell_1-\ell_2-\ell_3-\ell_4-\ell_5-\ell_6-2\ell_7}{x_3=x_0^2+x_1^3+\lambda
    x_0x_1x_2=0}{-y_0^2y_2-y_1^3-\lambda y_0y_1y_2}
}
\coxother{$(0)$-curve $E_{12}$}{
  \coximd{12}{\ell_0-\ell_7}{x_0=x_1^2+x_2^3x_3=0}{y_0} 
}
\coxdynkinpp{$\pi(E_8\cap E_{11}\cap E_{12}) = (0:0:1:0)$; for $\lambda=0$,
  $E_8,E_{11}$ touch in one point where they also meet $E_{12}$ transversally;
  for $\lambda=1$, $E_8,E_{11}$ meet transversally in two points; in one of
  them, they also meet $E_{12}$}
{\li{11} \ar@{-}[rrr] \ar@{=}[dr] \ar@{-}[dd] & & & \li{10} \ar@{-}[r] & \ex 1
  \ar@{-}[r] & \ex 2 \ar@{-}[dr]\\
  & \ex 8 \ar@{-}[r] & \li 9 \ar@{-}[r] & \ex 6 \ar@{-}[r] & \ex 5 \ar@{-}[r]
  & \ex 4 \ar@{-}[r] & \ex 3\\
  E_{12} \ar@{-}[ur] \ar@{-}[rrrrr] & & & & & \ex 7 \ar@{-}[ur]}
\coxonerel{\e_1, \dots, \e_{12}}
{\e_4\e_5^2\e_6^3\e_8\e_9^4 +
    \e_1^2\e_2\e_{10}^3\e_{11} + \e_7\e_{12}^2+ \lambda
    \e_1\e_2\e_3\e_4\e_5\e_6\e_7\e_9\e_{10}\e_{12}}
{3\ell_0-\ell_1-\ell_2-\ell_3-2\ell_7}
\coxanti{
  (\e_1^3\e_2^6\e_3^9\e_4^7\e_5^5\e_6^3\e_7^5\e_8\e_9\e_{12},
  \e_1^2\e_2^4\e_3^6\e_4^5\e_5^4\e_6^3\e_7^3\e_8\e_9^2,
  \e_1^2\e_2^3\e_3^4\e_4^3\e_5^2\e_6\e_7^2\e_{10}, \e_8\e_{11})}


\coxcase{i}{$\Eeight$}
\coxembed{\PP(3,2,1,1)}{x_0^2+x_1^3+x_2^5x_3+\lambda x_0x_1x_2=0}
\coxbiratinv{x_0:x_1x_2:x_2^3}
{y_0y_2^8:y_1y_2^5:y_2^3:-y_0^2y_2-y_1^3-\lambda y_0y_1y_2}
\coxtwoone{\Eeight}{(0:0:0:1)}{E_1,\dots, E_8}{
  \coximop{1}{\ell_1-\ell_2}{(0:0:0:1)}{(1:0:0)}
  \coximop{2}{\ell_2-\ell_3}{(0:0:0:1)}{(1:0:0)}
  \coximop{3}{\ell_3-\ell_4}{(0:0:0:1)}{(1:0:0)}
  \coximop{4}{\ell_4-\ell_5}{(0:0:0:1)}{(1:0:0)}
  \coximop{5}{\ell_5-\ell_6}{(0:0:0:1)}{(1:0:0)}
  \coximop{6}{\ell_6-\ell_7}{(0:0:0:1)}{(1:0:0)}
  \coximop{7}{\ell_7-\ell_8}{(0:0:0:1)}{(1:0:0)}
  \coximod{8}{\ell_0-\ell_1-\ell_2-\ell_3}{(0:0:0:1)}{y_2}
}
\coxonem{E_9}{
  \coximpd{9}{\ell_8}{x_2=x_0^2+x_1^3=0}{(1:0:0)}
}
\coxother{$(0)$-curve $E_{10}=$ and $(1)$-curves $E_{11}, E_{12}$}{
  \coxim{10}{\ell_0-\ell_1}{x_1=x_0^2+x_2^5x_3=0}{y_1} 
  \coximss{11}{3\ell_0-\ell_1-\ell_2-\ell_3-\ell_4-\ell_5-\ell_6-\ell_7-\ell_8}{x_3=x_0^2+x_1^3+\lambda x_0x_1x_2=0}{-y_0^2y_2-y_1^3-\lambda y_0y_1y_2} 
  \coximd{12}{\ell_0}{x_0=x_1^3+x_2^5x_3=0}{y_0} 
}
\coxdynkinpp{with $E_{10},E_{11},E_{12}$ meeting in $p$ with
  $\pi(p)=(0:1:0:0)$; the automorphism group is infinite precisely for $\lambda=0$}
{E_{10} \ar@{-}[rrrrr] \ar@{=}[dr] \ar@{-}[dd] & & & & & \ex 1 \ar@{-}[r] & \ex 2 \ar@{-}[dr]\\
  & E_{11} \ar@{-}[r] & \li 9 \ar@{-}[r] & \ex 7 \ar@{-}[r] & \ex 6 \ar@{-}[r] & \ex 5 \ar@{-}[r] & \ex 4 \ar@{-}[r] & \ex 3\\
  E_{12} \ar@{-}[rrrrrr] \ar@3{-}[ur] & & & & & & \ex 8 \ar@{-}[ur]}
\coxonerel{\e_1, \dots, \e_{12}}
{\e_1^2\e_2\e_{10}^3 +
    \e_4\e_5^2\e_6^3\e_7^4\e_9^5\e_{11} +
    \e_8\e_{12}^2+\lambda \e_1\e_2\e_3\e_4\e_5\e_6\e_7\e_8\e_9\e_{10}\e_{12}}
{3\ell_0-\ell_1-\ell_2-\ell_3}
\coxanti{
  (\e_1^5\e_2^{10}\e_3^{15}\e_4^{12}\e_5^9\e_6^6\e_7^3\e_8^8\e_{12},
  \e_1^4\e_2^7\e_3^{10}\e_4^8\e_5^6\e_6^4\e_7^2\e_8^5\e_{10},
  \e_1^2\e_2^4\e_3^6\e_4^5\e_5^4\e_6^3\e_7^2\e_8^3\e_9, \e_{11})}


\bibliographystyle{alpha}

\bibliography{cox_hypersurface_3}

\begin{thebibliography}{TVAV11}

\bibitem[ADHL]{arXiv:1003.4229}
I.~Arzhantsev, U.~Derenthal, J.~Hausen, and A.~Laface.
\newblock Cox rings.
\newblock {\em arXiv:1003.4229,
  http://www.mathematik.uni-tuebingen.de/\~{}hausen/}.

\bibitem[AN06]{MR2227002}
V.~Alexeev and V.~V. Nikulin.
\newblock {\em Del {P}ezzo and {$K3$} surfaces}, volume~15 of {\em MSJ
  Memoirs}.
\newblock Mathematical Society of Japan, Tokyo, 2006.

\bibitem[Bau13]{bauer_diplom}
S.~Bauer.
\newblock Die {M}anin-{V}ermutung f\"ur eine del-{P}ezzo-{F}l\"ache, Master
  thesis, Universit\"at M\"unchen, 2013.

\bibitem[BB07]{MR2320172}
R.~de~la Bret{\`e}che and T.~D. Browning.
\newblock On {M}anin's conjecture for singular del {P}ezzo surfaces of degree
  4. {I}.
\newblock {\em Michigan Math. J.}, 55(1):51--80, 2007.

\bibitem[BB11]{MR2838351}
R.~de~la Bret{\`e}che and T.~D. Browning.
\newblock Manin's conjecture for quartic del {P}ezzo surfaces with a conic
  fibration.
\newblock {\em Duke Math. J.}, 160(1):1--69, 2011.

\bibitem[BB13]{arXiv:1011.3434}
S.~Baier and T.~D. Browning.
\newblock Inhomogeneous cubic congruences and rational points on del {P}ezzo
  surfaces.
\newblock {\em J. reine angew. Math.}, 680:69--151, 2013.

\bibitem[BBD84]{MR767407}
D.~Bindschadler, L.~Brenton, and D.~Drucker.
\newblock Rational mappings of del {P}ezzo surfaces, and singular
  compactifications of two-dimensional affine varieties.
\newblock {\em Tohoku Math. J. (2)}, 36(4):591--609, 1984.

\bibitem[BBD07]{MR2332351}
R.~de~la Bret{\`e}che, T.~D. Browning, and U.~Derenthal.
\newblock On {M}anin's conjecture for a certain singular cubic surface.
\newblock {\em Ann. Sci. \'Ecole Norm. Sup. (4)}, 40(1):1--50, 2007.

\bibitem[BBP12]{MR2874644}
R.~de~la Bret{\`e}che, T.~D. Browning, and E.~Peyre.
\newblock On {M}anin's conjecture for a family of {C}h\^atelet surfaces.
\newblock {\em Ann. of Math. (2)}, 175(1):297--343, 2012.

\bibitem[BD09a]{MR2520769}
T.~D. Browning and U.~Derenthal.
\newblock Manin's conjecture for a cubic surface with {$D_5$} singularity.
\newblock {\em Int. Math. Res. Not. IMRN}, (14):2620--2647, 2009.

\bibitem[BD09b]{MR2543667}
T.~D. Browning and U.~Derenthal.
\newblock Manin's conjecture for a quartic del {P}ezzo surface with {$A_4$}
  singularity.
\newblock {\em Ann. Inst. Fourier (Grenoble)}, 59(3):1231--1265, 2009.

\bibitem[BD12]{arXiv:1205.0373}
S.~Baier and U.~Derenthal.
\newblock {Quadratic congruences on average and rational points on cubic
  surfaces}, arXiv:1205.0373, 2012.

\bibitem[BF04]{MR2099200}
R.~de~la Bret{\`e}che and {\'E}.~Fouvry.
\newblock L'\'eclat\'e du plan projectif en quatre points dont deux
  conjugu\'es.
\newblock {\em J. reine angew. Math.}, 576:63--122, 2004.

\bibitem[BM90]{MR1032922}
V.~V. Batyrev and Yu.~I. Manin.
\newblock Sur le nombre des points rationnels de hauteur born\'e des
  vari\'et\'es alg\'ebriques.
\newblock {\em Math. Ann.}, 286(1-3):27--43, 1990.

\bibitem[Bou09]{MR2573192}
D.~Bourqui.
\newblock Comptage de courbes sur le plan projectif \'eclat\'e en trois points
  align\'es.
\newblock {\em Ann. Inst. Fourier (Grenoble)}, 59(5):1847--1895, 2009.

\bibitem[Bou11]{MR2809202}
D.~Bourqui.
\newblock Fonction z\^eta des hauteurs des vari\'et\'es toriques non
  d\'eploy\'ees.
\newblock {\em Mem. Amer. Math. Soc.}, 211(994):viii+151, 2011.

\bibitem[Bou12]{arXiv:1205.3573}
D.~Bourqui.
\newblock {Exemples de comptage de courbes sur les surfaces}.
\newblock {\em Math. Ann., to appear}, arXiv:1205.3573, 2012.

\bibitem[BP04]{MR2029863}
V.~V. Batyrev and O.~N. Popov.
\newblock The {C}ox ring of a del {P}ezzo surface.
\newblock In {\em Arithmetic of higher-dimensional algebraic varieties (Palo
  Alto, CA, 2002)}, volume 226 of {\em Progr. Math.}, pages 85--103.
  Birkh\"auser Boston, Boston, MA, 2004.

\bibitem[Bre02]{MR1909606}
R.~de~la Bret{\`e}che.
\newblock Nombre de points de hauteur born\'ee sur les surfaces de del {P}ezzo
  de degr\'e 5.
\newblock {\em Duke Math. J.}, 113(3):421--464, 2002.

\bibitem[Bro06]{MR2250046}
T.~D. Browning.
\newblock The density of rational points on a certain singular cubic surface.
\newblock {\em J. Number Theory}, 119(2):242--283, 2006.

\bibitem[Bro07]{MR2362193}
T.~D. Browning.
\newblock An overview of {M}anin's conjecture for del {P}ezzo surfaces.
\newblock In {\em Analytic number theory}, volume~7 of {\em Clay Math. Proc.},
  pages 39--55. Amer. Math. Soc., Providence, RI, 2007.

\bibitem[Bro09]{MR2559866}
T.~D. Browning.
\newblock {\em Quantitative arithmetic of projective varieties}, volume 277 of
  {\em Progress in Mathematics}.
\newblock Birkh\"auser Verlag, Basel, 2009.

\bibitem[BT98]{MR1620682}
V.~V. Batyrev and Yu. Tschinkel.
\newblock Manin's conjecture for toric varieties.
\newblock {\em J. Algebraic Geom.}, 7(1):15--53, 1998.

\bibitem[BW79]{MR80f:14021}
J.~W. Bruce and C.~T.~C. Wall.
\newblock On the classification of cubic surfaces.
\newblock {\em J. London Math. Soc. (2)}, 19(2):245--256, 1979.

\bibitem[Cay69]{cayley}
A.~Cayley.
\newblock A memoir on cubic surfaces.
\newblock {\em Philos. Trans. Roy. Soc. London}, 159:231--326, 1869.

\bibitem[CLT02]{MR1906155}
A.~Chambert-Loir and Yu. Tschinkel.
\newblock On the distribution of points of bounded height on equivariant
  compactifications of vector groups.
\newblock {\em Invent. Math.}, 148(2):421--452, 2002.

\bibitem[Cox95]{MR95i:14046}
D.~A. Cox.
\newblock The homogeneous coordinate ring of a toric variety.
\newblock {\em J. Algebraic Geom.}, 4(1):17--50, 1995.

\bibitem[CT88]{MR89f:11083}
D.~F. Coray and M.~A. Tsfasman.
\newblock Arithmetic on singular {D}el {P}ezzo surfaces.
\newblock {\em Proc. London Math. Soc. (3)}, 57(1):25--87, 1988.

\bibitem[CTS76]{MR54:2657}
J.-L. Colliot-Th{\'e}l{\`e}ne and J.-J. Sansuc.
\newblock Torseurs sous des groupes de type multiplicatif; applications \`a
  l'\'etude des points rationnels de certaines vari\'et\'es alg\'ebriques.
\newblock {\em C. R. Acad. Sci. Paris S\'er. A-B}, 282(18):Aii, A1113--A1116,
  1976.

\bibitem[CTS77a]{MR0447250}
J.-L. Colliot-Th{\'e}l{\`e}ne and J.-J. Sansuc.
\newblock La descente sur une vari\'et\'e rationnelle d\'efinie sur un corps de
  nombres.
\newblock {\em C. R. Acad. Sci. Paris S\'er. A-B}, 284(19):A1215--A1218, 1977.

\bibitem[CTS77b]{MR0447246}
J.-L. Colliot-Th{\'e}l{\`e}ne and J.-J. Sansuc.
\newblock Vari\'et\'es de premi\`ere descente attach\'ees aux vari\'et\'es
  rationnelles.
\newblock {\em C. R. Acad. Sci. Paris S\'er. A-B}, 284(16):A967--A970, 1977.

\bibitem[CTS80]{MR605344}
J.-L. Colliot-Th{\'e}l{\`e}ne and J.-J. Sansuc.
\newblock La descente sur les vari\'et\'es rationnelles.
\newblock In {\em Journ\'ees de G\'eometrie Alg\'ebrique d'Angers, Juillet
  1979/Algebraic Geometry, Angers, 1979}, pages 223--237. Sijthoff \&
  Noordhoff, Alphen aan den Rijn, 1980.

\bibitem[CTS87]{MR89f:11082}
J.-L. Colliot-Th{\'e}l{\`e}ne and J.-J. Sansuc.
\newblock La descente sur les vari\'et\'es rationnelles. {II}.
\newblock {\em Duke Math. J.}, 54(2):375--492, 1987.

\bibitem[Der07]{arXiv:0710.1583}
U.~Derenthal.
\newblock {Manin's conjecture for a quintic del Pezzo surface with
  $\mathbf{A}_2$ singularity}, arXiv:0710.1583, 2007.

\bibitem[Der09]{MR2520770}
U.~Derenthal.
\newblock Counting integral points on universal torsors.
\newblock {\em Int. Math. Res. Not. IMRN}, (14):2648--2699, 2009.

\bibitem[DF13a]{arXiv:1302.6151}
U.~Derenthal and C.~Frei.
\newblock {Counting imaginary quadratic points via universal torsors},
  arXiv:1302.6151, 2013.

\bibitem[DF13b]{arXiv:1304.3352}
U.~Derenthal and C.~Frei.
\newblock {Counting imaginary quadratic points via universal torsors, II}.
\newblock {\em Math. Proc. Cambridge Philos. Soc., to appear}, arXiv:1304.3352,
  2013.

\bibitem[DJT08]{MR2377367}
U.~Derenthal, M.~Joyce, and Z.~Teitler.
\newblock The nef cone volume of generalized del {P}ezzo surfaces.
\newblock {\em Algebra Number Theory}, 2(2):157--182, 2008.

\bibitem[DL10]{MR2753646}
U.~Derenthal and D.~Loughran.
\newblock Singular del {P}ezzo surfaces that are equivariant compactifications.
\newblock {\em Zap. Nauchn. Sem. S.-Peterburg. Otdel. Mat. Inst. Steklov.
  (POMI)}, 377(Issledovaniya po Teorii Chisel. 10):26--43, 241, 2010.

\bibitem[Dol12]{MR2964027}
I.~V. Dolgachev.
\newblock {\em Classical algebraic geometry}.
\newblock Cambridge University Press, Cambridge, 2012.
\newblock A modern view.

\bibitem[DP80]{MR579026}
M.~Demazure and H.~C. Pinkham, editors.
\newblock {\em S\'eminaire sur les {S}ingularit\'es des {S}urfaces}, volume 777
  of {\em Lecture Notes in Mathematics}.
\newblock Springer, Berlin, 1980.
\newblock Held at the Centre de Math\'ematiques de l'\'Ecole Polytechnique,
  Palaiseau, 1976--1977.

\bibitem[DT07]{MR2290499}
U.~Derenthal and Yu. Tschinkel.
\newblock Universal torsors over del {P}ezzo surfaces and rational points.
\newblock In {\em Equidistribution in number theory, an introduction}, volume
  237 of {\em NATO Sci. Ser. II Math. Phys. Chem.}, pages 169--196. Springer,
  Dordrecht, 2007.

\bibitem[DV34]{duval}
P.~Du~Val.
\newblock {On isolated singularities of surfaces which do not affect the
  conditions of adjunction (Parts II and III)}.
\newblock {\em Proc. Cambridge Philos. Soc.}, 30:460--465, 483--491, 1934.

\bibitem[FMT89]{MR89m:11060}
J.~Franke, Yu.~I. Manin, and Yu. Tschinkel.
\newblock Rational points of bounded height on {F}ano varieties.
\newblock {\em Invent. Math.}, 95(2):421--435, 1989.

\bibitem[Fre12]{arXiv:1204.0383}
C.~Frei.
\newblock {Counting rational points over number fields on a singular cubic
  surface}.
\newblock {\em Algebra Number Theory, to appear}, arXiv:1204.0383, 2012.

\bibitem[GM05]{MR2115006}
C.~Galindo and F.~Monserrat.
\newblock The total coordinate ring of a smooth projective surface.
\newblock {\em J. Algebra}, 284(1):91--101, 2005.

\bibitem[Has09]{MR2498061}
B.~Hassett.
\newblock Rational surfaces over nonclosed fields.
\newblock In {\em Arithmetic geometry}, volume~8 of {\em Clay Math. Proc.},
  pages 155--209. Amer. Math. Soc., Providence, RI, 2009.

\bibitem[Hau13]{arXiv:1106.0854}
J.~Hausen.
\newblock Three lectures on {C}ox rings.
\newblock In {\em Torsors, {\'E}tale Homotopy and {A}pplications to {R}ational
  {P}oints}, volume 405 of {\em London Math. Soc. Lecture Note Ser.}, pages
  3--60. Cambridge Univ. Press, Cambridge, 2013.

\bibitem[HB03]{MR2075628}
D.~R. Heath-Brown.
\newblock The density of rational points on {C}ayley's cubic surface.
\newblock In {\em Proceedings of the Session in Analytic Number Theory and
  Diophantine Equations}, volume 360 of {\em Bonner Math. Schriften}, page~33,
  Bonn, 2003. Univ. Bonn.

\bibitem[HK00]{MR2001i:14059}
Y.~Hu and S.~Keel.
\newblock Mori dream spaces and {GIT}.
\newblock {\em Michigan Math. J.}, 48:331--348, 2000.
\newblock Dedicated to William Fulton on the occasion of his 60th birthday.

\bibitem[HP52]{MR13:972c}
W.~V.~D. Hodge and D.~Pedoe.
\newblock {\em Methods of algebraic geometry. {V}ol. {II}. {B}ook {III}:
  {G}eneral theory of algebraic varieties in projective space. {B}ook {IV}:
  {Q}uadrics and {G}rassmann varieties}.
\newblock Cambridge, at the University Press, 1952.

\bibitem[HS10]{MR2671185}
J.~Hausen and H.~S{\"u}{\ss}.
\newblock The {C}ox ring of an algebraic variety with torus action.
\newblock {\em Adv. Math.}, 225(2):977--1012, 2010.

\bibitem[HT04]{MR2029868}
B.~Hassett and Yu. Tschinkel.
\newblock Universal torsors and {C}ox rings.
\newblock In {\em Arithmetic of higher-dimensional algebraic varieties (Palo
  Alto, CA, 2002)}, volume 226 of {\em Progr. Math.}, pages 149--173.
  Birkh\"auser Boston, Boston, MA, 2004.

\bibitem[HW81]{MR646042}
F.~Hidaka and K.~Watanabe.
\newblock Normal {G}orenstein surfaces with ample anti-canonical divisor.
\newblock {\em Tokyo J. Math.}, 4(2):319--330, 1981.

\bibitem[LB12a]{MR2990624}
P.~Le~Boudec.
\newblock Manin's conjecture for a cubic surface with {$2\bold A_2+\bold A_1$}
  singularity type.
\newblock {\em Math. Proc. Cambridge Philos. Soc.}, 153(3):419--455, 2012.

\bibitem[LB12b]{MR2961294}
P.~Le~Boudec.
\newblock Manin's conjecture for a quartic del {P}ezzo surface with {$\bold
  A_3$} singularity and four lines.
\newblock {\em Monatsh. Math.}, 167(3-4):481--502, 2012.

\bibitem[LB12c]{MR2853047}
P.~Le~Boudec.
\newblock Manin's conjecture for two quartic del {P}ezzo surfaces with {$3{\bf
  A}_1$} and {${\bf A}_1+{\bf A}_2$} singularity types.
\newblock {\em Acta Arith.}, 151(2):109--163, 2012.

\bibitem[LB12d]{arXiv:1207.2685}
P.~Le~Boudec.
\newblock Affine congruences and rational points on a certain cubic surface,
  arXiv:1207.2685, 2012.

\bibitem[Lou10]{MR2769338}
D.~Loughran.
\newblock Manin's conjecture for a singular sextic del {P}ezzo surface.
\newblock {\em J. Th\'eor. Nombres Bordeaux}, 22(3):675--701, 2010.

\bibitem[Lou12]{MR2980925}
D.~Loughran.
\newblock Manin's conjecture for a singular quartic del {P}ezzo surface.
\newblock {\em J. Lond. Math. Soc. (2)}, 86(2):558--584, 2012.

\bibitem[LV09]{MR2529093}
A.~Laface and M.~Velasco.
\newblock Picard-graded {B}etti numbers and the defining ideals of {C}ox rings.
\newblock {\em J. Algebra}, 322(2):353--372, 2009.

\bibitem[Mol12]{moll_diplom}
S.~Moll.
\newblock {Rationale Punkte auf einer singul\"aren del Pezzo Fl\"ache},
  Diplomarbeit, Universit\"at M\"unchen, 2012.

\bibitem[Pey98]{MR1679842}
E.~Peyre.
\newblock Terme principal de la fonction z\^eta des hauteurs et torseurs
  universels.
\newblock {\em Ast\'erisque}, (251):259--298, 1998.
\newblock Nombre et r\'epartition de points de hauteur born\'ee (Paris, 1996).

\bibitem[Sal98]{MR1679841}
P.~Salberger.
\newblock Tamagawa measures on universal torsors and points of bounded height
  on {F}ano varieties.
\newblock {\em Ast\'erisque}, (251):91--258, 1998.
\newblock Nombre et r\'epartition de points de hauteur born\'ee (Paris, 1996).

\bibitem[Sch63]{schlaefli}
L.~Schl\"afli.
\newblock On the distribution of surfaces of the thrid order into species, in
  reference to the absence or presence of singular points, and the relaity of
  their lines.
\newblock {\em Philos. Trans. Roy. Soc. London}, 153:193--241, 1863.

\bibitem[Sko93]{MR1260765}
A.~N. Skorobogatov.
\newblock On a theorem of {E}nriques-{S}winnerton-{D}yer.
\newblock {\em Ann. Fac. Sci. Toulouse Math. (6)}, 2(3):429--440, 1993.

\bibitem[STV07]{MR2358614}
M.~Stillman, D.~Testa, and M.~Velasco.
\newblock Gr\"obner bases, monomial group actions, and the {C}ox rings of del
  {P}ezzo surfaces.
\newblock {\em J. Algebra}, 316(2):777--801, 2007.

\bibitem[TVAV09]{MR2579393}
D.~Testa, A.~V{\'a}rilly-Alvarado, and M.~Velasco.
\newblock Cox rings of degree one del {P}ezzo surfaces.
\newblock {\em Algebra Number Theory}, 3(7):729--761, 2009.

\bibitem[TVAV11]{MR2824848}
D.~Testa, A.~V{\'a}rilly-Alvarado, and M.~Velasco.
\newblock Big rational surfaces.
\newblock {\em Math. Ann.}, 351(1):95--107, 2011.

\bibitem[Ura83]{MR713283}
T.~Urabe.
\newblock On singularities on degenerate del {P}ezzo surfaces of degree {$1,$}
  {$2$}.
\newblock In {\em Singularities, {P}art 2 ({A}rcata, {C}alif., 1981)},
  volume~40 of {\em Proc. Sympos. Pure Math.}, pages 587--591. Amer. Math.
  Soc., Providence, R.I., 1983.

\bibitem[Ye02]{MR1933881}
Q.~Ye.
\newblock On {G}orenstein log del {P}ezzo surfaces.
\newblock {\em Japan. J. Math. (N.S.)}, 28(1):87--136, 2002.

\end{thebibliography}

\end{document}